\documentclass[10pt]{article}
\usepackage{amsthm,amsfonts,amsmath,amscd,amssymb,epsfig,enumerate,array}
\theoremstyle{plain}
\newtheorem{theorem}{Theorem}[section]
\newtheorem{thm}[theorem]{Theorem}

\newtheorem{cor}[theorem]{Corollary}

\newtheorem{prop}[theorem]{Proposition}
\newtheorem{lemma}[theorem]{Lemma}
\newtheorem{definition}[theorem]{Definition}
\newtheorem*{theorem*}{Theorem}

\theoremstyle{remark}

\newtheorem{remark}[theorem]{Remark}

\newcommand{\id}{{{\mathchoice {\rm 1\mskip-4mu l} {\rm 1\mskip-4mu l}
{\rm 1\mskip-4.5mu l} {\rm 1\mskip-5mu l}}}}

\newcommand{\p}{\partial}
\newcommand{\pb}{\bar\partial}
\newcommand{\om}{\omega}
\newcommand{\Om}{\Omega}
\newcommand{\eps}{\varepsilon}
\newcommand{\into}{\hookrightarrow}
\newcommand{\la}{\langle}
\newcommand{\ra}{\rangle}
\newcommand{\N}{{\mathbb{N}}}
\newcommand{\Z}{{\mathbb{Z}}}
\newcommand{\R}{{\mathbb{R}}}
\newcommand{\C}{{\mathbb{C}}}
\newcommand{\Q}{{\mathbb{Q}}}

\newcommand{\J}{{\bf J}}
\newcommand{\z}{{\bf z}}
\newcommand{\f}{{\bf f}}
\newcommand{\w}{{\bf w}}
\newcommand{\I}{{\bf I}}
\newcommand{\GW}{{\rm GW}}
\newcommand{\Aut}{{\rm Aut}}
\newcommand{\ind}{{\rm ind}}
\newcommand{\im}{{\rm im }}        
\newcommand{\st}{{\rm st}}
\newcommand{\ds}{{\rm ds}}

\renewcommand{\min}{{\rm min}}
\renewcommand{\max}{{\rm max}}

\newcommand{\reg}{{\rm reg}}
\newcommand{\codim}{{\rm codim}}

\newcommand{\ev}{{\rm ev}}

\newcommand{\End}{{\rm End}}

\newcommand{\Hom}{{\rm Hom}}

\newcommand{\Coh}{{\rm Coh}}

\newcommand{\nc}{{\rm nc}}
\newcommand{\EE}{\mathcal{E}}
\newcommand{\BB}{\mathcal{B}}
\newcommand{\JJ}{\mathcal{J}}
\newcommand{\MM}{\mathcal{M}}

\renewcommand{\AA}{\mathcal{A}}

\newcommand{\ZZ}{\mathcal{Z}}

\newcommand{\II}{\mathcal{I}}
\newcommand{\KK}{\mathcal{K}}
\title{Symplectic hypersurfaces and transversality in Gromov-Witten
  theory}
\author{K.~Cieliebak and K.~Mohnke}
\date{April 2008}
\begin{document}
\maketitle
\begin{abstract}
We present a new method to prove transversality for holomorphic curves in
symplectic manifolds, and show how it leads to a definition of genus zero
Gromov-Witten invariants. The main idea is to introduce additional marked
points that are mapped to a symplectic hypersurface of high degree in
order to stabilize the domains of holomorphic maps.
\end{abstract}
\tableofcontents

\section{Introduction}\label{sec:intro}

Donaldson's construction of symplectic hypersurfaces~\cite{Do} and its
extensions lead to a new method to prove transversality for
holomorphic curves in symplectic manifolds. In this paper, we
illustrate this method for holomorphic spheres in closed
symplectic manifolds, and show how it leads to a definition of genus
zero Gromov-Witten invariants.

The transversality question in Gromov-Witten theory has previously
been discussed by many authors, see e.g.~\cite{FO,LiT,LiuT,RT,Si,MS}
and further references therein. Transversality for $J$-holomorphic
disks with Lagrangian boundary conditions is discussed
in~\cite{FOOO}. A very general approach to transversality is being
developed by Hofer, Wysocki and Zehnder~\cite{HWZ-trans}. Their
``polyfold Fredholm theory'' provides a framework to
study transversality for solution spaces of elliptic PDEs, which
may have applications far beyond symplectic geometry. By contrast
our approach, being specifically designed for holomorphic curves,
appears to be much simpler in this case. In particular, we use
only standard tools from functional analysis and our moduli spaces
(in the genus zero case) are smooth finite-dimensional {\em
manifolds}.
\medskip

We now state our results. Consider a closed symplectic
manifold $(X^{2n},\om)$ of dimension $2n$ such that $\om$
represents an integral cohomology class $[\om]\in H^2(X;\Z)$. 
Let $J$ be an almost complex structure tamed by $\om$, i.e.~$\om(v,Jv)>0$ for
all tangent vectors $v$. Our construction
uses as auxiliary datum a {\em hypersurface} (i.e.~a closed
codimension two submanifold) $Y\subset X$, representing the homology
class Poincar\'e dual to $D[\om]$ for some {\em degree} $D\in\N$,
whose K\"ahler angle (see Section~\ref{sec:hyp}) with respect to
$(\om,J)$ satisfies $\theta(Y;\om,J)<\theta_3$.  
Here $0<\theta_3<\theta_1<1$ are universal constants depending only on
the dimension. In particular, the hypersurface $Y$ is symplectic. 
By a result of Donaldson~\cite{Do}, for a given $J$ such a hypersurface $Y$
exists for each sufficiently large $D\in\N$. We will refer to such
$(J,Y)$ as a {\em Donaldson pair} if in addition $D\geq D^*(X,\om,J)$
for a suitable positive constant $D^*$ depending on $(X,\om,J)$ (see
Definition~\ref{def:D-pair}). 

Given a Donaldson pair $(J,Y)$, we denote by $\JJ(X,Y;J,\theta_1)$ the
space of $\om$-tamed almost complex structures on $X$ which leave $TY$
invariant and are $\theta_1$-close to $J$ in the $C^0$-norm. For
$\ell\geq 3$ we define a class
$$
   \JJ_{\ell+1}(X,Y;J,\theta_1) \subset
   C^\infty\Bigl(\bar\MM_{\ell+1},\JJ(X,Y;J,\theta_1)\Bigr) 
$$
of almost complex structures depending smoothly on points in the
Deligne-Mumford space $\bar\MM_{\ell+1}$ of stable genus zero
curves with $\ell+1$ marked points $z_0,\dots,z_\ell$. Moreover, we
require elements in $\JJ_{\ell+1}(X,Y;J,\theta_1)$ to be {\em
  coherent} under gluing of strata in $\bar\MM_{\ell+1}$, see
Section~\ref{sec:coh} for the precise definition.

Via the canonical projection
$\pi:\bar\MM_{\ell+1}\to\bar\MM_\ell$ forgetting the marked point
$z_0$, we think of $K\in\JJ_{\ell+1}(X,Y;J,\theta_1)$
as a collection
of almost complex structures $K_\z\in\JJ(X,Y;J,\theta_1)$, parametrized by
$\z\in\bar\MM_\ell$,
depending smoothly on $z_0$ thought of as a point on
the nodal Riemann surface $\pi^{-1}(\z)$ defined by $\z$. To
$J\in\JJ_{\ell+1}(X,Y;J,\theta_1)$ and a stable map $(\z,\f)$ with $\ell$
marked points we associate the $(0,1)$-form
$\pb_{K_{\st(\z)}}\f$, where $\st(\z)\in\bar\MM_\ell$ is the
stabilization of the underlying nodal curve $\z$. We call
$(\z,\f)$ $K$-holomorphic if $\pb_{K_{\st(\z)}}\f=0$.

Denote by $c_1$ the first Chern class of $X$. For a homology class
$A\in H_2(X;\Z)$ set $\ell:=D\om(A)$. In addition, fix an integer
$k\geq 0$. Every $K\in\JJ_{\ell+1}(X,Y;J,\theta_1)$ induces
elements $\pi_\ell^*K\in\JJ_{\ell+k+1}(X,Y;J,\theta_1)$ by
composition with the map
$\pi_\ell:\bar\MM_{k+\ell+1}\to\bar\MM_{\ell+1}$ that forgets the
first $k$ marked points (and stabilizes). Using this, we define the
moduli space $\MM_{k+\ell}(A,K;Y)$ of (smooth)
$\pi_\ell^*K$-holomorphic spheres in the class $A$ with $k+\ell$
marked points and mapping the last $\ell$ marked points to $Y$.

More generally, for
a $(k+\ell)$-labelled tree $T$ (see Section~\ref{sec:curves}) denote
by $\MM_T(\!A,\!K\!;\!Y)$ the space of stable $K$-holomorphic spheres,
modelled over $T$ and mapping the last $\ell$ marked points to $Y$. We
will assume that every {\em ghost tree}, i.e.~a maximal subtree
$T'\subset T$ corresponding to components of the nodal surface on
which the map $\f$ is constant, contains at most one of the last
$\ell$ marked points. Our first result is

\begin{thm}\label{thm:moduli}
Let $(X,\om)$ be a closed symplectic manifold such that $[\om]\in
H^2(X;\Z)$. Then for every Donaldson pair $(J,Y)$ of
degree $D\geq D^*(X,\om,J)$ there exist nonempty 
sets
$\JJ^\reg_{\ell+1}(X,Y;J,\theta_1)\subset\JJ_{\ell+1}(X,Y;J,\theta_1)$,
$\ell\geq 3$, such that for $K\in\JJ_{\ell+1}^\reg(X,Y;J,\theta_1)$ the
following holds: Let 
$A\in H_2(X;\Z)$ be a homology class with $\ell=D\om(A)$. For $k\geq
0$ let $T$ be a $(k+\ell)$-labelled tree with $e(T)$ edges such that
every ghost tree contains at most one of the last $\ell$ marked points. 
Then the moduli space $\MM_T(A,K;Y)$ of stable $K$-holomorphic spheres
in the class $A$, modelled over $T$ and mapping the last $\ell$ marked
points to $Y$, is a smooth manifold of real dimension
$$
   \dim_\R\MM_T(A,K;Y) = 2\Bigl(n-3+k+c_1(A)-e(T)\Bigr).
$$
\end{thm}

To interpret this result, recall from~\cite{MS} the definition of a
pseudocycle. We say that a subset $\Om\subset X$ of a manifold $X$ has
{\em dimension at most $d$} if $\Om$ is contained in the image of a
smooth map $g:N\to X$ from a manifold of dimension $\dim N\leq d$.
A smooth map $f:M\to X$ from an oriented $d$-manifold $M$
to a manifold $X$ is called {\em $d$-dimensional pseudocycle} if
$f(M)$ has compact closure and its "omega limit set"
$$
   \Om_f := \bigcap_{K\subset M \text{ compact}}\overline{f(M\setminus
   K)}
$$
has dimension at most $d-2$. Thus intuitively, $f(M)$ can be
compactified by adding strata of codimension at least $2$.

Our second result, proved in Section~\ref{sec:proof}, will be a
consequence of Theorem~\ref{thm:moduli} and an appropriate version of
Gromov compactness. 

\begin{thm}\label{thm:pseudo}
In the notation of Theorem~\ref{thm:moduli}, for every $k\geq 1$ the
evaluation map at the first $k$ points
$$
   \ev^k:\MM_{k+\ell}(A,K;Y)\mapsto X^k
$$
represents a pseudocycle $\ev^k(J,Y;K)$ of dimension
$$
   2d := 2\Bigl(n-3+k+c_1(A)\Bigr).
$$
\end{thm}

The third result concerns the dependence of this pseudocycle on
the auxiliary data $(J,Y;K)$. Two pseudocycles $f_i:M_i\to X$,
$i=0,1$ are called {\em cobordant} if there exists a smooth map
$F:W\to X$ from an oriented $(d+1)$-manifold $W$ with boundary $\p
W=M_1-M_0$ such that $F|_{M_i}=f_i$ and $\dim\Om_F\leq d-1$. Let a
{\em rational pseudocycle} be a pseudocycle multiplied with a
positive rational weight.
For pseudocycles $f_0$, $f_1$ we say that $f_0$ and $\ell f_1$
are {\em equal as currents} if there exists
a covering map $\phi:M_0\to M_1$ of degree $\ell$ such that
$f_0=f_1\circ\phi$.
We denote by {\em rational cobordism} the
equivalence relation on rational pseudocycles generated by
equality as currents and cobordism of pseudocycles.

The following result, proved in Section~\ref{sec:indep}, will be a
consequence of the asymptotic uniqueness~\cite{Aur} of
Donaldson hypersurfaces.

\begin{thm}\label{thm:indep}
Up to rational cobordism, the rational pseudocycle
$\frac{1}{\ell!}\ev^k(J,Y;K)$ in Theorem~\ref{thm:pseudo} does not
depend on the choice of the auxiliary data $(J,Y;K)$.
\end{thm}

Following the strategy of~\cite{MS}, Theorems~\ref{thm:pseudo}
and~\ref{thm:indep} allow us to
define Gromov-Witten invariants for arbitrary closed symplectic
manifolds $(X,\om)$ as follows. Let $A\in H_2(X;\Z)$ and let
$\alpha_1,\dots,\alpha_k$ be non-torsion integral cohomology classes
on $X$ of total degree
$$
   \sum_{i=1}^k\deg(\alpha_i) = 2d = 2n-6+2k+2c_1(A).
$$
Represent the Poincar\'e dual of
$\pi_1^*\alpha_1\cup\dots\cup\pi_k^*\alpha_k\in H^*(X^k)$ by a
cycle $a$ in $X^k$ that is {\em strongly transverse} in the sense
of~\cite{MS} to the evaluation map $\ev^k(J,Y;K):
\MM_{k+\ell}(A,J;Y)\mapsto X^k$ from Theorem~\ref{thm:pseudo}. Define
the {\em Gromov-Witten invariant} of holomorphic spheres passing
through cycles dual to $\alpha_1\!,\!\dots\!,\!\alpha_k$ as the intersection
number
$$
   \GW_A(\alpha_1,\dots,\alpha_k) := \frac{1}{\ell!}\ev^k(J,Y;K)\cdot
   a \in \Q,
$$
with $\ell=D\om(A)$ as in Theorem~\ref{thm:moduli}. In view of
Theorem~\ref{thm:indep}, this number does not depend on the choice
of $(J,Y;K)$ and by Theorem~\ref{thm:pseudo} not on the choice of $a$.

\begin{remark}
The method described in this paper is applicable to prove
transversality for more general holomorphic curves. 
The extension to genus zero holomorphic curves with boundary on a
Lagrangian submanifold as in Floer homology~\cite{FOOO} relies on the
existence of approximately holomorphic hypersurfaces in the complement
of a Lagrangian submanifold due to Auroux, Gayet and Mohsen~\cite{AGM},
and the extension to genus zero holomorphic curves with punctures
asymptotic to closed Reeb orbits as in symplectic field
theory~\cite{EGH} is based on the existence of contact open books due
to Giroux~\cite{Gi}. Both cases also require good descriptions of the
relevant Deligne-Mumford spaces of stable genus zero curves (in the
spirit of Section~\ref{sec:coh}): curves with marked points and
boundary in the first case, and curves with marked points and
directions at the marked points in the second case. Note that in both
cases the moduli spaces have codimension one boundary, so one needs to
prove appropriate gluing results. 
These extensions may be treated in subsequent papers.
\end{remark}

\begin{remark}
This paper is clearly related to the work of Ionel and Parker on
relative Gromov-Witten invariants~\cite{IP1,IP2}. Indeed, moduli
spaces of holomorphic curves with tangency conditions to one or more
symplectic hypersurfaces are the central objects of study in their as
well as our work. However, the motivations and goals are very
different. In particular, we do {\em not} define relative
Gromov-Witten invariants (although this should be possible 
using the techniques developed in this paper). 
\end{remark}

\begin{remark}
One may wonder how the Gromov-Witten invariants defined in this paper
are related to other definitions. First, we claim that if the
symplectic manifold $(X,\om)$ is semipositive in the sense
of McDuff and Salamon~\cite{MS}, then our definition agrees with the
one given in~\cite{MS}. This essentially follows from the fact that in
that case we can choose our perturbation $K$ {\em domain-independent}.   
On the other hand, the comparison of our definition with ones
involving abstract perturbations and virtual moduli cycles as
in~\cite{FO,LiT,LiuT,RT,Si} lies beyong the scope of this paper, whose
purpose is exactly to avoid abstract perturbations. However, taking
virtual moduli cycles for granted, our invariant seems to coincide with
the relative Gromov-Witten invariant in~\cite{IP1}, which in turn
equals the absolute Gromov-Witten invariant according to
Proposition 14.9 in~\cite{IP2}. 
\end{remark}

{\em Acknowledgements. }
The ideas for this paper emerged in discussions with Y.~Eliashberg
back in 2000. Over the subsequent years we benefited from discussions
with many people including D.~Auroux, D.~McDuff, O.~Fabert, E.~Ionel,
J.~Latschev, D.~Salamon, and C.~Stamm. 
Finally, we thank the referee for pointing out several gaps and errors
in earlier versions.

\section{Nodal curves}\label{sec:curves}

In this section we describe the space $\bar\MM_k$ of stable curves
of genus zero with $k$ marked points. We adopt the approach
and notation of~\cite{MS}.

A {\em $k$-labelled tree} is a triple $T=(T,E,\Lambda)$, where $(T,E)$
is a (connected) tree with set of vertices $T$ and edge relation
$E\subset T\times
T$, and $\Lambda=\{\Lambda_\alpha\}_{\alpha\in T}$ is a decomposition
of the index set $\{1,\dots,k\}=\amalg_{\alpha\in T}\Lambda_\alpha$. We
write $\alpha E\beta$ if $(\alpha,\beta)\in E$. Note that the
labelling $\Lambda$ defines a unique map $\{1,\dots,k\}\to T$,
$i\mapsto\alpha_i$ by the requirement $i\in\Lambda_{\alpha_i}$.
Let
$$
   e(T) = |T|-1
$$
be the number of edges. A {\em tree homomorphism} $\tau:T\to\tilde T$
is a map which collapses some subtrees of $T$ to vertices of $\tilde
T$ (cf.~\cite{MS}). A tree homomorphism $\tau$ is called {\em tree
  isomorphism} if it is bijective and $\tau^{-1}$ is a tree
homomorphism.

A tree $T$ is called {\em stable} if for each $\alpha\in T$,
$$
   n_\alpha := \#\Lambda_\alpha + \#\{\beta\mid \alpha E\beta\} \geq
   3.
$$
Note that for $k<3$ every $k$-labelled tree is unstable. For $k\geq
3$, a $k$-labelled tree $T$ can be stabilized in a canonical way to a
stable $k$-labelled tree $\st(T)$ by deleting the vertices with
$n_\alpha<3$ and modifying the edges in the obvious way
(cf.~\cite{MS}).

A {\em nodal curve of genus zero with $k$ marked points
  modelled over the tree $T=(T,E,\Lambda)$} is a tuple
$$
   \z = (\{z_{\alpha\beta}\}_{\alpha E\beta},\{z_i\}_{1\leq i\leq k})
$$
of points $z_{\alpha\beta},z_i\in S^2$ such that for each $\alpha\in
T$ the {\em special points}
$$
   SP_\alpha:=\{z_{\alpha\beta}\mid \alpha E\beta\}\cup \{z_i\mid
   \alpha_i=\alpha\}
$$
are pairwise distinct. For $\alpha\in T$ and $i\in\{1,...,k\}$
denote by $z_{\alpha i}$ either the point $z_i$ if
$i\in\Lambda_\alpha$ or the point $z_{\alpha\beta_1}$ if
$z_i\in\Lambda_{\beta_r}$ and
$(\alpha,\beta_1),(\beta_1,\beta_2),...,(\beta_{r-1},\beta_r)\in
E$. Note
that $n_\alpha=\# SP_\alpha$. We associate to $\z$ the {\em nodal
Riemann surface}
$$
   \Sigma_\z := \coprod_{\alpha\in T}S_\alpha\Bigr/z_{\alpha\beta}\sim
   z_{\beta\alpha},
$$
obtained by gluing a collection of standard spheres
$\{S_\alpha\}_{\alpha\in
  T}$ at the points $z_{\alpha\beta}$ for $\alpha E\beta$, with marked
points $z_i\in S_{\alpha_i}$, $i=1,\dots,k$. Note that $\z$ can be
uniquely recovered from $\Sigma_\z$, so we will sometimes not
distinguish between the two. A nodal curve $\z$ is called {\em
stable} if the underlying tree is stable, i.e., every sphere
$S_\alpha$ carries at least $3$ special points. Stabilization of
trees induces a canonical stabilization of nodal curves
$\z\mapsto\st(\z)$.

We will usually omit the
genus zero from the notation. Denote the space of all nodal curves
(of genus zero) with $k$ marked points modelled over $T$ by
$$
   \tilde\MM_T\subset (S^2)^E\times (S^2)^k.
$$
Note that this is an open subset of the product of spheres.
A {\em
morphism} between nodal curves $\z,\tilde\z$ modelled over trees
$T,\tilde T$ is a tuple
$$
   \phi=(\tau,\{\phi_\alpha\}_{\alpha\in T}),
$$
where $\tau:T\to\tilde T$ is a tree homomorphism and
$\phi_\alpha:S^2\cong S_\alpha\to S_{\tau(\alpha)}\cong S^2$ are
(possibly constant) holomorphic maps such that
\begin{align*}
   \tilde
   z_{\tau(\alpha)\tau(\beta)} &= \phi_\alpha(z_{\alpha\beta}) \mbox{ if }
   \tau(\alpha)\neq\tau(\beta),\\
   \phi_\alpha(z_{\alpha\beta}) &= \phi_\beta(z_{\beta\alpha}) \mbox{ if }
   \tau(\alpha)= \tau(\beta),\\
   \tilde\alpha_i & =\tau(\alpha_i),\qquad \tilde
   z_i=\phi_{\alpha_i}(z_i)
\end{align*}
for $i=1,\dots,k$ and $\alpha E\beta$. A morphism
$\phi:\z\to\tilde\z$ induces a natural holomorphic map
$\Sigma_\z\to\Sigma_{\tilde\z}$ (i.e., a continuous map that is
holomorphic on each component $S_\alpha$). A morphism
$(\tau,\{\phi_\alpha\})$ is called {\em isomorphism} if $\tau$ is a
tree isomorphism and each $\phi_\alpha$ is biholomorphic.

Isomorphisms from $\z$ to itself are called {\em automorphisms}.
If $\z$ is stable its only automorphism is the identity
(see~\cite{MS}, discussion after Definition D.3.4). Thus for a stable
tree $T$ we have a free and proper holomorphic action
$$
   G_T\times\tilde\MM_T\to\tilde\MM_T
$$
of the group $G_T$ of isomorphisms fixing $T$.

\begin{remark}
Properness follows directly from the characterization of sequences
of M\"obius transforms which do not have a uniformly convergent
subsequence in~\cite{MS}, Lemma~D.1.2. Indeed, if $\z^\nu\to \z$ and
$\phi^\nu(\z^\nu)\to \z'$ in $\tilde\MM_T$ as $\nu\to\infty$ and
$\{\phi^\nu\}$ has no uniformly convergent subsequence in $G_T$, then
there is an $\alpha\in T$ and a subsequence $\phi^{\nu'}$ which converges
to the constant map $S_\alpha\to \{y\}$ on compact subsets of
$S_\alpha\setminus\{x\}$ for some points $x,y\in S_\alpha$. In
particular, the limit of the sets $\{\phi^{\nu'}_\alpha(z^{\nu'}_i),
\phi^{\nu'}_\alpha(z^{\nu'}_{\alpha\beta}) \mid z_i\in \Lambda_\alpha,
\alpha E\beta\}$ consists of at most two points, contradicting
$\phi^\nu(\z^\nu)\to \z'$ and the stability of $\z'$.
\end{remark}

Hence the quotient
$$
   \MM_T :=\tilde\MM_T/G_T
$$
is a complex manifold of dimension
$$
   \dim_\C\MM_T = k+2e(T)-3|T| = k-3-e(T).
$$
For $k\geq 3$, denote by $\MM_k=\tilde\MM_k/G$ the moduli space of
stable curves modelled over the $k$-labelled tree with one
vertex. As a set, the {\em Deligne-Mumford space (of
genus zero) with $k$ marked points} is given by
$$
   \bar\MM_k := \coprod_T\MM_T,
$$
where the union is taken over the (finitely many) isomorphism
classes of stable $k$-labelled trees. However, $\bar\MM_k$ is
equipped with the topology of Gromov convergence which makes it a
compact connected metrizable space (see~\cite{MS}). As the
notation suggests, $\bar\MM_k$ is the compactification of $\MM_k$
in the Gromov topology. For a stable $k$-labelled tree $T$, the
closure of $\MM_T$ in $\bar\MM_k$ is given by
$$
   \bar\MM_T = \coprod_{\tilde T}\MM_{\tilde T},
$$
where the union is taken over all isomorphism classes of stable
$k$-labelled trees $\tilde T$ for which there exists a surjective tree
homomorphism $\tau:\tilde T\to T$ with $\tau(\tilde\alpha_i)=\alpha_i$
for $i=1,\dots,k$.

We have the following result of Knudsen (cf.~\cite{MS}).

\begin{thm}\label{thm:DM}
For $k\geq 3$, the Deligne-Mumford space $\bar\MM_k$ is a compact
complex manifold of dimension $\dim_\C\bar\MM_k=k-3$. Moreover, for
each stable $k$-labelled tree $T$, the space
$\bar\MM_T\subset\bar\MM_k$ is a compact complex submanifold of
codimension $\codim_\C\bar\MM_T=e(T)$.
\end{thm}

\begin{remark}
In this paper we do not use the complex structure on $\bar\MM_k$, but
only its structure as a smooth manifold.
\end{remark}

We have a canonical projection $\pi:\bar\MM_{k+1}\to\bar\MM_k$ by
forgetting the $(k+1)$-st marked point and stabilizing. The map $\pi$
is holomorphic and the fibre $\pi^{-1}([\z])$ is naturally
biholomorphic to $\Sigma_\z$. The projection
$\pi:\bar\MM_{k+1}\to\bar\MM_k$ is called the {\em universal curve}.

\begin{lemma}\label{lem:fibre}
For $[\z]\in\bar\MM_k$, every component $S_\alpha$ of
$\Sigma_\z\cong\pi^{-1}([\z])$ is an embedded holomorphic sphere in
$\bar\MM_{k+1}$.
\end{lemma}

\begin{proof}
This follows from the construction of coordinate charts on
$\bar\MM_{k+1}$ in terms of
cross ratios in~\cite{MS} and the fact that any nonconstant cross
ratio maps $S_\alpha$ biholomorphically onto $\C P^1$.
\end{proof}

The following two lemmas follow directly from the construction of
the stabilization.

\begin{lemma}\label{lem:stab}
Let $\z$ be a nodal curve modelled over the tree $T$ and $\st(\z)$ its
stabilization. Then the stabilization map induces a morphism
$$
   \st=(\tau,\{\phi_\alpha\}_{\alpha\in
   T}):\Sigma_\z\mapsto\Sigma_{\st(\z)}
$$
with the following property: There exists a collection of subtrees
$T'\subset T$
such that $\phi_\alpha$ is constant for $\alpha\in T'$ and biholomorphic
otherwise. Moreover, each sphere $S_\alpha$ with $\alpha\notin T'$
carries at least $3$ special points.
\end{lemma}
\begin{lemma}\label{lem:proj}
The projection
$\pi_\ell:=\pi\circ\dots\circ\pi:\bar\MM_{k+\ell}\to\bar\MM_\ell$
induces for each $[\z]\in\bar\MM_{k+\ell}$ and
$[\tilde\z]:=\pi_\ell[\z]$ a morphism
$$
   \pi_\ell=(\tau,\{\phi_\alpha\}_{\alpha\in
   T}):\Sigma_\z\mapsto\Sigma_{\tilde\z}
$$
with the following property: There exists a collection of subtrees
$T'\subset T$
such that $\phi_\alpha$ is constant for $\alpha\in T'$ and biholomorphic
otherwise.
\end{lemma}

\section{Coherent almost complex structures}\label{sec:coh}

{\bf Coherent maps. }
Now we consider the Deligne-Mumford space $\bar\MM_{k+1}$ with $k+1$
marked points $z_0,\dots,z_k$. In the following discussion, the point
$z_0$ plays a special role (it will be the variable for holomorphic
maps in later sections).

Given a stable $(k+1)$-labelled tree $T$, we 
define an equivalence relation on
$\{0,\dots,k\}$ by $i\sim j$ iff $z_{\alpha_0
i}=z_{\alpha_0j}$. The equivalence classes yield a decomposition
$$
   \{0,\dots,k\} = I_0\cup\dots\cup I_\ell.
$$
Note that the marked points on $S_{\alpha_0}$ correspond to
equivalence classes consisting of one element; in particular, we may
put $I_0:=\{0\}$. Stability implies $\ell+1=n_{\alpha_0}\geq
3$. Conversely, we call a decomposition $\I=(I_0,\dots,I_\ell)$ of
$\{0,\dots,k\}$ {\em stable} if $I_0=\{0\}$ and $|\I|:=\ell+1\geq 3$.  
We will always order the $I_j$ such that the integers
$i_j:=\min\{i\mid i\in I_j\}$ satisfy
$$
   0=i_0<i_1<\dots<i_\ell. 
$$ 
Denote by $\MM_\I=\MM_{(I_0,\dots,I_\ell)}\subset\bar\MM_{k+1}$ the
union over those stable trees that give rise to the stable
decomposition $\I=(I_0,\dots,I_\ell)$. The
$\MM_\I$ are submanifolds of $\bar\MM_{k+1}$ with
$$
   \bar\MM_{k+1} = \bigcup_\I\MM_\I,
$$
and the closure of $\MM_\J$ is a union of certain
strata $\MM_\I$ with $|\I|\leq|\J|$. The above ordering of the $I_j$ 
determines a projection
$$
   p_\I:\MM_\I\to\MM_{|\I|},
$$
sending a stable curve $\z$ to the special points on the component
$S_{\alpha_0}$.  

\begin{definition}\label{def:coherent}
Let $Z$ be a Banach space and $k\geq 3$. We call a continuous map 
$F:\bar\MM_{k+1}\to Z$ {\em coherent} if it satisfies
the following two conditions:

(a) $F\equiv 0$ in a neighbourhood of those $\MM_\I$ with $|\I|=3$.

(b) For every stable decomposition $\I$ with $|\I|\geq 4$ there
exists a smooth map $F_\I:\MM_{|\I|}\to Z$ such that
$$
   F|_{\MM_\I}=F_\I\circ p_\I:\MM_\I\to Z.
$$
More generally, let $Z^*\subset Z$ be an open neighbourhood of $0$ and
let $\II$ be a
collection of stable decompositions. Then we call a continuous map 
$F:\cup_{\I\in\II}\MM_\I\to Z^*$ coherent if it satisfies conditions
(a) and (b) and in addition

(c) The image of $F$ is contractible in $Z^*$. 
\end{definition}

We denote the space of coherent maps from $\bar\MM_{k+1}$ to $Z$
resp.~$\cup_{\I\in\II}\MM_\I\to Z^*$ by
$$
    \Coh(\bar\MM_{k+1},Z) \quad\text{resp.}\quad
    \Coh(\cup_{\I\in\II}\MM_\I,Z^*). 
$$
It is equipped with the $C^0$-topology on $\bar\MM_{k+1}$ and the
$C^\infty$-topology on each $\MM_\I$ via the projection $p_\I$. 

\begin{remark}
Condition (a) ensures that almost complex structures
which depend in a coherent way on the domain will coincide with a
structure which is independent of the domain near double points.
It is included to keep the analysis of gluing holomorphic curves as
simple as possible (although gluing is not studied in this paper). 
\end{remark}

Our constructions of coherent maps will be based on the following

\begin{prop}\label{prop:coh}
Let $Z$ be a Banach space and $k\geq 3$. For each stable decomposition
$\I$ with $|\I|\geq 4$ fix an open subset
$U_\I\subset\MM_{|\I|}$ with compact closure. Then there exist
continuous linear extension maps 
$$
   E_\I:C_0^\infty(U_\I,Z)\to\Coh(\bar\MM_{k+1},Z)
$$
such that 
$$
   (E_\I\xi)|_{\MM_\I} = \xi\circ p_\I
$$
for every $\xi\in C_0^\infty(U_\I,Z)$. Moreover, the extension
operators can be chosen such that $E_\I\xi$ and $E_\J\eta$ have
disjoint supports for all $\xi\in
C_0^\infty(U_\I,Z)$ and $\eta\in C_0^\infty(U_\J,Z)$ with
$\I\neq\J$. 
If $Z^*\subset Z$ is open and star-shaped with respect to
the origin, then $E_\I$ maps $C_0^\infty(U_\I,Z^*)$ to
$\Coh(\bar\MM_{k+1},Z^*)$. 
\end{prop}

The proof of Proposition~\ref{prop:coh} will be given below. We first
give a reformulation of the coherency condition.

{\bf Cross ratios. }
For pairwise distinct points $z_0,z_1,z_2,z_3\in S^2=\C\cup\{\infty\}$
we have their {\em cross ratio}
$$
   w(z_0,z_1,z_2,z_3) := \frac{(z_1-z_2)(z_3-z_0)}{(z_0-z_1)(z_2-z_3)}
   \in S^2\setminus\{0,1,\infty\}.
$$
It extends to the cases where two (but no three) of the $z_i$ are
equal with the values $0,1,\infty$. For $0\leq i<j<m<n\leq k$ define
cross ratios $w_{ijmn}(\z):=w(z_i,z_j,z_m,z_n)$ on $\MM_{k+1}$. 
We will need the following two facts, of which the first one is
obvious and the second and third ones are proved in~\cite{MS}. 

{\bf Fact 1. }The cross ratios $w_{0,1,2,j}$ with $3\leq j\leq k$
define an embedding 
$$
   \phi_{k+1}:\MM_{k+1}\into (S^2)^{k-2}.
$$
{\bf Fact 2. }Each cross ratio $w_{ijmn}$, $0\leq i<j<m<n\leq k$,
extends to a continuous function $w_{ijmn}:\bar\MM_{k+1}\to S^2$; on a
stable curve $\z$ modelled over a tree $T$ it is given
by $w_{ijmn}(\z):=w(z_{\alpha i},z_{\alpha j},z_{\alpha m},z_{\alpha
n})$, where $\alpha\in T$ is the unique vertex for which no three of
the points $z_{\alpha i},z_{\alpha j},z_{\alpha m},z_{\alpha n}$ are
equal. 

{\bf Fact 3. }$\bar\MM_{k+1}\setminus\MM_{k+1}$ is precisely the locus
where some cross ratio $w_{0ijm}$, $1\leq i<j<m\leq k$, takes values
in $\{0,1,\infty\}$. 

\begin{remark}
In fact, it is shown in~\cite{MS} that the cross ratios
$w_{ijmn}$ for $0\leq i<j<m<n\leq k$ define an
embedding $\bar\MM_{k+1}\into (S^2)^{\binom{k+1}{4}}$. 
But we will not use this result in this paper. 
\end{remark}

Let $\I=(I_0,\dots,I_\ell)$ be a stable decomposition.
Recall that $i,j$ are equivalent iff they belong to the
same $I_a$ for some $1\leq a\leq\ell$.
We say that with respect to the decomposition $\I$ a triple $(i,j,m)$ 
with $1\leq i<j<m\leq k$ is of
\begin{itemize}
\item type I: if $i,j,m$ are pairwise nonequivalent;
\item type II: if exactly two of the $i,j,m$ are equivalent;
\item type III: if $i,j,m$ are all equivalent.
\end{itemize}
We define the type of a cross ratio $\w_{0,i,j,m}$ as the type of
the triple $(i,j,m)$. For two stable decompositions $\I,\J$ we say
that $\J$ is a {\em refinement} of $\I$ if $\J$-equivalence implies
$\I$-equivalence, in other words every $J_b$ is contained in some
$I_a$. Note that $\J$ is a refinement of $\I$ if and only if
$\MM_\I\subset\bar\MM_\J$.

\begin{lemma}\label{lem:types}
(a) The type of a triple $(i,j,m)$ does not decrease in the Gromov
limit on $\bar\MM_{k+1}$.

(b) A map $F:\MM_{k+1}\to Z$ satisfies condition (b) in
Definition~\ref{def:coherent} if and only if 
$F|_{\MM_\I}$ is a smooth function of the type I cross ratios for
every stable decomposition $\I$. 

(c) For stable decompositions $\I\neq \J$ such that $\J$ is a
refinement of $\I$ there exists a cross ratio that is of type I for
$\J$ and type II for $\I$. 

(d) For stable decompositions $\I\neq \J$ such that  $\J$ is not a
refinement of $\I$ and $|\I|\geq 4$ there exists a cross ratio that is
of type I for $\I$ and type II for $\J$. 
\end{lemma}

\begin{proof}
(a) follows immediately from the definition Gromov compactness,
see~\cite{MS}. 

For (b) first suppose that $F$ satisfies condition (b) in
Definition~\ref{def:coherent} and $\I$ is a stable decomposition with
$|\I|=\ell+1\geq 4$, so $F|_{\MM_\I}=F_\I\circ p_\I:\MM_\I\to Z$ for a
smooth map $F_\I:\MM_{\ell+1}\to Z$. Denote points in $\MM_{\ell+1}$
by $(z_0,\dots,z_\ell)$. 
In view of Fact 1 above, $F_\I$ can be uniquely written as a 
smooth function of the cross ratios $\w_{0,1,2,j}$ with $3\leq
j\leq\ell$. By definition of $p_\I$ this means that $F|_{\MM_\I}$ is a
smooth function of the type I cross ratios $\w_{i_0,i_1,i_2,i_j}$,
where $i_j$ is the minimal element in $I_j$. The converse follows
similarly. 

For (c) 
pick indices $i\neq j$ that are equivalent for $\I$ and not
equivalent for $\J$. Let $m$ be an index which is not equivalent to
$0$ or $i$ for $\I$. Then the cross ratio $w_{0,i,j,m}$ is of type I
for $\J$ and type II for $\I$.  

For (d) pick a $J_a$ that is not contained in any $I_b$. Write
$\I=(I_0,\dots,I_\ell)$ with $\ell\geq 3$. Set
$B:=\{b\in\{1,\dots,\ell\}\mid J_a\cap I_b\neq\emptyset\}$ and denote
by $\{J_a\cap I_b\}_{b\in B}$ the induced decomposition of $J_a$. Note
that $|B|\geq 2$. Now we distinguish two cases. 
If $B\neq\{1,\dots,\ell\}$ pick $d\in\{1,\dots,\ell\}\setminus B$ and
$b,c\in B$ with $b\neq c$, and pick $m\in I_d$, $i\in J_a\cap I_b$ and
$j\in J_a\cap I_c$. 
If $B=\{1,\dots,\ell\}$ pick pairwise distinct $b,c,d\in B$ such that
$I_d$ is not contained in $J_a$ (this is possible since $\ell\geq 3$
and $|\J|\geq 3$), and pick $m\in I_d\setminus J_a$, $i\in J_a\cap
I_b$ and $j\in J_a\cap I_c$. In both cases the cross ratio
$w_{0,i,j,m}$ is of type I for $\I$ and type II for $\J$ (with $i\sim
j$). 
\end{proof}

\begin{proof}[Proof of Proposition~\ref{prop:coh}]
Fix a union $W$ of three disjoint disks around $0,1,\infty$ in $S^2$
and a cutoff function $\chi:S^2\to[0,1]$ with support in $W$ which
equals $1$ near $0,1,\infty$. In view of Fact 3 above, we can (and
will) choose $W$ disjoint from the compact sets $w_{0,i,j,m}(\bar
U_\I)$ for all type I cross ratios $w_{0,i,j,m}$, $1\le i,j,m\le |\I|$
associated to stable decompositions $\I$.  

For $3\leq\ell\leq k$ denote by $\phi_{\ell+1}:\MM_{\ell+1}\into
(S^2)^{\ell-2}$ the embedding given by the cross ratios $w_{0,1,2,j}$
with $3\leq j\leq\ell$ (see Fact 1 above). In view of Fact 2 above,
each stable decomposition $\I=(I_0,\dots,I_\ell)$ induces a continuous
map
$$
   \phi_\I:\bar\MM_{k+1}\into (S^2)^{\ell-2}
$$
given by the cross ratios $w_{0,i_1,i_2,i_j}$ with $3\leq
j\leq\ell$. By construction we have 
$\phi_\I|_{\MM_\I}=\phi_{\ell+1}\circ p_\I$
and $\phi_{\ell+1}(\bar U_\I)\subset (S^2\setminus W)^{\ell-2}$. The
situation is summarized in the following diagram, where the first
vertical arrow is the inclusion and the second one an embedding:
\begin{equation*}
\begin{CD}
   \bar\MM_{k+1} @>{\phi_\I}>> (S^2)^{\ell-2}\supset (S^2\setminus
   W)^{\ell-2} \\ 
   @AAA @AA{\phi_{\ell+1}}A \\
   \MM_\I @>{p_\I}>> \MM_{|\I|}\supset\bar U_\I 
\end{CD}
\end{equation*}
Composition yields a continuous linear embedding
$$
   C_0^\infty(U_\I)\to C_0^0(\bar\MM_{k+1}),\qquad \xi\mapsto
   \xi\circ\phi_{\ell+1}^{-1}\circ\phi_\I
$$
satisfying $\xi\circ\phi_{\ell+1}^{-1}\circ\phi_\I|_{\MM_\I} =
   \xi\circ p_\I$. 
Define a continuous linear extension map
$E_\I:C_0^\infty(U_\I,Z)\to C^0(\bar\MM_{k+1},Z)$ by
$$
   E_\I\xi :=
   \left[\prod_{(i,j,m)\text{ type II }\atop \text{ with respect to
   }\I}\chi\circ w_{0,i,j,m}\right] \cdot
   \xi\circ\phi_{\ell+1}^{-1}\circ\phi_\I. 
$$
By construction, it satisfies  
$$
   (E_\I\xi)|_{\MM_\I} = \xi\circ\phi_{\ell+1}^{-1}\circ\phi_\I =
   \xi\circ p_\I,
$$
so this restriction is continuous with respect to the
$C^\infty$-topology via $p_\I$. Moreover, $E_\I\xi$ vanishes on the
locus where some type I cross ratio with respect to $W$ takes values
in $W$, so we could also write it as
$$
   E_\I\xi =
   \left[\prod_{(i,j,m)\text{ type II }\atop \text{ with respect to
   }\I}\chi\circ w_{0,i,j,m}\right] \cdot
   \left[\prod_{(i,j,m)\text{ type I }\atop \text{ with respect to 
   }\I}(1-\chi)\circ w_{0,i,j,m}\right] \cdot
   \xi\circ\phi_{\ell+1}^{-1}\circ\phi_\I. 
$$
Consider now a stable decomposition $\J\neq\I$. Suppose first that $\J$ is
not a refinement of $\I$. By Lemma~\ref{lem:types} (d) there exists a
cross ratio $w_{0,i,j,k}$ that is of type I for $\I$ and type II for
$\J$. It follows that $w_{0,i,j,k}\in\{0,1,\infty\}$ on $\MM_\J$, so 
by the preceding discussion $E_\I\xi$ vanishes on $\MM_\J$. Next
suppose that $\J$ is a refinement of $\I$,
i.e.~$\MM_\I\subset\bar\MM_\J$. Note that $E_\I\xi$ depends only on
cross ratios of  
types I and II with respect to $\I$. By Lemma~\ref{lem:types} (a),
these cross ratios do not become type III with respect to $\J$, so
$E_\I\xi|_{\MM_\J}$ is a smooth function of type I and II cross ratios
with respect to $\J$. But type II cross ratios are constant on
$\MM_\J$, so $E_\I\xi|_{\MM_\J}$ is a function of type I cross ratios
and thus coherent by Lemma~\ref{lem:types} (b). 

Next consider $\xi\in C_0^\infty(U_\I,Z)$ and $\eta\in
C_0^\infty(U_\J,Z)$ with $\I\neq\J$. Suppose first that $\J$ is
not a refinement of $\I$. By Lemma~\ref{lem:types} (d) there exists a
cross ratio $w_{0,i,j,k}$ that is of type I for $\I$ and type II for
$\J$. By construction, the support of $E_\J\eta$ is
contained in the set $w_{0,i,j,m}^{-1}(W)$ and $E_\I\xi$ vanishes on
this set, so $E_\I\xi$ and $E_\J\eta$ have disjoint support. 
Next suppose that $\J$ is a refinement of $\I$
By Lemma~\ref{lem:types} (c) 
there exists a cross ratio $w_{0,i,j,m}$ that is of type I for $\J$
and type II for $\I$. By construction, the support of $E_\I\xi$ is
contained in the set $w_{0,i,j,m}^{-1}(W)$ and $E_\J\eta$ vanishes on
this set, so again this shows that $E_\I\xi$ and $E_\J\eta$ have
disjoint support.

Finally, the formula for the extension operator $E_\I$ in terms of
cutoff functions shows that it preserves a star-shaped set
$\Z^*$. This concludes the proof of Proposition~\ref{prop:coh}.  
\end{proof}

Next we discuss coherent maps on certain subsets of $\bar\MM_{k+1}$. 
Denote by $\pi:\bar\MM_{k+1}\to\bar\MM_{k}$ the projection forgetting
the marked point $z_0$. Let $T$ be a stable $k$-labelled tree and
$S\subset T$ any subset. Denote by $\pi^{-1}\bar\MM_S\subset
\pi^{-1}\bar\MM_T\subset\bar\MM_{k+1}$ the closure of the subset of
$\pi^{-1}\bar\MM_T$ where the point $z_0$ lies on a component
belonging to $S$. Then  
$$
   \pi^{-1}\bar\MM_T = \cup_{\I\in\II_T}\MM_\I,\qquad 
   \pi^{-1}\bar\MM_S = \cup_{\I\in\II_S}\MM_\I
$$
for suitable collections $\II_S\subset\II_T$ of stable decompositions,
so we can speak of coherent maps on these sets. Clearly,
restriction defines a continuous linear map 
$$
   R_S:\Coh(\bar\MM_{k+1},Z)\to
   \Coh(\pi^{-1}\bar\MM_S,Z).  
$$ 

\begin{prop}\label{prop:coh2}
For every stable $k$-labelled tree $T$ and subset $S\subset T$ there
exists a continuous linear extension map 
$$
   E_S:\Coh(\pi^{-1}\bar\MM_S,Z)\to
   \Coh(\bar\MM_{k+1},Z)  
$$
such that 
$$
   R_S\circ E_S=\id.
$$
Moreover, the extension operator can be chosen such that its image
consists of maps with support in a small neighbourhood of
$\pi^{-1}\bar\MM_S$ and vanishes on $\pi^{-1}\bar\MM_T\setminus
\pi^{-1}\bar\MM_S$.  
If $Z^*\subset Z$ is open and star-shaped with respect to the origin,
then the extension operator preserves maps with values in $Z^*$. 
\end{prop}

\begin{proof}
The proof is very similar to the one of Proposition~\ref{prop:coh},   
so we only indicate the necessary modifications.
First note that the vanishing on $\pi^{-1}\bar\MM_T\setminus
\pi^{-1}\bar\MM_S$ allows us to extend the restrictions to vertices
$\alpha\in S$ one at a time, so it suffices to consider the case
that $S=\{\alpha\}$ consists of a single vertex $\alpha\in T$. Let
$\I_\alpha$ be the stable decomposition corresponding to the
$(k+1)$-labelled tree obtained from $T$ by adding the point $z_0$ on
the vertex $\alpha$. Then any $\xi\in\Coh(\pi^{-1}\bar\MM_{\{\alpha\}},Z)$
can be written as a unique function of suitable type I cross ratios
with respect to $\I_\alpha$ and we denote its canonical extension to
$\bar\MM_{k+1}$ by the same letter. 

Consider a vertex $\beta\in T$ with $\beta\neq\alpha$. Fix an
$m\in\{1,\dots,k\}$ with $z_{\beta m}\neq z_{\beta\alpha}$ and
consider all pairs of $i,j\in\{1,\dots,k\}$ such that $z_{\alpha
  i},z_{\alpha j},z_{\alpha m}$ are pairwise distinct. Then all the
corresponding cross ratios $w_{0,i,j,m}$ vanish on
$\pi^{-1}\bar\MM_{\{\beta\}}$. On the other hand, the locus in
$\pi^{-1}\bar\MM_{\{\alpha\}}$ where all these $w_{0,i,j,m}$ vanish
corresponds to configurations where $z_0$ belongs to a new component
$\gamma$ with precisely 3 special points
$z_0,z_{\gamma\alpha},z_{\gamma\beta}$, on which $\xi$ vanishes by
property (a) in Definition~\ref{def:coherent}. This proves that the
extended function $\xi$ vanishes on $\pi^{-1}\bar\MM_T\setminus
\pi^{-1}\bar\MM_{\{\alpha\}}$. Now as in the proof of
Proposition~\ref{prop:coh} we multiply $\xi$ by cutoff functions
of the type II cross ratios with respect to $\I_\alpha$ to obtain the
desired extension operator $\xi\mapsto E_S\xi$. 
\end{proof}

\begin{cor}\label{cor:coh2}
For every stable $k$-labelled tree $T$ and subset $S\subset T$ the
preimage of a Baire set (countable intersection of open dense sets)
under the restriction map $R_S:\Coh(\bar\MM_{k+1},Z)\to
\Coh(\pi^{-1}\bar\MM_S,Z)$ is again a Baire set.    
\end{cor}

\begin{proof}
This follows from Proposition~\ref{prop:coh2} and the following
general fact: If $f:X\to Y$ and is a continuous map with continuous
right inverse $g:Y\to X$, then preimages of Baire sets under $f$ are
Baire sets. Indeed, since preimages of open sets are open by
continuity, we only need to show that the preimage of a dense set
$A\subset Y$ is dense in $X$. So let $U\subset X$ be open. Then
$f(U)=g^{-1}(U)$ is open, thus $A\cap f(U)$ is nonempty, and
surjectivity of $f$ implies that $f^{-1}\bigl(A\cap
f(U)\bigr)=f^{-1}(A)\cap U$ is also nonempty. 
 \end{proof}

\begin{remark}[Collapsing subtrees]\label{rem:collapsing}
We will apply the Corollary~\ref{cor:coh2} in the following situation.  
Let $T$ be a stable $k$-labelled tree
and $T'\subset T$ a subtree. Let $\bar T$ be a $\bar k$-labelled
tree obtained by removing all but one of the marked points on $T'$ and
stabilizing. 
Let $\pi:\bar\MM_{k+1}\to\bar\MM_{k}$ and  $\bar\pi:\bar\MM_{\bar
 k+1}\to\bar\MM_{\bar k}$ denote the two projections. Then each
$F\in\Coh(\pi^{-1}\bar\MM_T,Z)$ induces an
$\bar F\in\Coh(\bar\pi^{-1}\bar\MM_{\bar T},Z)$ which equals the
restriction of $F$ on $\pi^{-1}\bar\MM_{T\setminus T'}$ and zero on
$\pi^{-1}\bar\MM_{T'}$. 
This defines a coherent map $\bar F$ due to condition (a) of
Definition~\ref{def:coherent}.  
\end{remark}

{\bf Almost complex structures. } For a symplectic vector space
$(V,\om)$, denote by $\J(V,\om)$ the space of $\om$-tamed complex
structures, i.e., endomorphisms $J:V\to V$ such that $J^2=-\id$ and
$\om(v,Jv)>0$ for all vectors $v$. We equip $V$ with a euclidean structure $g$.
Then the space
$\J(V,\om)$ is a Riemannian manifold with tangent space
$$
   T_{J}\J(V,\om) = \{Y\in\End(V)\mid \ JYJ=Y\}.
$$
and metric given by $\mbox{trace}(Y^tY)$ for instance.
The exponential map defines a diffeomorphism
$$
\exp_J:B(0,\rho(g,J))\subset T_{J}\J(V,\om)\longrightarrow \J(V,\om)
$$
from the open ball of radius $\rho(g,J)>0$ which continuously depends on $g$ and
$J$.



Following~\cite{Fl}, for a vector bundle $E\to X$ over a closed
manifold $X$ we denote the space of Floer's $C^\eps$-sections in $E$ by
$$
   C^\eps(X,E) := \{s\in\Om(X,E)\mid \sum_{i=1}^\infty
   \eps_i\|s\|_{C^i}<\infty\}.  
$$
Here $\eps=(\eps_i)_{i\in\N}$ is a fixed sequence of positive numbers
and $\|\ \|_{C^i}$ is the $C^i$-norm with respect to some connection
on $E$. It is shown in Lemma 5.1 of~\cite{Fl} that if the $\eps_i$ converge
sufficiently fast to zero, then $C^\eps(X,E)$ is a Banach space
consisting of smooth sections and containing sections with support in
arbitrarily small sets in $X$. 

Now let $(X,\om)$ be a closed symplectic manifold. Fix a tamed
almost complex structure $J_0$ on $(X,\om)$, i.e.~a smooth section in
the bundle $\J(TX,\om)\to X$ with fibres $\J(T_xX,\om_x)$. Define
a Riemannian metric on $X$ via 
$$
g(v,w):= \om(v,Jw) + \om(w,Jv).
$$

Let $T_{J_0}\J(TX,\om)\to X$ be the vector bundle with fibres
$T_{J_0(x)}\J(T_xX,\om_x)$ and set
\begin{gather*}
   T_{J_0}\JJ := C^\eps\bigl(X,T_{J_0}\J(TX,\om)\bigr), \cr
   \JJ := \JJ(X,\om) := \exp_{J_0}(B) 
\end{gather*}
where $B:= \{Y\in T_{J_0}\JJ \mid Y(x)\in B(0,\rho(g(x),J_0(x)))\}$.
Thus $\JJ$ is the space of tamed almost complex structures of
$(X,\om)$ that are of class $C^\eps$ over $J_0$ via $\exp_{J_0}$. 
We think of $\JJ$ as a Banach manifold with a single chart
$\exp_{J_0}$. 

\begin{remark}
We cannot replace the space $C^\eps$ by $C^\infty$ because the latter
one is not a Banach space and we need to apply the implicit function
theorem below. However, we could replace $C^\eps$ by $C^N$ for
sufficiently large $N$ (which amounts to choosing $\eps_i=0$ for
$i>N$); then the moduli spaces of holomorphic spheres would not be
smooth but of class $C^N$ as well, which suffices for the
definition of Gromov-Witten invariants. 
\end{remark}

More generally, let $P$ be any manifold. A {\em tamed almost
complex structure on $(X,\om)$ parametrized by $P$} is a smooth
section in the pullback bundle $\J(TX,\om)\to P\times X$. Fix $J_0$ as
above.
Let $T_{J_0}\J_P(TX,\om)\to P\times X$ be the vector bundle
with fibres $T_{J_0(p,x)}\J(T_xX,\om_x)$ and set
\begin{gather*}
   T_{J_0}\JJ_P := C^\eps\bigl(P\times X,T_{J_0}\J_P(TX,\om)\bigr), \cr
   \JJ_P := \JJ_P(X,\om) := \exp_{J_0}(B_P),
\end{gather*}
where $B_P:= \{Y\in T_{J_0}\JJ_P \mid Y(p,x)\in B(0,\rho(g(x),J_0(x)))\}$.
Note that we may think of $J\in\JJ_P$ as a map $P\to\JJ$. 
For an open subset $U\subset
P$, we denote by $T_J\JJ_U\subset T_J\JJ_P$ the subspace of those
sections that have compact support in $U$. We will be interested in
the spaces
$$
   \JJ_{S^2}=\JJ_{S^2}(X,\om)\quad\text{and}\quad
   \JJ_{\bar\MM_{k+1}} = \JJ_{\bar\MM_{k+1}}(X,\om)
$$
parametrized by the Riemann sphere and the Deligne-Mumford space
with $k+1$ marked points, respectively. To better understand the
second case, recall that we have a holomorphic projection
$\pi:\bar\MM_{k+1}\to\bar\MM_k$. Fix a stable curve
$\z=(\{z_{\alpha\beta}\}_{\alpha E\beta}, \{\z_i\}_{1\leq i\leq
k})\in\bar\MM_k$, modelled over the $k$-labelled tree $T$. Recall
that $\pi^{-1}[\z]$ is naturally identified with the nodal Riemann
surface $\Sigma_\z$ which is a union of $|T|$ copies $S_\alpha$ of
$S^2$, glued together at the points $z_{\alpha\beta}$. Restriction
of $J\in\JJ_{\bar\MM_{k+1}}$ to $\pi^{-1}[\z]$ yields a continuous map
$$
   J_\z:\Sigma_\z\to\JJ
$$
which is smooth on each $S_\alpha$.
We define the space $\JJ_{k+1}$ of {\em coherent almost complex
  structures} by 
\begin{gather*}
   T_{J_0}\JJ_{k+1} := \Coh(\bar\MM_{k+1},T_{J_0}\JJ)\subset
   T_{J_0}\JJ_{\bar\MM_{k+1}}, \cr
   \JJ_{k+1} := \exp_{J_0}(T_{J_0}\JJ_{k+1}) \subset
   \JJ_{\bar\MM_{k+1}}, 
\end{gather*}
where $\Coh(\bar\MM_{k+1},T_{J_0}\JJ)\subset
T_{J_0}\JJ_{\bar\MM_{k+1}}$ denotes the subspace of those maps
$Y:\bar\MM_{k+1}\to T_{J_0}\JJ$ satisfying the conditions of
Definition~\ref{def:coherent}. 

\begin{lemma}\label{lem:J-smooth}
For $\z\in\bar\MM_k$ modelled over the tree $T$ and $J\in\JJ_{k+1}$
the following holds. 

(a) The restriction of $J_\z$ to each component $S_\alpha$ of
$\Sigma_\z$ defines a smooth tamed almost complex structure
$J_\alpha\in\JJ_{S_\alpha}$. Similarly, any $Y\in T_J\JJ_{k+1}$
restricts on $S_\alpha$ to a smooth section $Y|_{S_\alpha}\in
T_{J_\alpha}\JJ_{S_\alpha}$.

(b) Conversely, let $U_\alpha\subset S_\alpha\setminus SP_\alpha$
be open subsets (for $\alpha\in T$), where $SP_\alpha$ is the set of
special points on $S_\alpha$ (except $z_0$).
Then for any collection of sections $Y_\alpha\in
T_{J_\alpha}\JJ_{U_\alpha}$ with compact support in $U_\alpha$,
there exists a $Y\in T_J\JJ_{k+1}$ with $Y|_{S_\alpha}=Y_\alpha$.
Moreover, this extension is given by an injective continuous linear
map 
$$
\oplus_\alpha T_{J_\alpha}\JJ_{U_\alpha}
\longrightarrow T_J\JJ_{k+1}
$$
\end{lemma}

\begin{proof}
(a) follows immediately from the fact (Lemma~\ref{lem:fibre}) that
the $S_\alpha$ are submanifolds of $\bar\MM_{k+1}$. The same fact
also implies (b):
Since $U_\alpha\times X\subset\MM_{n_\alpha+1}\times X$ is a
submanifold, the section $Y_\alpha\in
T_{J_\alpha}\JJ_{U_\alpha}\subset T_{J_\alpha}\JJ_{S_\alpha}$ of
the vector bundle $T_J\J_{S_\alpha}(TX,\om)\to S_\alpha\times X$
extends to a $C^\eps$-section $\tilde Y_\alpha$ of the vector bundle
$T_J\J_{\MM_{n_\alpha+1}}(TX,\om)\to \MM_{n_\alpha+1}\times X$
with compact support. Now the construction proceeds exactly as in
the proof of Proposition~\ref{prop:coh}, choosing the cutoff function
$\chi$ to be a $C^\eps$-function. 
\end{proof}

\section{Holomorphic maps}\label{sec:holo}

{\bf $J$-holomorphic maps. }
To a map $f:S^2\to X$ and a tamed almost complex structure
$J\in\JJ_{S^2}$ we associate the $(0,1)$--form
$$
   \bar\p_Jf := \frac{1}{2}\Bigl(df+J(z,f)\circ df\circ
   j\Bigr)
$$
which at the point $z\in S^2$ is given by
$$
   \bar\p_Jf(z) := \frac{1}{2}\Bigl(df(z)+J
   \bigl(z,f(z)\bigr)\circ df(z)\circ j(z)\Bigr).
$$
We call $f$ $J$-holomorphic if $\bar\p_Jf=0$. Then its {\em energy}
is given by
$$
   E(f) := \frac{1}{2}\int_{S^2}|df|^2{\rm dvol} = \int_{S^2}f^*\om.
$$
Fix an integer $m\geq 1$ and a real number $p>1$ such that $mp>2$
and define the Banach manifold
$$
   \BB := \BB^{m,p}:= W^{m,p}(S^2,X)
$$
of maps of Sobolev class $W^{m,p}$ (which are continuous because
$mp>2$). Let
$$
   \EE := \EE^{m-1,p}\to\BB
$$
be the Banach bundle whose fibre at $f\in\BB$ is given by
$\EE_f=W^{m-1,p}\bigl(S^2\!,\!\Om^{0,1}\!(f^*T\!X)\bigr)$. Then the
Cauchy-Riemann operator defines a smooth section
$$
   \pb_J:\BB\to\EE,\qquad f\mapsto\pb_Jf.
$$
More generally, the {\em universal Cauchy-Riemann operator} is the
section
$$
   \pb:\BB\times\JJ_{S^2}\to\EE,\qquad (f,J)\mapsto\pb_Jf.
$$
Now fix pairwise disjoint points $z_1,\dots,z_k\in S^2$.
Evaluation at these points defines a smooth map
$$
   \ev^k:\BB\to X^k,\qquad
   f\mapsto \bigl(f(z_1),\dots,f(z_k)\bigr).
$$

\begin{lemma}\label{lem:surj-nonconst}
Let $k\in\N$, $f\in\BB$ and $J\in\JJ_{S^2}$ satisfy $\pb_Jf=0$. If
$f$ is nonconstant, then the linearization of
$$
   (\pb,\ev^k):\BB\times\JJ_{S^2}\to\EE\times X^k
$$
at $(f,J)$ is surjective. More precisely, for any nonempty open
subset $U\subset S^2\setminus\{z_1,\dots,z_k\}$, the restriction of
the linearization to the subspace
$$
   T_f\BB\oplus T_J\JJ_U\subset T_f\BB\oplus T_J\JJ_{S^2}
$$
of sections with support in $U$ is surjective.
\end{lemma}

\begin{proof}
The proof is very similar to that of Proposition 3.4.2 in~\cite{MS},
so we only sketch the argument. The linearization of $(\pb,\ev^k)$ at
$(f,J)$, restricted to the subspace $T_f\BB^{m,p}\oplus T_J\JJ_U$,
is the linear operator
\begin{gather*}
   L:T_f\BB^{m,p}\oplus T_J\JJ_U\to\EE_f\oplus\bigoplus_{i=1}^kT_{f(z_i)}X,
   \cr
   (\xi,Y)\mapsto\Bigl(D_f\xi + \frac{1}{2}Y(z,f)\circ df\circ
   j,\xi(z_1),\dots,\xi(z_k)\Bigr),
\end{gather*}
where $D_f:T_f\BB^{m,p}\to\EE_f$ is the linearization of $\pb_J$ at
$f$. We will prove surjectivity of the operator
\begin{gather*}
   L_0:B_0^{m,p}\oplus T_J\JJ_U\to\EE_f, \cr
   (\xi,Y)\mapsto D_f\xi + \frac{1}{2}Y(z,f)\circ df\circ j,
\end{gather*}
where $B_0^{m,p}:=\{\xi\in T_f\BB^{m,p}\mid
\xi(z_1)=\dots=\xi(z_k)=0\}$. Surjectivity of $L_0$ easily implies
surjectivity of $L$: Given $(\eta,v_1,\dots,v_k)\in
\EE_f\oplus\bigoplus_{i=1}^kT_{f(z_i)}X$, pick a $\xi_1\in T_f
\BB^{m,p}$ with $\xi_1(z_i)=v_i$ for all $i$ and find
$(\xi_2,Y)\in B_0^{m,p}\oplus T_J\JJ_U$ with
$L_0(\xi,Y)=\eta-D_f\xi_1$, hence
$L(\xi_1+\xi_2,Y)=(\eta,v_1,\dots,v_k)$.

To prove surjectivity of $L_0$, first consider the case $m=1$. Since
$D_f$ is a Fredholm operator, the image of $L_0$ is closed. So it
suffices to show that every element $\eta\in\EE_f^*$ orthogonal to
$\im(L_0)$ is zero. Note that $\EE_f$ is a space of $L^p$-sections,
so its dual space $\EE_f^*$ can be identified with the appropriate
space of $L^q$-sections with $1/p+1/q=1$ via the $L^2$ inner
product. Thus orthogonality to $\im(L_0)$ amounts to the equations
\begin{gather*}
   \la D_f\xi,\eta\ra_{L^2} = 0,\cr
   \la Y(z,f)\circ df\circ j,\eta\ra_{L^2} = 0
\end{gather*}
for all $(\xi,Y)\in B_0^{1,p}\oplus T_J\JJ_U$. By Lemma 3.4.7
of~\cite{MS}, the first equation implies that $\eta$ is smooth on
$S^2\setminus\{z_0,\dots,z_k\}$ and satisfies on this set the
Cauchy-Riemann type equation $D_f^*\eta=0$. Since $f$ is
nonconstant, the set $U^*:=\{z\in U\mid df(z)\neq 0\}$ is open and
dense in $U$. Suppose that $\eta(z)\neq 0$ for some $z\in U^*$.
Then by Lemma 3.2.2 of~\cite{MS} and cutting off near $z$, we find
a $Y\in T_J\JJ_U$ such that $\la Y(z',f)\circ df(z')\circ
j,\eta(z')\ra\geq 0$ for all $z'$ and $>0$ at $z'=z$, in contradiction
to the second equation. (Here we have used the fact that $Y$ may
depend on points in $U$, which allows us to cut off $Y$ near $z$!)
It follows that $\eta(z)=0$ for all $z\in U^*$, hence by unique
continuation (see Lemma 3.4.7 of~\cite{MS}) $\eta\equiv 0$. This
proves surjectivity of $L_0$ in the case $m=1$. The case of
general $m\geq 1$ follows from the $m=1$ case by elliptic
regularity: Given $\eta\in\EE_f^{m-1,p}$, by the $m=1$ case we
find $(\xi,Y)\in B_0^{1,p}\oplus T_J\JJ_U$ with $L_0(\xi,Y)=\eta$,
and elliptic regularity (Theorem C.2.3 of~\cite{MS}) yields
$\xi\in B_0^{m,p}$. This concludes the proof of
Lemma~\ref{lem:surj-nonconst}.
\end{proof}

For constant $f$ we have a similar result with $k=1$.

\begin{lemma}\label{lem:surj-const}
If $f\in\BB$ is constant and $J\in\JJ_{S^2}$, then the linearization of
$$
   (\pb_J,\ev^1):\BB\to\EE\times X
$$
at $f$ is surjective.
\end{lemma}

\begin{proof}
Let $f\equiv x\in X$ and identify $T_xX$ with $\R^{2n}$. Then the
linearization of $\pb_J$ at $f$ is the Cauchy-Riemann operator
$$
   L:W^{m,p}(S^2,\R^{2n})\to
   W^{m-1,p}\bigl(\Om^{0,1}(S^2,\R^{2n})\bigr),\qquad
   \xi\mapsto\pb_J\xi
$$
for a family of linear complex structures $J(z)$ on $\R^{2n}$, $z\in
S^2$. Elements in the kernel of $L$ have zero energy (since they are
homologically trivial) and are therefore constant. Since $L$ has
index $2n$, it follows that $L$ is surjective and its kernel
consists of the constant functions. Hence the linearized evaluation
map $\xi\mapsto\xi(z_1)$ is surjective on the kernel of $L$ and the
result follows.
\end{proof}

Lemma~\ref{lem:surj-nonconst} and Lemma~\ref{lem:surj-const}
immediately imply

\begin{prop}\label{prop:surj-smooth}
The {\em universal solution space}
$$
   \tilde\MM(\JJ_{S^2}) := \{(f,J)\in\BB\times\JJ_{S^2}\mid \pb_Jf=0\}
$$
is a Banach manifold and the $1$-point evaluation
map
$$
   \ev^1:\tilde\MM(\JJ_{S^2})\to X
$$
is a submersion. Moreover, for any $k\in\N$ the $k$-point evaluation
map
$$
   \ev^k:\tilde\MM^\nc(\JJ_{S^2})\to X^k
$$
is a submersion on the open component
$\tilde\MM^\nc(\JJ_{S^2})\subset\tilde\MM(\JJ_{S^2})$ of pairs $(f,J)$
with $f$ nonconstant.
\end{prop}

For a homology class $A\in H_2(X;\Z)$, $J\in\JJ_{S^2}$ and a smooth
submanifold $Z\subset X^k$ define
\begin{gather*}
   \tilde\MM(A,J) := \{f:S^2\to X\mid \pb_Jf=0,\ [f]=A\}, \cr
   \tilde\MM(A,J;Z) := \{f\in\tilde\MM(A,J)\mid
   \bigl(f(z_1),\dots,f(z_k)\bigr)\in Z\}.
\end{gather*}
Taking regular values of the projection $(f,J)\mapsto J$, we obtain

\begin{cor}\label{cor:surj-smooth}
For every smooth submanifold $Z\subset X$ there exists a Baire set
$\JJ_{S^2}^\reg(Z)\subset\JJ_{S^2}$ such that for every homology
class $A\in H_2(X;\Z)$ and every $J\in\JJ_{S^2}^\reg(Z)$ the
solution space
$$
   \tilde\MM(A,J;Z) = \{f\in\tilde\MM(A,J)\mid
   f(z_1)\in Z\}
$$
is a smooth manifold of dimension
$$
   \dim\tilde\MM(A,J;Z) = 2n + 2c_1(A) - \codim_\R Z.
$$
More generally, for every smooth submanifold $Z\subset X^k$ there
exists a Baire set $\JJ_{S^2}^\reg(Z)\subset\JJ_{S^2}$ such that for 
every nonzero homology class $0\neq A\in H_2(X;\Z)$ and every
$J\in\JJ_{S^2}^\reg(Z)$ the solution space
$$
   \tilde\MM(A,J;Z) = \{f\in\tilde\MM(A,J)\mid
   \bigl(f(z_1),\dots,f(z_k)\bigr)\in Z\}
$$
is a smooth manifold of dimension
$$
   \dim\tilde\MM(A,J;Z) = 2n + 2c_1(A) - \codim_\R Z.
$$
\end{cor}

\section{Nodal maps}\label{sec:nodal}

Throughout this section, we fix a symplectic manifold $(X,\om)$ of
dimension $2n$, not necessarily closed. 

{\bf Weighted trees. }
A {\em weighted $k$-labelled tree} $(T,\{A_\alpha\})$ is a
$k$-labelled tree $T$ with a
collection of homology classes $A_\alpha\in H_2(X;\Z)$, $\alpha\in
T$. We call $(T,\{A_\alpha\})$ {\em weighted stable} if each component
$\alpha\in T$ with $A_\alpha=0$ carries at least 3 special
points. Note that stability implies weighted stability but not vice
versa. We call a component $\alpha$ with $A_\alpha=0$ a {\em ghost
component}, and a maximal subtree consisting of ghost components a
{\em ghost tree}. The {\em reduced index set} of $(T,\{A_\alpha\})$ is
the subset $R\subset\{1,\dots,k\}$ containing all marked points on
components with $A_\alpha\neq 0$, and the maximal marked point on each
ghost tree. 

{\bf Stable maps. }
Fix a $k$-labelled tree $T=(T,E,\Lambda)$ (not necessarily stable).
Let
$$
   \JJ_T := \prod_{\alpha\in T}\JJ_{S_\alpha}
$$
be the space of collections $\J=\{J_\alpha\}_{\alpha\in T}$ of
tamed almost complex structures $J_\alpha\in\JJ_{S_\alpha}$
parametrized by the spheres $S_\alpha$. Let
$\z=(\{z_{\alpha\beta}\}_{\alpha E\beta}, \{\z_i\}_{1\leq i\leq k})$
be a nodal curve modelled over $T$ and $\Sigma_\z=\cup_{\alpha\in
T}S_\alpha$ the associated nodal Riemann surface. A {\em continuous
map} $\f:\Sigma_\z\to X$ is a collection $\{f_\alpha\}_{\alpha\in
T}$ of continuous maps $f_\alpha:S_\alpha\to X$ such that
$f_\alpha(z_{\alpha\beta})=f_\beta(z_{\beta\alpha})$ whenever
$\alpha E\beta$. Note that each such map $\f$ induces homology
classes $A_\alpha=[f_\alpha]\in H_2(X;\Z)$ and thus weights on the
tree $T$. To $\J\in\JJ_T$ and $\f:\Sigma_\z\to X$ we
associate the Cauchy-Riemann operator
$$
   \bar\p_\J\f := \frac{1}{2}\Bigl(d\f+\J(z,\f)\circ d\f\circ
   j\Bigr)
$$
which at the point $z\in S_\alpha$ is given by
$$
   \bar\p_{J_\alpha}f_\alpha(z) := \frac{1}{2}\Bigl(df_\alpha(z)+J_\alpha
   \bigl(z,f_\alpha(z)\bigr)\circ df_\alpha(z)\circ j(z)\Bigr).
$$
A continuous map $\f$ is called {\em $\J$-holomorphic} if
$\bar\p_\J\f=0$ (hence each $f_\alpha$ is smooth by elliptic
regularity). A pair $(\z,\f)$ consisting of a nodal curve $\z$ and a
$\J$-holomorphic map $\f:\Sigma_\z\to X$ is called a {\em nodal
$\J$-holomorphic map of genus zero with $k$ marked points}, or simply a
{\em nodal map}. Note that a nodal map $(\z,\f)$ represents a
homology class
$$
   [\f]:=\sum_{\alpha\in T}[f_\alpha]\in H_2(X;\Z)
$$
and carries the energy
$$
   E(\f) := \sum_{\alpha\in T}E(f_\alpha).
$$
A nodal map $(\z,\f)$ is called {\em stable} if the underlying
weighted tree $(T,\{A_\alpha=[f_\alpha]\})$ is weighted stable,
i.e.~if each component $\alpha\in T$ for which $f_\alpha$ is constant
carries at least 3 special points. 

For a weighted tree $(T,\{A_\alpha\})$ and $\J\in\JJ_T$ let
$$
   \tilde\MM_T(\{A_\alpha\},\J)
$$
be the space of all nodal maps $(\z,\f)$ modelled over the tree $T$
with $[f_\alpha]=A_\alpha$ for all $\alpha\in T$. 
The space of {\em stable} maps modelled on the tree $T$ representing
the class $A$ is the disjoint union
$$
   \tilde\MM_T(A,\J) := \coprod_{\sum
   A_\alpha=A}\tilde\MM_T(\{A_\alpha\},\J),
$$
where the union is taken over all decompositions $\{A_\alpha\}$ of
$A$ such that $(T,\{A_\alpha\})$ is stable.

{\bf Domain-stable maps. }
Fix an integer $k\geq 3$ and consider the
space $\JJ_{k+1}$ of coherent tamed almost complex
structures on $(X,\om)$ parametrized by the Deligne-Mumford space
$\bar\MM_{k+1}$.  Fix a {\em stable curve}
$\z=(\{z_{\alpha\beta}\}_{\alpha E\beta}, \{\z_i\}_{1\leq i\leq
k})$, modelled over the $k$-labelled tree $T=(T,E,\Lambda)$.
Restriction of $J\in\JJ_{k+1}$ to
$\pi^{-1}[\z]\cong\Sigma_\z=\cup_{\alpha\in T}S_\alpha$ yields a
map
$$
   J_\z:\Sigma_\z\to\JJ.
$$
By Lemma~\ref{lem:J-smooth}, the restriction of $J_\z$ to each
component $S_\alpha$ of $\Sigma_\z$ is smooth, so $J_\z$ yields an
element in $\JJ_T=\prod_{\alpha\in T}\JJ_{S_\alpha}$. As above, to a
map $\f=\{f_\alpha\}_{\alpha\in T}:\Sigma_\z\to X$ and
$J\in\JJ_{k+1}$ we associate the Cauchy-Riemann operator
$$
   \bar\p_{J_\z}\f := \frac{1}{2}\Bigl(d\f+J_\z(z,\f)\circ d\f\circ
   j\Bigr)
$$
which at the point $z\in S_\alpha$ is given by
$$
   \bar\p_{J_\z}f_\alpha(z) := \frac{1}{2}\Bigl(df_\alpha(z)+J_\z
   \bigl(z,f_\alpha(z)\bigr)\circ df_\alpha(z)\circ j(z)\Bigr).
$$
If $\z$ is a stable curve and $\f$ is $J_\z$-holomorphic for
$J\in\JJ_{k+1}$, we call the pair $(\z,\f)$ a {\em domain-stable
$J$-holomorphic map of genus zero with $k$ marked points}, or simply
a {\em domain-stable map}. Note that every domain-stable map is
stable, but domain-stability is more restrictive because it requires
the underlying curve to be stable.

More generally, let $\z$ be a nodal curve with $k$ marked points,
not necessarily stable. By Lemma~\ref{lem:stab}, stabilization
yields a holomorphic map $\st:\Sigma_\z\to\Sigma_{\st(\z)}$. So
$J\in\JJ_{k+1}$ yields a map
$$
   J^\st_\z:\Sigma_\z\to\JJ,\qquad
   J^\st_\z(z):=J_{\st(\z)}\bigl(\st(z)\bigr)
$$
whose restriction to each sphere $S_\alpha$ is smooth. Note that
$J^\st_\z$ is constant on the spheres $S_\alpha$ with $\alpha$
in the collection of subtrees $T'\subset T$ from
Lemma~\ref{lem:stab}. This
construction allows us to define the Cauchy-Riemann operator
$\bar\p_{J^\st_\z}\f$ as above and speak about {\em nodal
$J$-holomorphic maps $(\z,\f)$} for $J\in\JJ_{k+1}$. As above, let
$$
   \tilde\MM_T(A,J) := \coprod_{\sum
   A_\alpha=A}\tilde\MM_T(\{A_\alpha\},J),
$$
be the space of stable $J$-holomorphic maps modelled on the
(not necessarily stable) $k$-labelled tree $T$ for $J\in\JJ_{k+1}$.

{\bf Moduli spaces. } Let $J\in\JJ_{k+1}$. An {\em isomorphism}
between two nodal $J$-holomorphic maps $(\z,\f)$ and
$(\tilde\z,\tilde\f)$ is an isomorphism
$(\tau,\{\phi_\alpha\}_{\alpha\in T})$ between the nodal curves $\z$
and $\tilde\z$ such that $\tilde
f_{\tau(\alpha)}\circ\phi_\alpha=f_\alpha$ for all $\alpha\in T$.
The following lemma shows that $J$-holomorphicity is preserved under
isomorphisms.

\begin{lemma}\label{lem:action}
Let $\z$, $\tilde\z$ be nodal curves with $k$ marked points and
$\f:\Sigma_\z\to X$, $\tilde\f:\Sigma_{\tilde\z}\to X$ be continuous
maps. Let $\phi=(\tau,\{\phi_\alpha\}_{\alpha\in T}):\z\to\tilde\z$
be an isomorphism between the nodal curves such that $\tilde
f_{\tau(\alpha)}\circ\phi_\alpha=f_\alpha$ for all $\alpha\in T$.
Let $J\in\JJ_{k+1}$. Then $\f$ is $J$-holomorphic if and only if
$\tilde\f$ is.
\end{lemma}

\begin{proof}
By definition of $\JJ_{k+1}$, $J$ is invariant under isomorphisms
$\psi$ in the sense that $J_{\psi(\w)}\bigl(\psi(w)\bigr) = J_\w(w)$
for $\w\in\tilde\MM_k$ and $w\in\Sigma_\w$. By Lemma~\ref{lem:stab},
the isomorphism $\phi:\z\to\tilde\z$ induces an isomorphism
$\phi^\st:\st(\z)\to\st(\tilde\z)$ between the stabilized curves such
that $\st\circ\phi=\phi^\st\circ\st$. It follows that
\begin{align}
   J^\st_{\phi(\z)}\bigl(\phi(z)\bigr) &=
   J_{\st\bigl(\phi(\z)\bigr)}\Bigl(\st\bigl(\phi(z)\bigr)\Bigr)
   = J_{\phi^\st\bigl(\st(\z)\bigr)}
   \Bigl(\phi^\st\bigl(\st(z)\bigr)\Bigr)\\
   &= J_{\st(\z)}\bigl(\st(z)\bigr) = J^\st_\z(z).
\end{align}
Together with holomorphicity of $\phi_\alpha$, this implies
$$
  \bar\p\tilde\f =
  \bar\p_{J^\st_{\phi(\z)}}(f_\alpha\circ\phi_\alpha^{-1})\bigl(\phi(z)\bigr)
   = \bar\p_{J^\st_\z}f_\alpha(z)\circ d\phi_\alpha(z)^{-1}
$$
for every $\alpha\in T$ and $z\in S_\alpha$, and the lemma follows.
\end{proof}

It follows that for every tree $T$ and homology class $A$ the group
$G_T$ of isomorphisms fixing $T$ (but not necessarily the identity on
$T$) acts on the space
$\tilde\MM_T(A,J)$ of stable maps. The action of $G_T$ is proper
(see Section 6.1 in~\cite{MS}; the argument is similar to the remark
in Section~\ref{sec:curves} above),
and due to the stability condition it has only finite isotropy
groups (see~\cite{MS}). Denote its quotient space by
$$
   \MM_T(A,J) := \tilde\MM_T(A,J)/G_T.
$$
The {\em moduli space of stable maps (of genus zero with $k$ marked
points representing the class $A$)} is the disjoint union
$$
   \bar\MM_k(A,J) = \bigcup_T\MM_T(A,J)
$$
over the isomorphism classes of $k$-labelled trees. It carries a
natural topology, the {\em Gromov topology} defined in~\cite{MS},
which makes it into a metrizable space. We call the $\MM_T(A,J)$ the
{\em strata} of $\bar\MM_k(A,J)$. In particular, we have the {\em
top stratum} $\MM_k(A,J):=\MM_{T_k}(A,J)$ of stable maps modelled
over the $k$-labelled tree $T_k$ with one vertex.

The {\em moduli space of domain-stable maps} is the subset
$$
   \bar\MM^\ds_k(A,J) := \bigcup_{T \text{ stable }}\MM_T(A,J)
   \subset \bar\MM_k(A,J)
$$
of those nodal maps whose underlying curve is stable. Note that if
the tree $T$ is stable, then the action of $G_T$ on
$\tilde\MM_T(A,J)$ is free and covers the trivial action on
$T$ (see~\cite{MS}). So in this case the quotient space has a natural
decomposition
$$
   \MM_T(A,J) = \coprod_{\sum
   A_\alpha=A}\MM_T(\{A_\alpha\},J), \qquad \MM_T(\{A_\alpha\},J) =
   \tilde\MM_T(\{A_\alpha\},J)/G_T.
$$
(If $T$ is unstable the corresponding decomposition need not descend
to the quotient because isomorphisms may interchange vertices of
$T$).

{\bf Compactness. }We now state Gromov's compactness theorem in
our context. It is proved in~\cite{MS} for
almost complex structures $J$ not depending on points in
$\bar\MM_{k+1}$, but the proof immediately carries over to our
situation.

\begin{thm}\label{thm:Gromov-comp}
Let $J\in\JJ_{k+1}$. Let $(\z^\nu,\f^\nu)$ be a sequence of stable
$J$-holomorphic maps (of genus zero with $k$ marked points) with
uniformly bounded energy $E(\f^\nu)\leq C$. Then a subsequence
converges in the Gromov topology to a stable $J$-holomorphic map
$(\z,\f)$ (of genus zero with $k$ marked points).
\end{thm}

For the precise definition of Gromov convergence see~\cite{MS}. It
implies, in particular, that after passing
to a subsequence the following properties hold:

(1) The $(\z^\nu,\f^\nu)$ all belong to the same moduli space
$\MM_{\tilde T}(\{A_{\tilde\alpha}\},J)$ for some stable weighted
tree $(\tilde T,\{A_{\tilde\alpha}\})$.

(2) The limit map belongs to some moduli space
$\MM_T(\{A_\alpha\},J)$ and there exists a surjective tree
homomorphism $\tau:T\to\tilde T$ such that
$$
   \tau(\alpha_i)=\tilde\alpha_i,\qquad
   \sum_{\alpha\in\tau^{-1}(\tilde\alpha)}A_{\alpha}=A_{\tilde\alpha}
$$
for $i=1,\dots,k$ and $\tilde\alpha\in\tilde T$.

(3) The stabilizations $\st(\z^\nu)$ converge to the stabilization
$\st(\z)$ and $\f$ is $J_{\st(\z)}$-holomorphic (cf.~\cite{MS},
Theorem 5.6.6).

(4) After suitable reparametrizations, the $\f^\nu$ converge
uniformly to $\f$.

\begin{remark}
By Theorem~\ref{thm:Gromov-comp} the moduli space $\bar\MM_k(A,J)$ of
stable maps is compact. We wish to emphasize, however, that the
subset $\bar\MM^\ds_k(A,J)\subset\bar\MM_k(A,J)$ of domain-stable maps is
in general not closed in the Gromov topology: The underlying nodal curve
of the Gromov limit of a sequence of domain-stable maps need not be stable.
However, in Section~\ref{sec:proof} we will encounter a
situation in which $\bar\MM^\ds_k(A,J)$ (with additional constraints)
{\em is} closed in $\bar\MM_k(A,J)$.
\end{remark}

{\bf Transversality. } The goal of this section is to show that for
a {\em stable} tree $T$, the space $\MM_T(A,J)$ is a manifold for
generic $J\in\JJ_{k+1}$, as well as various refinements of this
result. For this, we need some more notation. 

Fix a nodal curve $\z=(\{z_{\alpha\beta}\}_{\alpha E\beta},
\{\z_i\}_{1\leq i\leq k})$ with $k$ marked points, modelled over the
$k$-labelled tree $T=(T,E,\Lambda)$. Fix $mp>2$ and homology classes
$A_\alpha$, $\alpha\in T$, and set
$$
   \BB(\{A_\alpha\}) := \BB^{m,p}(\{A_\alpha\}):= \prod_{\alpha\in
   T}\BB^{m,p}(A_\alpha).
$$
We write elements in $\BB(\{A_\alpha\})$ as $\f=(f_\alpha)_{\alpha\in
T}$ with $f_\alpha:S_\alpha\to X$. Let $X^E$ be the product indexed by
pairs $\alpha\beta$ such that $\alpha E\beta$, and the {\em edge
diagonal} $\Delta^E\subset X^E$ be defined by
$x_{\alpha\beta}=x_{\beta\alpha}$ for $\alpha E\beta$. Evaluation at
the $z_{\alpha\beta}$ defines the {\em edge evaluation map}
$$
   \ev^E:\BB(\{A_\alpha\})\to X^E.
$$
Note that $\f\in\BB(\{A_\alpha\})$ induces a continuous map
$\Sigma_\z\to X$ iff $\ev^E(\f)\in\Delta^E$. For
$\J=(J_\alpha)\in\JJ_T$, the Cauchy-Riemann operators
$f_\alpha\mapsto\pb_{J_\alpha}f_\alpha$ define a smooth section
$$
   \pb_\J:\BB(\{A_\alpha\})\to\EE(\{A_\alpha\})
$$
in the bundle
$$
   \EE(\{A_\alpha\}):=\EE^{m-1,p}(\{A_\alpha\}):= \prod_{\alpha\in
   T}\EE^{m-1,p}(A_\alpha).
$$
Note that $(\z,\f)$ is a nodal $\J$-holomorphic map iff
$\f\in(\pb_J,\ev^E)^{-1}(0_\EE\times\Delta^E)\subset
\BB(\{A_\alpha\})$, where $0_\EE$ denotes the zero section of
$\EE(\{A_\alpha\})$. We also have the {\em $k$-point evaluation map}
$$
   \ev^k:\BB(\{A_\alpha\})\to X^k,\qquad
   \f\mapsto\bigl(f_{\alpha_1}(z_1),\dots,f_{\alpha_k}(z_k)\bigr).
$$
Similarly, the reduced index set $R\subset\{1,\dots,k\}$ gives rise to
the {\em reduced evaluation map}
$$
   \ev^R:\BB(\{A_\alpha\})\to X^R,\qquad
   \f\mapsto\{f_{\alpha_i}(z_i)\}_{i\in R}.
$$

We first consider the case $\BB(\{0\})$ in which $A_\alpha=0$ for all
$\alpha\in T$. Then the whole tree $T$ is a ghost tree and therefore
$R=\{1\}$. 

\begin{lemma}\label{lem:surj-const-tree}
Let $T$ be a 1-labelled tree and $\J\in\JJ_T$. Then for any $x\in X$,
$$
   (\pb_\J,\ev^E,\ev^1):\BB(\{0\})\to
   \EE(\{0\})\times X^E\times X
$$
is transverse to $0_\EE\times\Delta^E\times\{x\}$.
\end{lemma}

\begin{proof}
Let $\f\in(\pb_\J,\ev^E,\ev^1)^{-1}(0_\EE\times\Delta^E\times\{x\})
\subset\BB(\{0\}$. Since $\f$ has zero energy, this implies
$f_\alpha\equiv x$ for all $\alpha\in T$. The linearized operator at
$\f$ is given by
$$
   T_\f\BB\to\EE_\f\oplus (T_xX)^E\oplus T_xX,\qquad
   \{\xi_\alpha\}_{\alpha\in T}\mapsto\!\Bigl
   (\{\pb_\alpha\xi_\alpha\}_{\alpha\in
   T},\{\xi_\alpha(z_{\alpha\beta})\}_{\alpha
   E\beta},\xi_{\alpha_1}(z_1)\Bigr),
$$
where $\pb_\alpha:T_{f_\alpha}\BB\to \EE_{f_\alpha}$ is the linearized
Cauchy-Riemann operator. By Lemma \ref{lem:surj-const}, $\pb_\alpha$ is
surjective and its kernel consists of the constant functions
$\xi_\alpha:S_\alpha\to T_xX$. Thus it suffices to show that the
restriction to the kernel of the $\pb_\alpha$
$$
   (T_xX)^T\to (T_xX)^E\oplus T_xX,\qquad \{\xi_\alpha\}_{\alpha\in
   T}\mapsto(\{\xi_\alpha\}_{\alpha E\beta},\xi_{\alpha_1})
$$
is transverse to the subspace $T_{x^E}\Delta^E\oplus\{0\}\subset
(T_xX)^E\oplus T_xX$. This is equivalent to the surjectivity of the
operator
\begin{gather*}
   (T_xX)^T\oplus T_{x^E}\Delta^E \to (T_xX)^E\oplus T_xX, \cr
   (\{\xi_\alpha\}_{\alpha\in T},\{\eta_{\alpha\beta}\}_{\alpha
   E\beta}) \mapsto(\{\xi_\alpha+\eta_{\alpha\beta}\}_{\alpha
   E\beta},\xi_{\alpha_1}).
\end{gather*}
Thus for given $\{v_{\alpha\beta}\}_{\alpha E\beta}$ and $v_1$, we
need to find $\{\xi_\alpha\}_{\alpha\in T}$ and
$\{\eta_{\alpha\beta}\}_{\alpha
   E\beta}$ satisfying
\begin{equation}\label{eq:system}
   \xi_\alpha+\eta_{\alpha\beta}=v_{\alpha\beta},\qquad
   \xi_{\alpha_1}=v_1,\qquad \eta_{\alpha\beta}=\eta_{\beta\alpha}
\end{equation}
for all $\alpha E\beta$. This system can be solved by induction as
follows. Let $T_\ell$ be the set of vertices that have distance at
most $\ell$ from $\alpha_1$. For $\ell=0$ set
$\xi_{\alpha_1}:=v_1$. For the induction step, suppose that we have a
solution of the system~\eqref{eq:system} for all $\alpha E\beta$
with $\alpha,\beta\in T_\ell$. Since $T$ is a tree, for $\beta\in
T_{\ell+1}\setminus T_\ell$ there exists a unique $\alpha\in T_\ell$
with $\alpha E\beta$. Then
$\eta_{\alpha\beta}:=\eta_{\beta\alpha}:=v_{\alpha\beta}-\xi_\alpha$
and $\xi_\beta:=v_{\beta\alpha}-\eta_{\beta\alpha}$
solves~\eqref{eq:system} on $T_{\ell+1}$.
\end{proof}

More generally, let $k\geq 3$ and let $I$ be a subset of
$\{1,\dots,k\}$ with $|I|\geq 3$. Call a 
$k$-labelled tree $T$ {\em $I$-stable} if it is still stable
after removing the marked points not in $I$. Note that $I$-stability
implies stability. Denote by 
$$
   \pi_I:\bar\MM_{k+1}\to\bar\MM_{|I|+1}
$$
the obvious projection (forgetting the marked points outside the
set $I$ and stabilizing). 

Fix $J_0\in\JJ$ and a nonempty open subset $V\subset X$ and set
$$
   \JJ_{k+1}(V) := \{J\in\JJ_{k+1}\mid J=J_0 \text{ outside }V\}.
$$
The projection $\pi_I$ defines a pullback map
$$
   \pi_I^*:\JJ_{|I|+1}(V)\to \JJ_{k+1}(V)
$$
and we will often denote the pullback $\pi_I^*J$ again by $J$. 

Fix a {\em stable} curve $\z$ modelled over the $k$-labelled tree
$T$ and homology classes $A_\alpha$, $\alpha\in T$. We have the {\em
  universal Cauchy-Riemann operator} 
$$
   \pb:\BB(\{A_\alpha\})\times\JJ_{|I|+1}(V)\to\EE(\{A_\alpha\}),\qquad
   (\f,J)\mapsto \pb_{J_\z}\f.
$$
\begin{definition}\label{def:B}
Denote by $\BB^*(\{A_\alpha\})\subset\BB(\{A_\alpha\})$ the open
subset of maps $\f$ for which the image of each nonconstant $f_\alpha$
meets $V$. 
\end{definition}

\begin{lemma}\label{lem:surj-tree}
For any stable curve $\z$ and $x\in X^R$
$$
   (\pb,\ev^E,\ev^R):\BB^*(\{A_\alpha\})\times\JJ_{|I|+1}(V)\to
   \EE(\{A_\alpha\})\times X^E\times X^R
$$
is transverse to $0_\EE\times \Delta^E\times\{x\}$.
\end{lemma}

\begin{proof}
Let $\f\in\BB^*(\{A_\alpha\})$ and $J\in\JJ_{|I|+1}(V)$ satisfy
$\pb_{J_\z}\f=0$. The linearization of $(\pb,\ev^E,\ev^R)$ at
$(f,J)$ is the linear operator
\begin{gather*}
   \bigoplus_{\alpha\in T}T_{f_\alpha}\BB\oplus T_J\JJ_{|I|+1}
   \to\bigoplus_{\alpha\in T}\EE_{f_\alpha}\oplus \bigoplus_{\alpha
   E\beta}T_{f_\alpha(z_{\alpha\beta})}X\oplus
   \bigoplus_{i\in R}T_{f_{\alpha_i}(z_i)}X, \cr
   (\{\xi_\alpha\},Y)\mapsto\Bigl(\{D_{f_\alpha}\xi_\alpha +
   \frac{1}{2}Y_\z\circ df_\alpha\circ j\},
   \{\xi_\alpha(z_{\alpha\beta})\}, \{\xi_{\alpha_i}(z_i)\}\Bigr).
\end{gather*}
Note that this operator only depends on the restrictions of $Y$ to
the $S_\alpha$. Let $T^\nc := \{\alpha\in T\mid f_\alpha\neq$
constant$\}$. Since $\f\in\BB^*(\{A_\alpha\})$, we find nonempty open
sets $U_\alpha\subset f_\alpha^{-1}(V)$ for $\alpha\in T^\nc$ that are
disjoint in $\Sigma_\z$. Consider the space
$$
   \bigoplus_{\alpha\in T}T_{J_\alpha}\JJ_{U_\alpha},
$$
where $J_\alpha=J|_{S_\alpha}$ and for $\alpha\notin T^\nc$ we have set
$U_\alpha:=\emptyset$ and $T_{J_\alpha}\JJ_{U_\alpha}:=\{0\}$.
By Lemma~\ref{lem:J-smooth}, each element
$\{Y_\alpha\}\in\bigoplus_{\alpha\in T}T_{J_\alpha}\JJ_{U_\alpha}$
extends to an element $Y\in T_J\JJ_{|I|+1}$. By the choice of the
$U_\alpha$, the extension can be chosen so that $Y\in
T_J\JJ_{|I|+1}(V)$. Composing with the corresponding injective linear map
(see Lemma~\ref{lem:J-smooth}) it suffices to prove
transversality to $T(0_\EE\times \Delta^E\times\{x\})$ of the
linearized operator
\begin{gather*}
   L:\bigoplus_{\alpha\in T}\Bigl[T_{f_\alpha}\BB\oplus
   T_{J_\alpha}\JJ_{U_\alpha}\Bigr]
   \to\bigoplus_{\alpha\in T}\Bigl[\EE_{f_\alpha}\oplus
   \bigoplus_{\beta:\alpha
   E\beta}T_{f_\alpha(z_{\alpha\beta})}X\oplus
   \bigoplus_{i:\alpha_i=\alpha}T_{f_\alpha(z_i)}X\Bigr], \cr
   \Bigl\{(\xi_\alpha,Y_\alpha)\Bigr\}
   \mapsto\Bigl\{\Bigl(D_{f_\alpha}\xi_\alpha +
   \frac{1}{2}Y_\alpha\circ df_\alpha\circ j,
   \{\xi_\alpha(z_{\alpha\beta})\}, \{\xi_\alpha(z_i)\}\Bigr)\Bigr\}.
\end{gather*}
Note that this operator is a direct sum of operators $L_\alpha$
(for $\alpha\in T$) of the form considered in
Section~\ref{sec:holo}. For $\alpha\in T^\nc$, $L_\alpha$ is
surjective by Lemma~\ref{lem:surj-nonconst}. (Note that for
$\alpha\notin T^\nc$ Lemma~\ref{lem:surj-nonconst} is not applicable
because $U_\alpha=\emptyset$.)

Next consider a maximal subtree $T'=(T',E',\Lambda')$
contained in $T\setminus T^\nc$, i.e.~$A_\alpha=0$ for all
$\alpha\in T'$. Thus $T'$ is a ghost tree, so by definition of $R$ it
contains at most one marked point in $R$. If a point $z_i\in R$ on
$T'$ exists denote by $\ev^1_i$ the evaluation map at this point 
and set $x_i:=\f(z_i)$. By Lemma~\ref{lem:surj-const-tree}, the map
$(\pb,\ev^{E'},\ev^1_i)$ is transverse to
$0_\EE\times\Delta^{E'}\times\{x_i\}$ over $T'$. 
Combining this with the surjectivity of $L_\alpha$ for $\alpha\in
T^\nc$, the transversality of $L$ to $T(0_\EE\times
\Delta^E\times\{x\})$ follows.
\end{proof}

For $J\in\JJ_{|I|+1}(V)$, denote by
$$
   \tilde\MM_T^*(\{A_\alpha\},J) \subset \tilde\MM_T(\{A_\alpha\},J)
$$
the open subset of those $[\z,\f]$ for which
$\f\in\BB^*(\{A_\alpha\})$.

\begin{prop}\label{prop:surj-tree}
Let $T$ be an $I$-stable $k$-labelled tree, $A_\alpha\in H_2(X;\Z)$
for $\alpha\in T$, and $J_0\in\JJ$. Then the {\em universal solution
  space} 
$$
   \tilde\MM_T^*\bigl(\{A_\alpha\},\JJ_{|I|+1}(V)\bigr) :=
   \bigcup_{J\in\JJ_{|I|+1}(V)}\tilde\MM_T^*(\{A_\alpha\},J)
$$
is a Banach manifold and the reduced evaluation map
$$
   \ev^R:\tilde\MM_T^*\bigl(\{A_\alpha\},\JJ_{|I|+1}(V)\bigr)\to X^R
$$
is a submersion.
\end{prop}

\begin{proof}
By Lemma~\ref{lem:surj-tree}, for every $x\in X^R$ the linearization
of the map 
\begin{gather*}
   \tilde\MM_T\times\BB^*(\{A_\alpha\})\times\JJ_{|I|+1}(V)\to
   \EE(\{A_\alpha\})\times X^E\times X^R, \cr
   (\z,\f,J)\mapsto\bigl(\pb_{J_\z}\f,\ev^E(\f),\ev^R(\f)\bigr).
\end{gather*}
is transverse to $T(0_\EE\times \Delta^E\times\{x\})$. Thus the
preimage
$$
   (\pb,\ev^E)^{-1}(0_\EE\times\Delta^E) =
   \tilde\MM^*_T\bigl(\{A_\alpha\},\JJ_{|I|+1}(V)\bigr)
$$
is a Banach manifold and the reduced evaluation map
$$
   \ev^R:\tilde\MM_T^*\bigl(\{A_\alpha\},\JJ_{|I|+1}(V)\bigr)\to X^R
$$
is a submersion.
\end{proof}

In view of Lemma~\ref{lem:action} and~\cite{MS}, the group $G_T$
of isomorphisms fixing the stable tree $T$ acts freely and
properly on the space
$\tilde\MM^*_T\bigl(\{A_\alpha\},\JJ_{|I|+1}(V)\bigr)$. Hence the
{\em universal moduli space}
$$
   \MM^*_T\bigl(\{A_\alpha\},\JJ_{|I|+1}(V)\bigr) :=
   \tilde\MM^*_T\bigl(\{A_\alpha\},\JJ_{|I|+1}(V)\bigr)/G_T
$$
is a Banach manifold. Moreover, the reduced evaluation map
descends to a submersion
$$
   \ev^R:\MM_T^*\bigl(\{A_\alpha\},\JJ_{|I|+1}(V)\bigr)\to X^R.
$$
For $J\in\JJ_{|I|+1}(V)$ define the moduli space
$$
   \MM_T^*(\{A_\alpha\},J) := \tilde\MM_T^*(\{A_\alpha\},J)/G_T
$$
of $J$-holomorphic maps belonging to the set
$\BB^*(\{A_\alpha\})$. For a smooth submanifold $Z\subset X^R$ define
$$
   \MM_T^*(\{A_\alpha\},J;Z) :=
   \{[\z,\f]\in\MM^*_T(\{A_\alpha\},J)\mid
   \{f_{\alpha_i}(z_i)\}_{i\in R}\in Z\}.
$$
\begin{cor}\label{cor:surj-tree}
Let $T$ be an $I$-stable $k$-labelled tree, $A_\alpha\in H_2(X;\Z)$
for $\alpha\in T$, $J_0\in\JJ$, and $Z\subset X^R$ a smooth
submanifold. Then there exists a Baire set
$\JJ_{|I|+1}^\reg(V)\subset\JJ_{|I|+1}(V)$ (depending on $Z$, $T$ and
the $A_\alpha$) such that for all $J\in\JJ_{|I|+1}^\reg(V)$ the moduli
space $\MM^*_T(\{A_\alpha\},J;Z)$ is a smooth submanifold of
$\MM^*_T(\{A_\alpha\},J)$ of dimension
$$
   \dim\MM^*_T(\{A_\alpha\},J;Z) = 2n -6 +
   \sum_{\alpha\in T}2c_1(A_\alpha) + 2k - 2e(T) - \codim_\R Z.
$$
\end{cor}

\begin{proof}
By the preceding discussion, the universal moduli space
$$
   \MM^*_T\bigl(\{A_\alpha\},\JJ_{|I|+1}(V);Z\bigr) :=
   (\ev^R)^{-1}(Z)\subset \MM_T^*\bigl(\{A_\alpha\},\JJ_{|I|+1}(V)\bigr)
$$
is a smooth Banach submanifold of
$\MM_T^*\bigl(\{A_\alpha\},\JJ_{|I|+1}(V)\bigr)$ of codimension
$\codim_\R Z$ with a smooth projection
$$
   \pi:\MM_T^*\bigl(\{A_\alpha\},\JJ_{|I|+1}(V);Z\bigr)
   \to\JJ_{|I|+1}(V),\qquad (f,J)\mapsto J.
$$
By the Sard-Smale theorem, there exists a Baire set
$\JJ_{|I|+1}^\reg(V)\subset\JJ_{|I|+1}(V)$ of regular values of the
projection $\pi$.
It follows that for $J\in \JJ_{|I|+1}^\reg(V)$ the moduli space
$\MM^*_T(\{A_\alpha\},J)$ is a smooth manifold and
$\MM^*_T(\{A_\alpha\},J;Z)$ is a smooth submanifold of codimension
$\codim_\R Z$ in $\MM^*_T(\{A_\alpha\},J)$. Since, according to
Lemma~\ref{lem:surj-tree}, $\MM^*_T(\{A_\alpha\},J)$ is
transversely cut out by the edge evaluation map, its dimension
follows from the usual dimension formula for its components and
$|T|=e(T)+1$:
\begin{align*}
   \dim\MM^*_T(\{A_\alpha\},J) &= \sum_{\alpha\in
   T}\Bigl(2n-6+2c_1(A_\alpha)\Bigr) + 2k + 4e(T) -2ne(T) \cr
   &= 2n -6|T| + 2\sum_{\alpha\in T}c_1(A_\alpha) + 2k + 4e(T) \cr
   &= 2n -6 + \sum_{\alpha\in T}2c_1(A_\alpha) + 2k - 2e(T).
\end{align*}
This implies the dimension formula in Corollary~\ref{cor:surj-tree}.
\end{proof}

We conclude with an improvement of the dimension formula in
Corollary~\ref{cor:surj-tree} that will play a crucial role
Section~\ref{sec:proof}. If $|R|<k$ we have marked points on ghost
components not belonging to $R$. Variations of these points contribute
to the dimension formula but have no effect on the $J$-holomorphic
maps (due to the coherency of $J$). So the ``actual'' dimension of the
moduli space should be lower than the formula in
Corollary~\ref{cor:surj-tree}. 
This statement can be made rigorous as follows: 
Suppose that $Z=\prod_{i\in R}Z_i$ and $V\cap Z_i=\emptyset$ for all
$i\in R$. 
Evaluation at the points of $R$ induces a smooth evaluation map
$$
   \ev^R:\MM^*_T(\{A_\alpha\},J;Z)\to Z. 
$$
Let $\pi_R(T)$ be the stable $|R|$-labelled tree obtained from $T$ by
forgetting the marked points not belonging to $R$ and stabilizing. On
non-ghost components this does not affect $J$ due to the coherency
condition. 
On ghost components $J\equiv J_0$ because each ghost tree is mapped to
on the $Z_i$ $\alpha$ and $V\cap Z_i=\emptyset$. 
Hence forgetting the marked points outside $R$ yields a
transversely cut out moduli space $\MM^*_{\pi_R(T)}(\{A_\alpha\},J;Z)$.
Since the evaluation map descends to this space, we have shown 

\begin{cor}\label{cor:surj-tree2}
In the situation of Corollary~\ref{cor:surj-tree}, suppose that
$Z=\prod_{i\in R}Z_i$ and $V\cap Z_i=\emptyset$ for all $i\in R$. Then
the evaluation map $\ev^R$ factors as 
$$
   \ev^R:\MM^*_T(\{A_\alpha\},J;Z)\to
   \MM^*_{\pi_R(T)}(\{A_\alpha\},J;Z)\to Z 
$$
through a smooth manifold of dimension 
$$
   \dim\MM^*_{\pi_R(T)}(\{A_\alpha\},J;Z) = 2n -6 +
   \sum_{\alpha\in T}2c_1(A_\alpha) + 2|R| - 2e\bigl(\pi_R(T)\bigr) -
   \codim_\R Z.  
$$
\end{cor}

\begin{remark}\label{rem:add}
The proofs show that Proposition~\ref{prop:surj-tree},
Corollary~\ref{cor:surj-tree} and Corollary~\ref{cor:surj-tree2}
remain true if $\JJ_{|I|+1}(V)$ is replaced by the subset
$\JJ_{|I|+1}^*(V)\subset\JJ_{|I|+1}(V)$ of those coherent $J$ mapping
into some given open subset $\JJ^*\subset\JJ$.
\end{remark}

\section{Tangencies}\label{sec:tangencies}

In this section we define moduli spaces of $J$-holomorphic spheres
with prescribed orders of tangency to complex submanifolds and prove
their regularity.

We begin with the local situation. Suppose $A:\R^m\to\R^{m\times
m}$ is a smooth matrix valued function. Assume further that $A$
preserves a subspace $\R^k:=\R^k\times\{0\}\subset\R^m$ ($0\leq
k\leq m-1$) in the sense that $A(x)(\R^k\times\{0\})\subset
\R^k\times\{0\}$ for all $x\in\R^k$. Let $D\subset\C$ be the
unit disk and $f:D\to\R^m$ be a smooth map. Denote coordinates on
$\C$ by $z=s+it$ and partial derivatives by $f_s$, $f_{st}$ etc.
Denote by $d^\ell f$ the vector of all partial derivatives of
order $i$ with $1\leq i\leq\ell$. We say $d^\ell f\in\R^k$ if all
the partial derivatives in $d^\ell f$ lie in the subspace $\R^k$.
The following lemma is proved by a simple induction over $i+j$.

\begin{lemma}\label{lem:A}
For all integers $i,j\geq 0$ we have an equation
$$
   \frac{\p^{i+j}}{\p s^i\p t^j}\Bigl(A(f)f_s\Bigr) =
   A(f)\frac{\p^{i+j+1}f}{\p s^{i+1}\p t^j}+A_{i,j}(f)(d^{i+j}f),
$$
where $A_{i,j}(f)$ is a multilinear function in the partial
derivatives up to order $i+j$. Moreover,
$A_{i,j}(f)(d^{i+j}f)\in\R^k$ if $f\in\R^k$ and $d^{i+j}f\in\R^k$.
\end{lemma}

We first apply Lemma~\ref{lem:A} in the following situation. Let $J$
be an almost complex structure on $\C^n$. We say that {\em $J$
preserves $\C^k:=\C^k\times\{0\}\subset\C^n$} if $J(x):\C^k\to\C^k$
for all $x\in\C^k$. Let $f:\C\supset D\to\C^n$ be a $J$-holomorphic
map, i.e., $f_s+J(f)f_t=0$. Then $f_t=J(f)f_s$ and Lemma~\ref{lem:A}
(with $A=J$) yields for $i,j\geq 0$:
$$
   \frac{\p^{i+j+1}f}{\p s^i\p t^{j+1}} = J(f)\frac{\p^{i+j+1}f}{\p
   s^{i+1}\p t^j} + A_{i,j}(f)(d^{i+j}f).
$$
Multiplying by $J$ we find
$$
   \bigl(\p_s+J(f)\p_t\bigr)\frac{\p^{i+j}f}{\p s^i\p t^j} =
   J_{i,j}(f)(d^{i+j}f),
$$
where $J_{i,j}(f):=J(f)A_{i,j}(f)$ satisfies
$J_{i,j}(f)(d^{i+j}f)\in\C^k$ if $f\in\C^k$ and $d^{i+j}f\in\C^k$.
Thus partial derivatives of $f$ are $J$-holomorphic up to lower
order terms. This equation allows us to swap one $s$-derivative for
a $t$-derivative, and iterating this process we obtain

\begin{cor}\label{cor:lin-tangent}
Let $f:(D,0)\to(\C^n,0)$ be $J$-holomorphic with $J$ preserving
$\C^k$. If $d^{\ell-1}f(0)\in\C^k$ and $\frac{\p^\ell f}{\p
s^\ell}(0)\in\C^k$ for some $\ell\geq 1$, then $d^\ell f(0)\in\C^k$.
\end{cor}

If $f(0)=0$ and $d^{\ell-1}f(0)\in\C^k$ we define the equivalence
class
$$
   j^\ell_{\C^k}f(0) := \left[\frac{\p^\ell f}{\p
   s^\ell}(0)\right]\in\C^n/\C^k.
$$
By Corollary~\ref{cor:lin-tangent}, $j^\ell_{\C^k}f(0)=0$ implies $d^\ell
f(0)\in\C^k$. So $j^{\ell+1}_{\C^k}f(0)$ is defined whenever $f(0)=0$ and
$j^1_{\C^k}f(0)=\dots=j^\ell_{\C^k}f(0)=0$.

Next we study how $j^\ell_{\C^k}f(0)$ transforms under diffeomorphisms
$\phi:(\C^n,0)\to(\C^n,0)$ with $\phi(\C^k)\subset\C^k$. We have
$(\phi\circ f)_s=D\phi(f)f_s$, so Lemma~\ref{lem:A} (with $A=D\phi$)
yields
$$
   \frac{\p^{i+j+1}(\phi\circ f)}{\p s^{i+1}\p t^j} =
   D\phi(f)\frac{\p^{i+j+1}f}{\p
   s^{i+1}\p t^j} + A_{i,j}(f)(d^{i+j}f).
$$
Of course, the same argument applies to the $t$-derivative and so we
find for all $i+j\geq 1$:
$$
   \frac{\p^{i+j}(\phi\circ f)}{\p s^i\p t^j} =
   D\phi(f)\frac{\p^{i+j}f}{\p
   s^i\p t^j} + \phi_{i,j}(f)(d^{i+j-1}f),
$$
where $\phi_{i,j}(f)$ satisfies $\phi_{i,j}(f)(d^{i+j-1}f)\in\C^k$
if $f\in\C^k$ and $d^{i+j-1}f\in\C^k$. This proves

\begin{cor}\label{cor:lin-transform}
Let $f:(D,0)\to(\C^n,0)$ be a smooth map and
$\phi:(\C^n,0)\to(\C^n,0)$ a diffeomorphism preserving $\C^k$. If
$d^{\ell-1}f(0)\in\C^k$ for some $\ell\geq 1$, then
$$
   \left[\frac{\p^\ell(\phi\circ f)}{\p s^\ell}(0)\right] =
   \left[D\phi(0)\frac{\p^\ell
   f}{\p s^\ell}(0)\right] \in\C^n/\C^k.
$$
\end{cor}

Corollary~\ref{cor:lin-transform} shows that if $f:(D,0)\to(\C^n,0)$
is $J$-holomorphic and $j^1_{\C^k}f(0)=\dots=j^{\ell-1}_{\C^k}f(0)=0$,
then $j^\ell_{\C^k}f(0)$ transforms like a normal vector to the
submanifold $\C^k\subset\C^n$.

Now it is clear how to globalize these constructions. Let $Z$ be a
$k$-dimensional (almost) complex submanifold of an almost complex
manifold
$(X,J)$. Let $f:\C\supset D\to X$ be a $J$-holomorphic map with
$f(0)=x\in Z$. Pick a coordinate chart $\psi$ for $X$ mapping $x$ to
$0$ and $Z$ to $\C^k$. We say that {\em $f$ is tangent to order
$\ell\geq 0$ to $Z$ at $0$}, and write $d^\ell f(0)\in T_xZ$, if
$d^\ell(\psi\circ f)(0)\in\C^k$. By the discussion preceding
Corollary~\ref{cor:lin-transform} this definition is independent of
the coordinate chart $\psi$.

Suppose that for some $\ell\geq 1$ we have $d^{\ell-1}f(0)\in
T_xZ$. Then we define the {\em $\ell$-jet normal to $Z$ of $f$ at
$0$} as
$$
   j^\ell_Z f(0) := D\psi(x)^{-1}j^\ell(\psi\circ f)(0) \in T_xX/T_xZ.
$$
By Corollary~\ref{cor:lin-tangent}, $j^\ell_Z f(0)$ is defined
whenever $j^1_Zf(0)= \dots =j^{\ell-1}_Zf(0)=0$ and by
Corollary~\ref{cor:lin-transform}, the definition of $j^\ell_Zf(0)$
does not depend on the coordinate chart $\psi$. So we have proved

\begin{lemma}\label{lem:jets-loc}
Let $Z$ be a complex submanifold of an almost complex manifold
$(X,J)$ and $f:\C\supset D\to X$ be a $J$-holomorphic map with
$f(0)=x\in Z$. Then there is a sequence
$j^1_Zf(0),j^2_Zf(0),\dots\in T_xX/T_xZ$ of normal derivatives such
that the $\ell$-th normal derivative $j^\ell_Z f(0)$ is defined
whenever $j^1_Zf(0)=\dots=j^{\ell-1}_Zf(0)=0$. Moreover, $f$ is
tangent to order $\ell$ to $Z$ at $0$ if and only if
$j^1_Zf(0)=\dots=j^\ell_Zf(0)=0$.
\end{lemma}

Now let $(X,\om)$ be a symplectic manifold with a tamed almost
complex structure $J_0$ and let $Z_1,\dots,Z_k$ be $k$ complex
submanifolds of $(X,J_0)$. Fix a nonempty open subset $V\subset X$
with $V\cap(Z_1\cup\dots\cup Z_k)=\emptyset$ and denote by
$$
   \JJ(V) := \{J\in\JJ(X,\om)\mid J=J_0 \text{ on }X\setminus V\}
$$
the space of tamed almost complex structures that agree with
$J_0$ outside $V$. Note that $Z_1,\dots,Z_k$ are $J$-complex
submanifolds for all $J\in\JJ(V)$.

Fix nonnegative integers $\ell_1,\dots,\ell_k$ and pairwise
  disjoint points $z_1,\dots,z_k\in S^2$. Denote by $\z:=(z_1,...,z_k)$
their collection
Also choose holomorphic coordinate
charts for $S^2$ around the $z_i$ mapping $z_i$ to $0\in\C$. For
$m-2/p>\max\{\ell_i\}$ define the {\em universal solution space}
\begin{align*}
   \tilde\MM_\z:=\tilde\MM_\z\bigl(\JJ(V);\{Z_i,\ell_i\}\bigr) &:=
   \{(f,J)\in\BB^{m,p}\times\JJ(V)\mid \pb_Jf=0,\ \cr
   &\ \ \ d^{\ell_i}f(z_i)\in T_{f(z_i)}Z_i \text{ for }i=1,\dots,k\}
\end{align*}
of holomorphic maps tangent to order $\ell_i$ to $Z_i$ at
$z_i$. Notice that this 
definition does not depend on the choice of holomorphic coordinates
around the $z_i$. 
Denote by
$$
   \tilde\MM_\z^*\bigl(\JJ(V);\{Z_i,\ell_i\}\bigr) :=
   \{(f,J)\in\tilde\MM_\z\bigl(\JJ(V);\{Z_i,\ell_i\}\bigr)
   \mid f^{-1}(V)\neq\emptyset\}
$$
the subspace of solutions passing through $V$, and by
$$
   \tilde\MM_\z^s:=\tilde\MM_\z^s\bigl(\JJ(V);\{Z_i,\ell_i\}\bigr) :=
   \{(f,J)\in\tilde\MM_\z^*\bigl(\JJ(V);\{Z_i,\ell_i\}\bigr)
   \mid f \text{ simple }\}
$$
the subspace of simple solutions passing through $V$. Here a
holomorphic map is called {\em simple} if it is not a nontrivial
branched covering of another holomorphic map. We wish to show that
$\tilde\MM_\z^s$ is a Banach manifold. So consider $(f,J)\in
\tilde\MM_\z^s$. Pick coordinate charts $\phi_i$ for $X$ around
$f(z_i)$ mapping $f(z_i)$ to $0$ and $Z_i$ to $\C^{n_i}\subset\C^n$.
Using $\phi_i$, we view $f(z)$ for $z$ near $z_i$ as a point in
$\C^n$. Then the candidate for the tangent space of $\tilde\MM_\z^s$ at
$(f,J)$ is
\begin{align}\label{eq:tangent-space}
   T_{(f,J)}\tilde\MM_\z^s &:= \{(\xi,Y)\in T_f\BB^{m,p}\times
   T_J\JJ(V)\mid D_f\xi+\frac{1}{2}Y(f)\circ df\circ j=0,\ \cr
   &\ \ \ d^{\ell_i}\xi(z_i)\in\C^{n_i} \text{ for }i=1,\dots,k\},
\end{align}
where $\xi$ is viewed as a $\C^n$-valued function near $z_i$ via
$d\phi_i$.
If in addition $m-2/p>\ell_j+1$ for some $j$, then in view of
Lemma~\ref{lem:jets-loc} (with respect to the fixed coordinate chart for
$S^2$ near $z_j$) we have a {\em normal $(\ell_j+1)$-jet evaluation
map} to the normal bundle $TX/TZ_j$ of $Z_j$ in $X$,
$$
   j^{\ell_j+1}_{Z_j}\ev_j:
   \tilde\MM_\z^s\to TX/TZ,
   \qquad (f,J)\mapsto j^{\ell_j+1}_{Z_j}f(z_j).
$$
With respect to the trivializations above, the linearization of
$j^{\ell_j+1}_{Z_j}\ev_j$ at $(f,J)$ is the linear operator
\begin{equation}\label{eq:Lj}
   L_j:T_{(f,J)}\tilde\MM_\z^s \to \C^n/\C^{k_j},\qquad (\xi,Y)\mapsto
   j^{\ell_j+1}_{\C^{n_j}}\xi(z_j).
\end{equation}

\begin{lemma}\label{lem:jets-reg}
(a) For $m-2/p>\max\{\ell_i\}$ the space
$\tilde\MM_\z^s\bigl(\JJ(V);\{Z_i,\ell_i\}\bigr)$ is a smooth Banach
manifold whose tangent space at $(f,J)$ is given
by~\eqref{eq:tangent-space}.

(b) Suppose that in addition $m-2/p>\ell_j+1$ for some $j$ with
$\ell_j\geq\max\{\ell_i\}-1$. Then the linearized normal
$(\ell_i+1)$-jet evaluation map $L_j$ given by~\eqref{eq:Lj} is
surjective.
\end{lemma}

We will give the proof Lemma~\ref{lem:jets-reg} below.
As the proof of Lemma~\ref{lem:surj-nonconst}, it
is an adaptation of the proof of Proposition 3.4.2 in~\cite{MS}. Note
that besides the presence of higher order jets,
Lemma~\ref{lem:jets-reg} also differs from
Lemma~\ref{lem:surj-nonconst} in that it uses only domain-independent
almost complex structures (which forces us to restrict to simple
holomorphic maps).

The proof is based on the following lemma.
Here we view $\eta\in\EE^{m-1,p}_f$ near $z_i$ via the
trivialization near $f(z_i)$ as a map $\C\supset U_i\to \C^n$, the
corresponding $(0,1)$-form being given by $\eta\,ds-J(f)\eta\,dt$
(cf.~\cite{MS}, Section 3.1). This allows us to speak of the
$\ell$-jets $d^\ell\eta(z_i)$ at the $z_i$. Note that the spaces in
Lemma~\ref{lem:Laszlo} are well-defined by the Sobolev embedding
theorem. We point out that in Lemma~\ref{lem:Laszlo} we consider the
full $\ell$-jets (not just those normal to the $Z_i$) and the integer
$\ell$ is entirely unrelated to the $\ell_i$.

\begin{lemma}\label{lem:Laszlo}
Let $(f,J)\in\tilde\MM_\z^s\bigl(\JJ(V);\{Z_i,\ell_i\}\bigr)$. For
$\ell\geq 0$ and $m-2/p>\ell+1$ define the spaces 
\begin{gather*}
   B^{m,p}_0 := \{\xi\in T_f\BB^{m,p} \mid d^{\ell+1}\xi(z_i)=0
   \text{ for }i=1,\dots,k\}, \cr
   E^{m-1,p}_0 := \{\eta\in \EE_f^{m-1,p} \mid d^\ell\eta(z_i)=0
   \text{ for }i=1,\dots,k\}.
\end{gather*}
Then the linear operator
$$
   F_0:B^{m,p}_0\oplus T_J\JJ(V)\to E^{m-1,p}_0,\qquad
   (\xi,Y)\mapsto D_f\xi+\frac{1}{2}Y(f)\circ df\circ j
$$
is surjective.
\end{lemma}

\begin{proof}
The corresponding statement for $\ell=-1$ was shown in the proof of
Lemma 3.4.3 in~\cite{MS}, see also the proof of
Lemma~\ref{lem:surj-nonconst} above. However, unlike in the case
$\ell=-1$, for $\ell\geq 0$ we cannot reduce the result to the case
$m=1$ because the $(\ell+1)$-jet $d^{\ell+1}\xi(z_i)$ at a point is
not well-defined in the Sobolev space $W^{1,p}$. This forces us to
work with distributions instead of functions.

Recall that $T_f\BB^{m,p}=W^{m,p}(S^2,B)$ and
$\EE_f^{m-1,p}=W^{m-1,p}(S^2,E)$ are spaces of Sobolev sections in
the vector bundles $B:=f^*TX\to S^2$ and $E\!:=\!\Lambda^{0,1}(f^*TX)\!\to
S^2$.

Suppose first that $\ell+1<m-2/p<\ell+2$. Let us view $F_0$ as an
operator from $B^{m,p}_0$ to $\EE_f^{m-1,p}$. Since $D_f$ is a
Fredholm operator, the image $\im\,F_0$ of $F_0$ is closed in
$\EE_f^{m-1,p}$. We will prove the inclusion
$(\im\,F_0)^\perp\subset(E^{m-1,p}_0)^\perp$ in the dual space
$(\EE_f^{m-1,p})^*$. In view of reflexivity of the Sobolev space
$\EE_f^{m-1,p}$ (see e.g.~\cite{Ad}), closedness of $\im\,F_0$, and
the Hahn-Banach theorem, this implies the desired statement
$E^{m-1,p}_0\subset\im\,F_0$. Note that for $m>1$,
$(\EE_f^{m-1,p})^*$ is the space of distributions $\Lambda\in
C^\infty(S^2,E)^*$ that satisfy an inequality
$$
   |\Lambda(\phi)|\leq C\|\phi\|_{W^{m-1,p}}
$$
for all smooth sections $\phi:S^2\to E$. Suppose that $\Lambda$
vanishes on $\im(F_0)$, i.e.,
$$
   \Lambda\Bigl(D_f\xi+\frac{1}{2}Y(f)\circ df\circ j\Bigr) = 0
$$
for all $(\xi,Y)\in B^{m,p}_0$. In particular, this implies that
$\Lambda(D_f\xi)=0$ for all smooth sections $\xi:S^2\to B$ with
compact support in $S^*:=S^2\setminus\{z_1,\dots,z_k\}$. Elliptic
regularity for distributions (see~\cite{Ru}, Theorem 8.12 for the
case with constant coefficients and \cite{DunSch}, Chapter XIV,
Theorem 6.2., p. 1704 for the general case)
implies that the restriction of $\Lambda$ to $S^*$ is represented by
a smooth section $\eta:S^*\to E$ such that
$$
   \Lambda(\phi) = \la\phi,\eta\ra_{L^2}
$$
for all $\phi\in C^\infty_0(S^*,B)$. Moreover, $\eta$ satisfies the
Cauchy-Riemann type equation $D_f^*\eta=0$ and the equation
\begin{equation}\label{eq:Yf}
   \la Y(f)\circ df\circ j,\eta\ra = 0
\end{equation}
for all $Y\in T_J\JJ(V)$. (Note that $Y(f)$ has compact support in
$S^*$ because $V$ is disjoint from the $Z_i$).

By definition of $\tilde\MM_\z^s$ ($f$ is simple!) and Proposition
2.5.1 in~\cite{MS}, the set $U:=f^{-1}(V)$ is nonempty and the set
$U^*:=\{z\in U\mid df(z)\neq 0,\ f^{-1}\bigl(f(z)\bigr)=\{z\}\}$ of
injective points is open and dense in $U$. Suppose that $\eta(z)\neq
0$ for some $z\in U^*$. Then as in the proof of Lemma 3.4.3
of~\cite{MS}, we find a $Y\in T_J\JJ(V)$ such that $\la Y(f)\circ
df(z')\circ j,\eta(z')\ra\geq 0$ for all $z'$ and $>0$ at $z'=z$, in
contradiction to equation~\eqref{eq:Yf}. It follows that $\eta(z)=0$
for all $z\in U^*$, hence by unique continuation (see Lemma 3.4.7
of~\cite{MS}) $\eta\equiv 0$ on $S^*$.

This proves that the distribution $\Lambda$ has support in the
finite set $\{z_1,\dots,z_k\}$. By Theorem 6.25 in~\cite{Ru}, there
exist $\R$-linear functionals $c^i_\alpha:E_{z_i}\cong\C^n\to\R$
such that
\begin{equation}\label{eq:Lambda}
   \Lambda = \sum_{i=1}^k\sum_{|\alpha|\leq N}c^i_\alpha
   D^\alpha\delta_{z_i}.
\end{equation}
Here $\alpha=(\alpha_1,\alpha_2)\in\N_0^2$ are multi-indices,
$D^\alpha = \frac{\p^{\alpha_1+\alpha_2}}{\p s^{\alpha_1}\p
t^{\alpha_2}}$ are partial derivatives, $\delta_{z_i}$ is the delta
distribution at the point $z_i$, and
$D^\alpha\delta_{z_i}(\phi)=(-1)^{|\alpha|}D^\alpha\phi(z_i)$ for
$\phi\in C^\infty(S^2,E)$. The integer $N$ in~\eqref{eq:Lambda} is
the order of $\Lambda$, which is by definition the smallest integer
such that $|\Lambda(\phi)|\leq C\|\phi\|_{C^k}$ for all $\phi\in
C^\infty(S^2,E)$. For $\Lambda\in(\EE_f^{m-1,p})^*$ we clearly have
$N\leq m-1$. More precisely, $\Lambda\in(\EE_f^{m-1,p})^*$ implies
that for every $c^i_\alpha\neq 0$ in~\eqref{eq:Lambda} we must have
$$
   |D^\alpha\phi(z_i)|\leq C\|\phi\|_{W^{m-1,p}}
$$
for all $\phi\in C^\infty(S^2,E)$. It follows from sharpness of the
Sobolev embedding theorem that $|\alpha|\leq[m-1-2/p]$ (if
$|\alpha|>[m-1-2/p]$ there exist smooth functions $\phi$ with
$W^{m-1,p}$-norm $1$ and arbitrarily large $|D^\alpha\phi(z_i)|$).
The assumption $\ell+1<m-2/p<\ell+2$ made above implies
$[m-1-2/p]=\ell$, so we have shown
$$
   \Lambda = \sum_{i=1}^k\sum_{|\alpha|\leq\ell}c^i_\alpha
   D^\alpha\delta_{z_i}.
$$
But this means that
$\Lambda(\phi)=\sum_{i=1}^k\sum_{|\alpha|\leq\ell}c^i_\alpha
(-1)^{|\alpha|}D^\alpha\phi(z_i)=0$ for all $\phi\in E^{m-1,p}_0$,
so $\Lambda\in(E^{m-1,p}_0)^\perp$. This proves
$(\im\,F_0)^\perp\subset(E^{m-1,p}_0)^\perp$ and hence, in view of
the discussion above, the lemma in the special case
$\ell+1<m-2/p<\ell+2$.

The case of general $m,p$ with $m-2/p>\ell+1$ follows from the special
case by elliptic regularity: Pick $m'\leq m$ and $p'\leq p$ such
that $\ell+1<m'-2/p'<\ell+2$. Given $\eta\in E^{m-1,p}_0\subset
E^{m'-1,p'}_0$, by the special case we find $(\xi,Y)\in
B^{m',p'}_0\oplus T_J\JJ(V)$ with $F_0(\xi,Y)=\eta$, and elliptic
regularity (Theorem C.2.3 of~\cite{MS}) yields $\xi\in B^{m,p}_0$.
This concludes the proof of Lemma~\ref{lem:Laszlo}.
\end{proof}

\begin{proof}[Proof of Lemma~\ref{lem:jets-reg}]
Part (a) in Lemma~\ref{lem:jets-reg} in the case
$\ell_1=\dots=\ell_k=0$ follows from Proposition 3.4.2 in~\cite{MS}.
Moreover, in view of Lemma~\ref{lem:jets-loc}, Parts (a) and (b) of
Lemma~\ref{lem:jets-reg} for $\ell_1,\dots,\ell_k$ imply Part (a)
with $\ell_j$ replaced by $\ell_j+1$. Hence by induction (always
preserving the condition $\ell_j\geq\max\{\ell_i\}-1$), it suffices
to prove Part (b).

As above (cf.~\cite{MS}, Section 3.1), we identify $(0,1)$-forms
near the $f(z_i)$ with maps $\C\supset U_i\to\C^n$. With this
identification understood, the linearized Cauchy-Riemann operator at
$f$ near $z_i$ is given in Section 3.1 of~\cite{MS} as
\begin{align*}
   D_f\xi = \xi_s + J(z)\xi_t + A(z)\xi
\end{align*}
with the $z$-dependent linear maps $J(z)\!:=\!J\bigl(f(z)\bigr)$ and
$A(z):= DJ\bigl(f(z)\bigr)\bigl(\cdot,f_t(z)\bigr)$. Pick the
coordinate charts such that $z_j$ corresponds to $0\in\C$ and
$J(0)=i$ is the standard complex structure on $\C^n$.

Now fix $j$ as in Part (b) and let $v\in\C^n$ be given.
Define the holomorphic function 
$$
   \tilde\xi(z) := z^{\ell_j+1}v
$$
near $z=0$ in the chart around $z_j$. This function satisfies
$\tilde\xi_s+i\tilde\xi_t=0$ and 
$$
   d^{\ell_j}\tilde\xi(0)=0,\qquad j^{\ell_j+1}\tilde\xi(0) =
   \frac{\p^{\ell_j+1}\tilde\xi}{\p
   s^{\ell_j+1}}(0) = \frac{\p^{\ell_j+1}\tilde\xi}{\p z^{\ell_j+1}}(0) =
   v.
$$
The function $\tilde\eta := D_f\tilde\xi$ is given by
$$
   \tilde\eta = \tilde\xi_s+i\tilde\xi_t + \bigl(J(z)-i\bigr)\tilde\xi_t + A(z)\tilde\xi =
   B(z)\tilde\xi_t + A(z)\tilde\xi,
$$
where the linear map $B(z) := J(z)-i$ satisfies $B(0)=0$. Since
$d^{\ell_j}\tilde\xi(0)=0$, all terms in $d^{\ell_j}\tilde\eta(0)$ vanish
except possibly $B(0)j^{\ell_j+1}\tilde\xi(0)$, but this one vanishes as
well because $B(0)=0$. So we have $d^{\ell_j}\tilde\eta(z_j)=0$.

Now extend $\tilde\xi$ to $S^2$ such that it vanishes near the
$z_i$, $i\neq j$, and let $\tilde\eta:=D_f\tilde\xi$ on $S^2$. Then
\begin{align*}
   d^{\ell_j}\tilde\xi(z_j)=0,\qquad j^{\ell_j+1}\tilde\xi(z_j) =
   v \in\C^n, \cr
   d^{\ell_j+1}\tilde\xi(z_i) = 0 \text{ for all }i\neq j, \cr
   d^{\ell_j}\tilde\eta(z_i) = 0 \text{ for all }i=1,\dots,k.
\end{align*}
The last equation shows that $\tilde\eta$ belongs to the space
$E^{m-1,p}_0$ defined in Lemma~\ref{lem:Laszlo} (with
$\ell=\ell_j$), so by Lemma~\ref{lem:Laszlo} there exists
$(\xi',Y)\in B^{m,p}_0\oplus T_J\JJ(V)$ such that
$$
   D_f\xi'+\frac{1}{2}Y(f)\circ df\circ j = -\tilde\eta.
$$
It follows that the element $(\xi:=\xi'+\tilde\xi,Y)\in
T_f\BB^{m,p}\oplus T_J\JJ(V)$ satisfies
\begin{gather*}
   D_f\xi+\frac{1}{2}Y(f)\circ df\circ j = 0, \cr
   d^{\ell_j}\xi(z_j)=0, \qquad j^{\ell_j+1}\xi(z_j) = v,
   \cr
   d^{\ell_j+1}\xi(z_i)=0 \text{ for all }i\neq j.
\end{gather*}
Now we use the hypothesis $\ell_j\geq\max\{\ell_i\}-1$, or
$\ell_i\leq\ell_j+1$ for all $i\neq j$, to conclude
$d^{\ell_i}\xi(z_i)=0$ for all $i=1,\dots,k$. Hence $(\xi,Y)$ belongs
to the tangent space $T_{f,J}\tilde\MM_\z^s$ defined
in~\eqref{eq:tangent-space} and $L_j(\xi,Y)
j^{\ell_j+1}_{\C^{n_j}}\xi(z_j) = [v] \in\C^n/\C^{n_j}$.
This proves surjectivity of $L_j$ and hence Lemma~\ref{lem:jets-reg}.
\end{proof}

Lemma~\ref{lem:jets-reg} also includes the cases where some
$\ell_i=-1$, in which case we do not impose any condition at $z_i$.
In particular, for $m-2/p>0$ the universal solution space
$$
   \tilde\MM^s\bigl(\JJ(V)\bigr) :=
   \{(f,J)\in\BB^{m,p}\times\JJ(V)\mid \pb_Jf=0,\ f \text{ simple
   },\ f^{-1}(V)\neq\emptyset\}
$$
of simple holomorphic maps meeting $V$ is a Banach manifold (this
case is also covered by Proposition 3.2.1 in~\cite{MS}). By
Lemma~\ref{lem:jets-reg}, for $m-2/p>\max\{\ell_i\}$,
$$
   \tilde\MM_\z^s\bigl(\JJ(V);\{Z_i,\ell_i\}\bigr)\subset
   \tilde\MM^s\bigl(\JJ(V)\bigr)
$$
is a Banach submanifold that is cut out transversally by the normal
jet evaluation maps $j^\ell_{Z_i}f(z_i)$, $0\leq\ell\leq\ell_i$
(where $j^\ell_{Z_i}f(z_i)$ is defined whenever $f(z_i)\in Z_i$ and
$j^1_{Z_i}f(z_i)=\dots=j^{\ell-1}_{Z_i}(z_i)=0$). In particular,
$\tilde\MM_\z^s\bigl(\JJ(V);\{Z_i,\ell_i\}\bigr)$ has codimension
$$
   2\sum_{i=1}^k(\ell_i+1)\codim_\C Z_i
$$
in $\tilde\MM^s\bigl(\JJ(V)\bigr)$. By the Sard-Smale theorem, there
exists a Baire set\\ $\JJ^\reg(V;\{Z_i\})\subset\JJ(V)$ such that each
$J\in\JJ^\reg(V;\{Z_i\})$ is a regular value of all the projections
$$
   \pi:\tilde\MM_\z^s\bigl(\JJ(V);\{Z_i,\ell_i\}\bigr)\to\JJ(V),\qquad
   (f,J)\mapsto J
$$
for $\ell_1,\dots,\ell_k\geq -1$. Here the complex submanifolds
$Z_1,\dots,Z_k$ are fixed and $\JJ^\reg(V;\{Z_i\})$ depends on
them. Now fix
a homology class $A\in H_2(X;\Z)$ and define for
$J\in\JJ^\reg(V;\{Z_i\})$
\begin{align*}
   \tilde\MM_\z^s(A,J;\{Z_i,\ell_i\}) &:= \{f:S^2\to X\mid \pb_Jf=0,\
   [f]=A \cr
   &\ \ \ d^{\ell_i}f(z_i)\in T_{f(z_i)}Z_i \text{ for } i=1,\dots,k\}.
\end{align*}
By Theorem 3.1.5 in~\cite{MS}, the space $\tilde\MM^s(A,J)$
(corresponding to $\ell_i=-1$ for all $i$) is a smooth manifold of
dimension
$$
   \dim\tilde\MM^s(A,J) = 2n + 2c_1(A).
$$
Hence the preceding discussion shows

\begin{lemma}\label{lem:jets}
Let $Z_1,\dots,Z_k$ be $J_0$-complex manifolds of the
$2n$-dimensional symplectic manifold $(X,\om)$ and $V\subset X$ an
open subset with $V\cap(Z_1\cup\dots\cup Z_k)=\emptyset$. Then
there exists a Baire set $\JJ^\reg(V;\{Z_i\})\subset\JJ(V)$
of tamed almost complex structures that agree with
$J_0$ outside $V$ with the following property: For all integers
$\ell_1,\dots,\ell_k\geq -1$, homology classes $A\in H_2(X;\Z)$, and
$J\in\JJ^\reg(V;\{Z_i\})$, the solution space
$\tilde\MM_\z^s(A,J;\{Z_i,\ell_i\})$ of {\em simple} $J$-holomorphic
maps, passing through $V$ and tangent to $Z_i$ of order $\ell_i$ at
$z_i$, is a smooth submanifold of $\tilde\MM^s(A,J)$ of dimension
$$
   \dim\tilde\MM_\z^s(A,J;\{Z_i,\ell_i\}) = 2n + 2c_1(A) -
   2\sum_{i=1}^k(\ell_i+1)\codim_\C Z_i.
$$
\end{lemma}

\begin{remark}
No assumption is made on the mutual intersections of the complex
submanifolds $Z_i$.
\end{remark}

Lemma~\ref{lem:jets} deals with {\em parametrized}
holomorphic curves and {\em fixed} marked points $z_i$.
However, the vanishing of the $\ell$-th jet is a property
which does not depend on the chosen complex coordinate
on the Riemann surface. Hence we can define the moduli space
$$
   \MM_k^s(A,J;\{Z_i,\ell_i\}) :=
   \tilde\MM_k^s(A,J;\{Z_i,\ell_i\})/\Aut(S^2)
$$
of all (unparametrized) simple $J$-holomorphic maps that pass through
$V$ and are tangent to $Z_i$ of order $\ell_i$ at the (varying!)
marked points $z_i$.

\begin{prop}\label{prop:jets}
Let $Z_1,\dots,Z_k$ be $J_0$-complex submanifolds of the
$2n$-dimen\-sio\-nal symplectic manifold $(X,\om)$ and $V\subset X$ an
open subset with $V\cap(Z_1\cup\dots\cup Z_k)=\emptyset$. Then
there exists a Baire set $\JJ^\reg(V;\{Z_i\})\subset\JJ(V)$
of tamed almost complex structures that agree with
$J_0$ outside $V$ with the following property: For all integers
$\ell_1,\dots,\ell_k\geq -1$, homology classes $A\in H_2(X;\Z)$, and
$J\in\JJ^\reg(V;\{Z_i\})$, the moduli space
$\MM_k^s(A,J;\{Z_i,\ell_i\})$ is a smooth manifold of dimension
$$
   \dim\MM_k^s(A,J;\{Z_i,\ell_i\}) = 2n -6 + 2c_1(A) +2k -
   2\sum_{i=1}^k(\ell_i+1)\codim_\C Z_i.
$$
\end{prop}

\begin{proof}
(1) Consider a set ${\cal U}:= U_1\times U_2\times\dots\times
U_k\subset (S^2)^k$, where $U_1,\dots,U_k$ are disjoint complex
coordinate neighbourhoods of $S^2$. Let
$$
\tilde\MM_{\cal U}^s(A,J;\{Z_i,\ell_i\}):=
\bigcup_{\z\in{\cal U}}\tilde\MM^s_\z(A,J;\{Z_i,\ell_i\})
$$
be the space of simple $J$-holomorphic maps, passing through $V$ and
tangent to $Z_i$ of order $\ell_i$ at $z_i$, where each $z_i$ is
allowed to vary in the open set $U_i$. Fixing holomorphic coordinate
charts on the $U_i$, an argument similar to the proof of
Lemma~\ref{lem:jets} yields a
Baire set $\JJ^\reg_{\cal U}(V;\{Z_i\})\subset\JJ(V)$ such that for
$J\in\JJ^\reg_{\cal U}(V;\{Z_i\})$ this is a manifold of dimension
$$
\dim\tilde\MM_{\cal U}^s(A,J;\{Z_i,\ell_i\}) = 2n + 2c_1(A) +2k -
2\sum_{i=1}^{k}(\ell_i+1)\codim_{\C}Z_i.
$$
(2) Next note that we can cover the space $\tilde\MM_k$ of $k$-tuples
$(z_1,\dots,z_k)$ of distinct points in $S^2$ by a countable family of
such open sets ${\cal U}:= U_1\times\dots\times U_k$, were
$U_1,\dots,U_k\subset S^2$ are disjoint coordinate
neighbourhoods. Indeed, for any tuple of distinct points
$\z=(z_1,\dots,z_k)$ whose real and imaginary parts are rational or equal
to $\infty$, let $\rho:=\min\{d(z_i,z_j)\mid i\neq j\}$ and $U_i$ be the
open ball around $z_i$ of radius $\rho/2$. Here all distances are
measured with respect to the standard metric of the unit sphere. This
yields a countable family of open sets ${\cal
  U}(\z):=U_1\times\dots\times U_k$ of the required form. To show
that this is a covering, let $\z=(z_1,\dots,z_k)$ be any $k$-tuple of
distinct points in $S^2$ and let $\rho:=\min\{d(z_i,z_j)\mid i\neq j\}$
be the minimal distance. Choose a tuple $\z'=(z'_1,...,z'_k)$ of distinct
rational points such that $d(z_i,z'_i)<\rho/4$ for all $i$ and thus
$\rho':=\min\{d(z'_i,z'_j)\mid i\neq j\}>\rho/2$. Then we have $\z\in
{\cal U}(\z')= U_1'\times ... \times U_k'$, where $U_i'$ is the ball
around $z_i'$ of radius $\rho'/2>\rho/4$.

(3) Now for each set ${\cal U}$ of the countable covering in
part (2) let $\JJ^\reg_{\cal U}(V;\{Z_i\})\subset\JJ(V)$ be the Baire set
provided by part (1). Then the countable intersection $\JJ^\reg(V;\{Z_i\}):=
\cap_{\cal U}\JJ^\reg_{\cal U}(V;\{Z_i\})$ is again a Baire set and for
$J\in\JJ^\reg(V;\{Z_i\})$, the solution space
$\tilde\MM_k^s(A,J;\{Z_i,\ell_i\})$ is a manifold of dimension
$$
   \dim\tilde\MM_k^s(A,J;\{Z_i,\ell_i\}) = 2n + 2c_1(A) +2k -
   2\sum_{i=1}^k(\ell_i+1)\codim_\C Z_i.
$$
Since the automorphism group $\Aut(S^2)$ acts freely and properly on
this space (see~\cite{MS}), the quotient $\MM_k^s(A,J;\{Z_i,\ell_i\})$
is a manifold of dimension
$$
   \dim\MM_k^s(A,J;\{Z_i,\ell_i\}) = 2n -6 + 2c_1(A) + 2k -
   2\sum_{i=1}^k(\ell_i+1)\codim_\C Z_i.
$$
\end{proof}

If we replace $\JJ(V)$ in Lemma~\ref{lem:Laszlo} by the space
$\JJ_{S^2}(V)$ of almost complex structures depending on points in
the domain $S^2$, then we can drop the hypothesis that $f$ is simple.
Indeed, if $f$ is constant then the hypotheses
$f^{-1}(V)\neq\emptyset$ and $V\cap Z_i=\emptyset$ imply $\ell_i=-1$
for all $i$ and the result follows from Lemma~\ref{lem:surj-const}. 
So suppose that $f$ is nonconstant. Then in the proof of
Lemma~\ref{lem:Laszlo} we just need to replace the cutoff construction
for the section $Y$ in the target by a cutoff in the domain
$S^2\setminus\{z_1,...,z_k\}$ as in the proof of
Lemma~\ref{lem:surj-nonconst}. We can also make this cutoff depend
on $\z=(z_1,...,z_k)$ and vanish outside a sufficiently small neighborhood
to ensure positivity of the pairing with the $\eta$ as in the proof of
Lemma~\ref{lem:Laszlo} between equations (\ref{eq:Yf}) and (\ref{eq:Lambda}).
As an easy instance of Lemma~\ref{lem:J-smooth} we obtain a coherent section
$Y$ depending on $\pi_I(\z)\in \MM_{|I|+1}$ for any subset $I\subset
\{1,\dots,k\}$ with $|I|\geq 3$. With the modified
Lemma~\ref{lem:Laszlo}, the proofs of Lemma~\ref{lem:jets-reg} and 
Proposition~\ref{prop:jets} carry over with minor adaptations and we
obtain 

\begin{prop}\label{prop:jets1}
Let $Z_1,\dots,Z_k$ be $J_0$-complex manifolds of the
$2n$-dimensio\-nal symplectic manifold $(X,\om)$, $V\subset X$ an
open subset with $V\cap(Z_1\cup\dots\cup Z_k)=\emptyset$, and
$I\subset\{1,\dots,k\}$ with $|I|\geq 3$. Then
there exists a Baire set
$\JJ^\reg_{|I|+1}(V;\{Z_i\})\subset\JJ_{|I|+1}(V)$ 
of coherent domain-dependent tamed almost
complex structures that agree with $J_0$ outside $V$ with the
following property: For all integers $\ell_1,\dots,\ell_k\geq -1$,
homology classes $A\in H_2(X;\Z)$, and
$J\in\JJ^\reg_{|I|+1}(V;\{Z_i\})$, the
moduli space $\MM_k^*(A,J;\{Z_i,\ell_i\})$ of (not necessarily
simple!) $J$-holomorphic maps, passing through $V$ and tangent to
$Z_i$ of order $\ell_i$ at the (varying) points $z_i$, is a smooth
submanifold of $\MM^*(A,J)$ of dimension
$$
   \dim\MM_k^*(A,J;\{Z_i,\ell_i\}) = 2n - 6 + 2c_1(A) + 2k -
   2\sum_{i=1}^k(\ell_i+1)\codim_\C Z_i.
$$
\end{prop}

Applying this proposition to the special points on one vertex
$\alpha$ of an $I$-stable tree, 
we obtain

\begin{cor}\label{cor:jets-tree}
Let $T$ be an $I$-stable $k$-labelled tree with $A_\alpha\neq 0$ for a
unique vertex $\alpha$ (so the reduced index set $R$ labels 
the special points on $\alpha$), and $J_0\in\JJ$. 
Let $Z_i$, $i\in R$ be $J_0$-complex manifolds of $X$, $V\subset
X$ an open subset with $V\cap Z_i=\emptyset$ for all $i\in R$, and
$\ell_i\geq -1$ integers for $i\in R$. 
Then there exists a Baire set
$\JJ_{|I|+1}^\reg(V)\subset\JJ_{|I|+1}(V)$ (depending on $T,A_\alpha$
and the $Z_i,\ell_i$) such that for all 
$J\in\JJ_{|I|+1}^\reg(V)$ the moduli space
$\MM^*_T(A_\alpha,J;\{Z_i,\ell_i\})$ of (not necessarily
simple) stable $J$-holomorphic maps modelled over $T$, passing through
$V$ and tangent to $Z_i$ of order $\ell_i$ at the (varying) special points
$z_{\alpha i}$ on the component $\alpha$, is a
smooth submanifold of $\MM^*_T(A_\alpha,J)$ of dimension
$$
   \dim\MM^*_T(A_\alpha,J;\{Z_i,\ell_i\}) = 2n - 6 +
   2c_1(A_\alpha) + 2k - 2e(T) -
   2\sum_{i\in R}(\ell_i+1)\codim_\C Z_i.
$$
\end{cor}

As in Corollary~\ref{cor:surj-tree2}, the dimension formula can be
improved by forgetting the marked points outside $R$ and stabilizing
(which yields the tree with the unique vertex $\alpha$):

\begin{cor}\label{cor:jets-tree2}
In the situation of Corollary~\ref{cor:jets-tree}, 
the evaluation map $\ev^R$ factors as 
$$
   \ev^R:\MM^*_T(A_\alpha,J;\{Z_i,\ell_i\})\to
   \MM^*_{|R|}(A_\alpha,J;\{Z_i,\ell_i\})\to \prod_{i\in R}Z_i
$$
through a smooth manifold of dimension 
$$
   \dim\MM^*_{|R|}(A_\alpha,J;\{Z_i,\ell_i\}) = 2n - 6 +
   2c_1(A_\alpha) + 2|R| -
   2\sum_{i\in R}(\ell_i+1)\codim_\C Z_i.
$$
\end{cor}

\section{Intersection numbers}\label{sec:int}

A {\em complex hypersurface} is an (almost) complex submanifold $Y$ in an
almost complex manifold $(X,J)$ of complex codimension one. Let
$f:(\Sigma,j)\to(X,J)$ be a holomorphic map from a Riemann surface.
Suppose that $f(z)\in Y$ is an {\em isolated intersection} of $f$
with $Y$. This means that there exist a closed disk $D\subset\Sigma$
around $z$ and a closed $(2n-2)$-ball $B\subset Y$ around $f(z)$
such that $f^{-1}(B)\cap D=\{z\}$. Then the {\em local
intersection number}
$$
   \iota(f,Y;z) := (f|_D)\cdot B.
$$
of $f$ with $Y$ at $z$ is defined by smoothly perturbing $f$ and
counting with signs (and independent of a sufficiently small
perturbation). More generally, suppose that $\Sigma$ and $Y$ are
compact and connected with (possibly empty) boundaries such that
$f^{-1}(Y)\cap\p\Sigma=f^{-1}(\p Y)\cap \Sigma=\emptyset$. Then we
have a well-defined {\em intersection number}
$$
   \iota(f,Y) := (f|_\Sigma)\cdot Y.
$$
of $f$ with $Y$.

\begin{prop}\label{prop:int}
Let $f:(\Sigma,j)\to(X,J)$ be holomorphic, $Y\subset X$ a complex
hypersurface such that $f^{-1}(Y)\cap\p\Sigma=f^{-1}(\p Y)\cap
\Sigma=\emptyset$ and $f(\Sigma)\not\subset Y$.
Then the set $f^{-1}(Y)$ is finite and
$$
   \iota(f,Y) = \sum_{z\in f^{-1}(Y)}\iota(f,Y;z).
$$
At each intersection point $z$, $f$ is tangent to $Y$ of some finite
order $\ell\geq 0$ with
$$
   \iota(f,Y;z) = \ell + 1.
$$
In particular, each local intersection number $\iota(f,Y;z)$ is
positive.
\end{prop}

\begin{proof}
We derive a local representation for $f$ in the direction normal to
$Y$ near an intersection point $z$. Pick local holomorphic
coordinates for $\Sigma$ near $z$ mapping $z$ to $0\in\C$ and local
coordinates for $X$ near $f(z)\in Y$ mapping $f(z)$ to $0$ and $Y$
to $\C^{n-1}\times\{0\}\subset\C^n$. Moreover, we can achieve that
$\{0\}\times\C$ are complex directions at points of
$\C^{n-1}\times\{0\}$. Denote coordinates on $\C^n=\C^{n-1}\oplus\C$
by $x=(y,u)$. Then $J$ corresponds to an almost complex structure on
$\C^n$ of the form
$$
   J(x) = \left(\begin{matrix}
   K(x) & a(x) \\ b(x) & j(x)\end{matrix}\right),\qquad
      J(y,0) = \left(\begin{matrix}
   K(y,0) & 0 \\ 0 & j(y,0)\end{matrix}\right),
$$
so it satisfies
$$
   a(y,0)=0,\qquad b(y,0)=0,\qquad K(y,0)^2=-\id,\qquad j(y,0)^2=-\id.
$$
The map $f$ corresponds to a map $f\!=\!(g,h)\!:\!\C\supset U\!\to\!\C^n\!=\!\C^{n-1}\!\oplus\!\C$
satisfying $f_s+J(f)f_t=0$. So its last component $h:\C\supset U\to\C$
satisfies
$$
   h_s+j(g,h)h_t+b(g,h)g_t = 0.
$$
Since $b(g,0)=0$, we can write
$$
   j(g,h)=j(g,0)+c(g,h)(h),\qquad b(g,h)=d(g,h)(h)
$$
with smooth maps
$$
   c:\C^n\to\Hom_\R\bigl(\C,\Hom_\R(\C,\C)\bigr), \quad
   d:\C^n\to\Hom_\R\bigl(\C,\Hom_\R(\C^{n-1},\C)\bigr).
$$
Then 
$$
   j(g,h)h_t = J_0(z)h_t + C(z)(h),\qquad b(g,h)g_t = D(z)(h)
$$
with the $\R$-linear maps $J_0(z):=j\bigl(g(z),0\bigr)$ and
$$
   C(z)(\cdot) :=
   c\bigl(g(z),h(z)\bigr)(\cdot)h_t(z),\quad D(z)(\cdot) :=
   d\bigl(g(z),h(z)\bigr)(\cdot)g_t(z).
$$
This shows that $h$ satisfies the equation
$$
   h_s + J_0(z)h_t + A(z)(h) = 0
$$
with $\R$-linear maps $J_0(z)$ and $A(z):=C(z)+D(z)$ such that
$J_0(z)^2\!=\!-\id$ for all $z$. By the version of the Carleman similarity
principle proved in~\cite{FHS}, there exists a continuous matrix
valued function $\Phi:\C\supset V\to GL(2,\R)$ and a holomorphic
function $\sigma:\C\supset V\to\C$ such that
$$
   h(z)=\Phi(z)\sigma(z),\qquad \sigma(0)=0,\qquad
   J_0(z)\Phi(z)=\Phi(z)i.
$$
Thus $h$ looks locally like the holomorphic function $\sigma$. The
proposition now is an easy consequence of this representation
(see~\cite{FHS}): If $\sigma\equiv 0$ then, by unique continuation,
we would conclude $f(\Sigma)\subset Y$, which is excluded by
hypothesis. Thus $\sigma(z)=z^{\ell+1}\tau(z)$ for an integer
$\ell\geq 0$ and a holomorphic function $\tau$ satisfying
$\tau(0)\neq 0$. This shows that $0$ is an isolated intersection
point of $f$ with $Y$. Moreover, $f$ is tangent to $Y$ at $0$ of
finite order $\ell$. 
Finally, the local intersection number of $f$
with $Y$ equals the local mapping degree of $h$ at $0$. Since
$J_0(z)\Phi(z)=\Phi(z)i$, the matrix $\Phi(z)$ is orientation
preserving and hence homotopic to the identity in $GL(2,\R)$. So the
local mapping degree of $h$ is given by the winding number of $\sigma$,
which is just $\ell+1$. 
\end{proof}

Now let $\f=(f_\alpha):\Sigma_\z\to X$ be a nonconstant genus zero nodal
$J$-holomorphic curve with $k\geq 1$ marked points modelled over a
tree $T$. Let $Y\subset(X,J)$ be a closed complex hypersurface.
Suppose that no nonconstant component $f_\alpha$ is contained in
$Y$. Then we define the local intersection number $\iota(f,Y;z_i)$
at a marked point $z_i$ as follows: If $f_{\alpha_i}$ is nonconstant
set $\iota(\f,Y;z_i):=\iota(f_{\alpha_i},Y;z)$. If $f_{\alpha_i}$ is
constant let $T_1\subset T$ be the maximal ghost tree containing
$\alpha_i$. Let $T_2\subset T\setminus T_1$ be the set of vertices
adjacent to $T_1$, so $f_\beta$ is nonconstant for all $\beta\in
T_2$. Recall that $z_{\beta i}\in S_\beta$ is the nodal point
connecting $\beta$ to $T_1$. Then set
$$
   \iota(\f,Y;z_i):=\sum_{\beta\in
   T_2}\iota(f_\beta,Y;z_{\beta i}).
$$
Note that in this definition we allow for constant components
$f_\alpha$ to be contained in $Y$. Also note that
$\iota(\f,Y;z_i)=\iota(\f,Y;z_j)$ if $z_i$ and $z_j$ belong to the
same ghost tree.

\begin{lemma}\label{lem:lim-int}
Let $Y\subset (X,J)$ be a closed complex hypersurface. Let
$(\z^\nu,\f^\nu)$ be a sequence of nonconstant genus zero nodal
$J$-holomorphic
maps with $k\geq 1$ marked points converging to $(\z,\f)$ in the
Gromov topology. Suppose that no nonconstant components of $\f^\nu$
and $\f$ are contained in $Y$.

(a) If $z_1$ lies on a nonconstant component of $f$, then
$$
   \iota(\f,Y;z_1) \geq
   \limsup_{\nu\to\infty}\iota(\f^\nu,Y;z_1^\nu).
$$
(b) If $z_1$ lies on a constant component of $f$, let $T_1\subset
T$ be the maximal ghost tree containing the vertex $\alpha_1$. Then
$$
   \iota(\f,Y;z_1) \geq
   \limsup_{\nu\to\infty}\sum_{\alpha_i\in T_1}\iota(\f^\nu,Y;z_i^\nu),
$$
where for each maximal ghost tree $T'\subset T^\nu$
at most one of the $z_i^\nu$ with $\alpha_i\in T'$ is
counted in the sum.
\end{lemma}

\begin{proof}
Pick a closed ball $B$ in $X$ around $\f(z_1)$ containing no other
intersection points of $Y$ with $\f$. Since each $\f^\nu$ intersects
$Y$ in only finitely many points, we have at most countably many
intersection points of $Y$ with all the $\f^\nu$. Thus we can choose
$B$ such that $\p B$ contains no intersection point of $Y$ with any
of the $\f^\nu$. Moreover, we can choose $B$ such that $\p B$ is
transverse to $\f$ and all the $\f^\nu$ (in particular, $\p B$
contains no nodal points of $\f$ or any of the $\f^\nu$). Let
$\Sigma\subset\Sigma_\z$ be the connected component of $\f^{-1}(B)$
containing $z_1$. Similarly, let $\Sigma^\nu\subset\Sigma_{\z^\nu}$
be the connected component of $(\f^\nu)^{-1}(B)$ containing
$z_1^\nu$. By definition of Gromov convergence, after suitable
reparametrization $\f^\nu|_{\Sigma^\nu}$ converges uniformly to
$\f|_\Sigma$. (Note that the domains $\Sigma^\nu$ and $\Sigma$ may
have different numbers of nodes due to pinching, but since the maps
are constant across the nodes uniform convergence makes sense.) So
by definition of the
intersection number,
$$
   \left(\f^\nu|_{\Sigma^\nu}\right)\cdot Y\to
   \left(\f|_\Sigma\right)\cdot Y
$$
as $\nu\to\infty$.

Now consider first case (a). Choosing the ball $B$ small enough, we
can ensure that $z_1$ is the only point in $\Sigma$ with $f(z_1)\in
Y$. This implies
$$
   \left(\f|_\Sigma\right)\cdot Y = \iota(\f,Y;z_1).
$$
By definition of Gromov convergence, for $\nu$ large the domain
$\Sigma^\nu$ contains the point $z_1^\nu$. It may also contain other
intersection points of $f^\nu$ with $Y$. But positivity of the local
intersection numbers in Proposition~\ref{prop:int} implies
$$
   \left(\f^\nu|_{\Sigma^\nu}\right)\cdot Y \geq
   \iota(\f^\nu,Y;z_1^\nu),
$$
and case (a) follows.

Next consider case (b). Choosing the ball $B$ small enough, we can
ensure that all points $z\in\Sigma$ with $f(z)\in Y$ belong to
the tree $T_1$. This implies
$$
   \left(\f|_\Sigma\right)\cdot Y = \iota(\f,Y;z_1).
$$
By definition of Gromov convergence, for $\nu$ large the domain
$\Sigma^\nu$ contains all points $z_i^\nu$ with $\alpha_i\in T_1$.
So as above, positivity of intersections implies
$$
   \left(\f^\nu|_{\Sigma^\nu}\right)\cdot Y \geq
   \sum_{\alpha_i\in T_1}\iota(\f^\nu,Y;z_i^\nu),
$$
where the sum is to be interpreted as explained in the statement of
the lemma. This proves case (b).
\end{proof}

\section{Symplectic hypersurfaces}\label{sec:hyp}

{\bf Existence and asymptotic uniqueness. }
Let $(X,\om)$ be a closed symplectic manifold such that $\om$
represents an integral cohomology class $[\om]\in H^2(X;\Z)$. By a
{\em hypersurface} we mean a closed codimension 2 submanifold
$Y\subset X$ whose
Poincar\'e dual in $H^2(X;\Z)$ equals $D[\om]$, where $D\in\N$ is
called the {\em degree} of $Y$. Fix a compatible almost complex
structure $J$ on $X$. We say that a hypersurface $Y\subset X$ of degree
$D$ is {\em approximately $J$-holomorphic} if its K\"ahler angle (see
the following paragraph for the precise definition) at each point 
is at most $CD^{-1/2}$, where $C$ is a constant independent of $D$. We say
that two hypersurfaces $Y,\bar Y\subset X$ of degrees $D,\bar D$
intersect {\em $\eps$-transversally} if they intersect transversally
and at each intersection point their tangent spaces have minimal
angle (also defined in the next paragraph) at least $\eps$ for an
$\eps>0$ independent of the degrees. 
The following theorem is an easy adaptation of the results
in~\cite{Do} and~\cite{Aur}.

\begin{theorem}[Donaldson~\cite{Do},
  Auroux~\cite{Aur}]\label{thm:Donaldson}
\hfill\\
{ (Existence). }For each sufficiently large integer $D$
there exists an approximately $J$-holomorphic hypersurface $Y\subset
X$ whose Poincar\'e dual in $H^2(X;\Z)$ equals $D[\om]$.

{ (Stabilization). }Let $Y\subset X$ be a $J$-holomorphic (not just
approximately!) hypersurface of degree $D$. Then there exists an
$\eps>0$ such that for each sufficiently large
integer $\bar D$ there exists an approximately
$J$-holomorphic hypersurface $\bar Y$ of degree $\bar D$ intersecting
$Y$ $\eps$-transversally.

{ (Uniqueness). }Let $Y_i$, $i=0,1$ be two approximately
$J_i$-holomorphic hypersurfaces arising from the construction in the
(Existence) part with the same degree $D$. Connect $J_0,J_1$ by a path of
compatible almost complex structures $J_t$. If $D$ is sufficiently
large, then $Y_0,Y_1$ are isotopic through approximately
$J_t$-holomorphic hypersurfaces. Moreover, this isotopy can be realized
through symplectomorphisms of $X$.
\end{theorem}

\begin{proof}
(Existence) is proved in~\cite{Do} and (Uniqueness) in~\cite{Aur}.
The (Stabilization) property is not stated explicitly in the
literature, but it follows from the proof in~\cite{Do}. For the sake
of completeness, let us indicate the required modifications. 

We follow the notation of~\cite{Aur}. 
Fix a Hermitian line bundle $L\to X$ with a unitary connection $A$ of
curvature $-i\om$. For $k\in\N$ we consider sections $s$ in the $k$-th
tensor power $L^k$ of $L$. We will measure all quantities with respect
to the rescaled metric $g_k:=kg$ on $X$, where
$g:=\om(\cdot,J\cdot)$. All constants are assumed to be independent of
$k$. Following~\cite{Do2}, for positive constants 
$C,\eta$ we say that $s$ is {\em $C$-bounded} if 
$$
   \|s\|_{C^3}<C,\qquad \|\pb s\|_{C^2}<Ck^{-1/2},
$$
and $s$ is {\em $\eta$-transverse} if $|s(x)|<\eta$ implies
$\nu(\nabla_xs)\geq\eta$. Here $\nu(\nabla_xs)^{-1}$ denotes the
minimal norm of a right inverse of $\nabla_xs$. Note that $s$ induces
a section $s|_Y$ in the bundle $L^k|_Y$. We say that $s$ is {\em
  $\eta$-transverse over $Y$} if $s|_Y$ is $\eta$-transverse. 

{\em Claim. }There exist $C,\eta>0$ and a sequence of sections $s_k$
of $L^k$ that are $C$-bounded, $\eta$-transverse, and
$\eta$-transverse over $Y$. 

From this claim the (stabilization) property easily follows: As
in~\cite{Do}, $C$-bounded\-ness and $\eta$-transversality imply that for
$k$ sufficiently large, $\bar Y:=s_k^{-1}(0)$ is a (transversally cut
out) submanifold whose K\"ahler angle is arbitrarily small; see
Lemma~\ref{lem:angle} (b) and ~\cite{Do}, remark after Theorem
5. Next, $\eta$-transversality over $Y$ and $C$-boundedness  
imply that $\bar Y$ and $Y$ intersect $\eps$-transversally for 
$\eps=\eta/C$, see Lemma~\ref{lem:angle3} (a) below. 
Finally, $\eps$-transversality of $Y$
and $\bar Y$ and smallness of their K\"ahler angles ($Y$ has K\"ahler
angle zero) imply that their intersection is positive and there
exists a compatible almost complex structure $\bar J$ which is
$C^0$-close to $J$ and makes both $Y$ and $\bar Y$ holomorphic; see
Lemma~\ref{lem:angle3} (b) below. 

To prove the claim, we follow the scheme in~\cite{Do}. First
recall the construction of local $C$-bounded sections $\sigma_{k,x}$
supported near a point $x\in X$. Pick a Darboux chart
$\psi_x:\C^n\supset U\to V\subset X$  with $\psi_x(0)=x$ whose
differential at $0$ is complex linear, and a trivialization of
$\psi_x^*L$ in which the connection $A$ is given by
$$
   A  = \frac{1}{4}\sum_{i=1}^n(z_id\bar z_i-\bar z_idz_i). 
$$
Then $f(z):=\exp(-k|z|^2/4)$ defines a holomorphic section of $L^k$
over $U\subset\C^n$. We multiply $f$ by a cutoff function $\beta$ and
transfer it to $X$ via 
$$
   \sigma_{k,x} := (\beta f)\circ\tilde\psi^{-1},\qquad
   \tilde\psi(z):=\psi(k^{-1/2}z). 
$$
By choosing the $\psi_x$ for all $x$ with uniformly bounded
$C^3$-norm, one can achieve that the sections $\sigma_{k,x}$ are
$C$-bounded.  

For $x\in Y$ we can choose $\psi_x$ to map the subspace $\C^{n-1}$ to
$Y$, so the $(\psi_x)|_{\C^{n-1}}$ define Darboux charts for $Y$ with the
same properties. This implies that the $\sigma_{k,x}$ restrict to
$C$-bounded sections on $Y$. 

Now we globalize the construction. For each $k$ pick finitely many
points $x_i\in X$ such that the $g_k$-balls $B_1(x_i)$ of
radius $1$ around the $x_i$ cover $X$. The required number of
points is of order $k^n$, where $\dim X=2n$. We can arrange that some
of the points $x_i$ lie on $Y$ and their unit balls cover $Y$. 
The desired section $s$ is obtained in the form
$$
   s(x) = \sum_{x_i}a_{k,x_i}\sigma_{k,x_i}
$$
for suitable complex numbers $a_{k,x_i}$. Here the coefficients
$a_{k,x_i}$ are defined inductively to achieve $\eta_i$-transversality
(for suitable $\eta_i>0$)   
on $B_1(x_i)$ without destroying it on the previous balls.
The proof
based on Theorem 12 in~\cite{Do} shows that for $x_i\in Y$ we can
choose $a_{k,x_i}$ such that $\eta_i$-transversality holds on
$B_1(x_i)$ as well as on $B_1(x_i)\cap Y$. Hence the resulting
$C$-bounded section $s$ is $\eta$-transverse (with
$\eta=\min\{\eta_i\}$) on $X$ as well as over
$Y$. This concludes the proof of the (stabilization) property. 
\end{proof}

\begin{remark}
It follows directly from the definitions that the hypersurfaces
constructed in Theorem~\ref{thm:Donaldson} have the following
properties:

{ (Existence). }$Y$ is symplectic and $\bar
J$-holomorphic for a compatible almost complex structure $\bar J$
arbitrarily $C^0$-close to $J$.

{ (Stabilization). }$\bar Y$ is symplectic and there exists a
compatible almost complex structure, $\bar J$, arbitrarily $C^0$-close
to $J$ and coinciding with $J$ on $Y$, such that $Y,\bar Y$ are both
$\bar J$-holomorphic.

{ (Uniqueness). }The $Y_t$ are symplectic and $\bar
J_t$-holomorphic for a path of compatible almost complex structures
$\bar J_t$ arbitrarily $C^0$-close to $J_t$.
\end{remark}

{\bf Linear algebra. }
Recall the linear algebra underlying Donaldson's construction
in~\cite{Do}. Consider a complex vector space $(V,J)$ of complex
dimension $n$ with Hermitian metric
$$
   (\cdot,\cdot)=\la\cdot,\cdot\ra-i\om.
$$
Define the {\em K\"ahler angle} $\theta(W)\in[0,\pi]$ of an
oriented $\R$-linear subspace $W\subset V$ of real dimension $2k$ with
respect to the Hermitian sructure $(\om,J)$ by
$$
   \theta(W):=\theta(W;\om,J) :=
   \cos^{-1}\left(\frac{\om^k|_W}{k!\Om_W}\right), 
$$
where $\Om_W$ is the volume form on $W$ defined by the metric and
orientation. We will usually drop $\om$ and/or $J$ in
$\theta(W;\om,J)$ when they are clear from the context.  
Note that $\theta(-W)=\pi-\theta(W)$, where $-W$ is the
subspace $W$ with the opposite orientation. Define the {\em angle}
$\angle\in[0,\pi/2]$ between two nonzero vectors $x,y\in V$ resp.~between
a nonzero vector $y$ and an $\R$-linear subspace $X$, and the {\em
  maximal angle} between two subspaces $X,Y\subset V$ by
\begin{gather*}
   \angle_M(x,y) := \cos^{-1}\left(\frac{|\la
   x,y\ra|}{|x|\,|y|}\right),\qquad
   \angle_M(X,y) := \inf_{0\neq x\in X}\angle_M(x,y), \cr
   \angle_M(X,Y) := \sup_{0\neq y\in Y}\angle_M(X,y) =
   \sup_{0\neq y\in Y}\inf_{0\neq x\in
   X}\angle_M(x,y).
\end{gather*}
Note that $\angle_M(X,Y)$ depends on the order of the subspaces $X,Y$.
It equals $\pi/2$ whenever $\dim Y>\dim X$ and $0$ iff $Y\subset X$. 
Denote by $X^\perp$ and $X^\om$ the orthogonal resp.~$\om$-orthogonal
complement. 

\begin{lemma}\label{lem:angle}
(a) $\theta(W)=0$ iff $W$ is complex linear and the orientation of $W$
agrees with the complex orientation, $\theta(W)=\pi$ iff $W$ is
complex linear and the orientation of $W$
is opposite to the complex orientation, $\theta(W)<\pi/2$ iff $W$
is symplectic, and $\theta(W)>\pi/2$ iff $-W$ is symplectic.

(b) Let $A:V\to\C$ be $\C$-linear and $B:V\to\C$ be $\C$-antilinear
such that $A+B:V\to\C$ is surjective. Then the K\"ahler angle of
$W:=\ker(A+B)$ (oriented via $A+B$) satisfies
$$
   \tan\theta(W) = \frac{2\sqrt{|A|^2|B|^2-|(A,\bar
   B)|^2}}{|A|^2-|B|^2}
$$
In particular, if $|B|<|A|$ then $W$ is symplectic.

(c) For two $\R$-linear subspaces $X,Y\in V$,
$$
   \angle_M(X,Y) = \angle_M(JY,JX) = \angle_M(Y^\perp,X^\perp) =
   \angle_M(Y^\om,X^\om).
$$
(d) The K\"ahler angle of an oriented $\R$-linear subspace $W$ of
dimension $2$ or $2n\!-\!2$ satisfies
$$
  \min\{\theta(W),\pi-\theta(W)\} = \angle_M(W,JW) = \angle_M(JW,W).
$$
(e) The K\"ahler angle of an oriented $\R$-linear subspace $W$ of
dimension $2$ or $2n-2$ satisfies
$$
   \theta(W) = \theta(JW) = \theta(W^\perp) = \theta(W^\om).
$$
\end{lemma}

\begin{proof}
(a) and (b) are proved in~\cite{Do}. The first equality in (c) is
clear because $J$ preserves the inner product $\la\cdot,\cdot\ra$, and
the third equality follows from the first two because
$X^\om=(JX)^\perp$. For the second equality, note that
$\angle_M(X,y)=\pi/2-\angle_M(X^\perp,y)$. Using this twice, we find
\begin{align*}
   \angle_M(X,Y) &= \sup_{y\in Y}\angle_M(X,y) = \sup_{y\in
   Y}\left(\frac{\pi}{2}-\angle_M(X^\perp,y)\right) \cr
   &= \frac{\pi}{2} - \inf_{y\in Y}\inf_{x^\perp\in
   X^\perp}\angle_M(x^\perp,y) = \frac{\pi}{2} - \inf_{x^\perp\in
   X^\perp}\angle_M(Y,x^\perp) \cr
   &= \frac{\pi}{2} - \inf_{x^\perp\in
   X^\perp}\left(\frac{\pi}{2}-\angle_M(Y^\perp,x^\perp)\right) =
   \sup_{x^\perp\in X^\perp}\angle_M(Y^\perp,x^\perp) \cr
   &= \angle_M(Y^\perp,X^\perp).
\end{align*}
For (d) first note that for a real subspace of arbitrary dimension
$\theta(W)=\pi/2$ iff $W\cap W^\om= W\cap(JW)^\perp\neq\{0\}$, or
equivalently $JW\cap W^\perp\neq\{0\}$. But this is equivalent to
$\angle_M(W,JW)=\pi/2$. For the remaining case $\theta(W)\neq\pi/2$
we may assume, after changing the orientation on $W$ if necessary,
that $W$ is symplectic, i.e.~$\theta(W)<\pi/2$.

Now assume $\dim W=2n-2$. Define the complex subspace
$W_0:=W\cap JW$. If $W_0=W$ all terms in the equation in (d) are
$0$. So assume $W_0\neq W$, thus $W_0$ has dimension $2n-4$. Denote
by $V_1$ its orthogonal complement in $V$ and set $W_1:=W\cap V_1$,
thus $\dim W_1=2$. For $0\neq x=x_0+x_1\in W=W_0\oplus W_1$ we have
$$
   \pi_W(Jx) = Jx_0+\pi_{W_1}Jx_1,
$$
where $\pi_W,\pi_{W_1}$ are the orthogonal projections onto $W,W_1$.
If $x_1=0$ we have $\pi_W(Jx)=Jx$ and thus $\angle_M(W,Jx)=0$. If
$x_1\neq 0$ let $y_1\in W_1$ be the unique vector such that
$|y_1|=|x_1|$, $y_1\perp x_1$, and $(x_1,y_1)$ is an oriented basis
of $W_1$. Since $\la x_1,\pi_{W_1}Jx_1\ra = \la x_1,Jx_1\ra = 0$, we
have $\pi_{W_1}Jx_1=\cos\theta\,y_1$ for a unique $\theta\in[0,\pi]$.
It follows that
$$
   \cos\theta = \frac{\la\pi_{W_1}Jx_1,y_1\ra}{|y_1|^2} =
   \frac{\om(x_1,y_1)}{|x_1|\,|y_1|} = \cos\theta(W_1) =
   \cos\theta(W),
$$
which is positive because $W$ is symplectic. On the other hand,
$$
   |\pi_WJx|^2 = |Jx_0+\cos\theta\,y_1|^2 =
   |x_0|^2+\cos^2\theta|y_1|^2 \geq \cos^2\theta|x|^2,
$$
with equality iff $x_0=0$ or $\cos^2\theta = 1$. The latter is
excluded since we assumed $W_0\neq W$.  This shows that
$$
   \cos\angle_M(W,Jx) = \frac{|\pi_WJx|}{|x|} \geq \cos\theta(W),
$$
and hence $\angle_M(W,Jx)\leq\theta(W)$, for all $x\in W$, with
equality iff $x_0=0$. Taking the supremum over $x\in W$, this proves
(d) in the case $\dim W=2n-2$. The case $\dim W=2$ follows by an
easier version of the same argument, setting $W_0=\{0\}$ and $W_1=W$.  
Part (e) is an immediate consequence of (c) and (d).
\end{proof}

\begin{remark}
Property (e) of the preceding lemma can be used to derive property
(b) as follows. Write
$$
   Az = (z,a),\qquad Bz=(b,z)
$$
for vectors $a,b\in V$. Thus under the identification of $V^*$ with
$V$ via the inner product, $A$ corresponds to $a$ and $B$ to $\bar
b$. Surjectivity of $A+B$ is equivalent to $a\neq\lambda b$ for
$\lambda\in S^1$. Then $W=\ker(A+B)$ is the orthogonal complement of
the 2-dimensional oriented subspace
$$
   W^\perp = {\rm span}\{a+b,Ja-Jb\}.
$$
To determine $\theta:=\theta(W)=\theta(W^\perp)$, we compute
\begin{gather*}
   \om(a+b,Ja-Jb) = |a|^2-|b|^2, \cr
   |a+b|^2|Ja-Jb|^2-\la a+b,Ja-Jb\ra^2 = (|a|^2+|b|^2)^2-4|(a,b)|^2.
\end{gather*}
If $|a|=|b|$ this shows that $\theta(W)=\pi/2$. If $|a|\neq|b|$ we
conclude (b) from
\begin{gather*}
   \cos\theta =
   \frac{|a|^2-|b|^2}{\sqrt{(|a|^2+|b|^2)^2-4|(a,b)|^2}}, \qquad
   \sin\theta =
   \frac{2\sqrt{|a|^2|b|^2-|(a,b)|^2}}{\sqrt{(|a|^2+|b|^2)^2-
   4|(a,b)|^2}}, \cr
   \tan\theta = \frac{2\sqrt{|a|^2|b|^2-|(a,b)|^2}}{|a|^2-|b|^2}\ .
\end{gather*}
\end{remark}

Let $(V,J,\om)$ be a Hermitian vector space as above. To a nonzero
skew-symmetric bilinear form $\sigma$ on $V$ we associate the
quantities
\begin{align*}
   \alpha(\sigma) &:= \inf_{x\in
   V,\sigma(x,Jx)>0}\frac{\sigma(x,Jx)}{|x|^2} - \inf_{x\in
   V,\sigma(x,-Jx)>0}\frac{\sigma(x,-Jx)}{|x|^2}, \cr
   \beta(\sigma) &:= \sup_{\dim W=2}\frac{\sigma|_W}{\Om_W} =
   \sup_{x,y}\frac{\sigma(x,y)}{\sqrt{|x|^2|y|^2-\la x,y\ra^2}}, \cr
   \gamma(\sigma) &:= \frac{\alpha(\sigma)}{\beta(\sigma)},
\end{align*}
where the infima are interpreted as zero over the empty set, the
first supremum is taken over all oriented 2-planes $W\subset V$, and
the second one over all ordered pairs of linearly independent
vectors $x,y\in V$. Note that $\alpha,\beta,\gamma$ depend on $J$
and the Hermitian metric. The following properties follow directly
from the definitions:
\begin{enumerate}
\item For each $\sigma$, at most one of the two terms in the definition
of $\alpha(\sigma)$ is nonzero. More precisely, the following
three cases can occur: If $\sigma(x,Jx)\geq 0$ for all $x$ the
second term vanishes; if $\sigma(x,Jx)\leq 0$ for all $x$ the
first term vanishes; and in the remaining case both terms vanish
(for if $\sigma(x,Jx)>0$ and $\sigma(y,Jy)<0$, then
$\sigma\bigl(x+ty,J(x+ty)\bigr)=0$ for some $t>0$).
\item $\beta(\sigma)>0$ for every nonzero $\sigma$ and
$|\alpha(\sigma)|\leq\beta(\sigma)$ (just take $y=Jx$ in the
definition of $\beta$), thus $|\gamma(\sigma)|\leq 1$.
Equality $\gamma(\sigma)=\pm 1$ holds iff
$\sigma = \pm c\omega$ for some $c>0$.
\item $\gamma(\sigma)>0$ iff $\sigma$ is a symplectic form taming $J$ and
$\gamma(\sigma)<0$ iff $\sigma$ is a symplectic form taming $-J$.
\item $\gamma(t\sigma)=\gamma(\sigma)$ for every $t>0$ and
$\gamma(-\sigma)=-\gamma(\sigma)$. Moreover,
$\gamma(\sigma+\sigma')\geq\min\{\gamma(\sigma),\gamma(\sigma')\}$
if $\gamma(\sigma)$ and $\gamma(\sigma')$ are both positive.
\end{enumerate}

Here the last property follows from
$\beta(\sigma+\sigma')\leq\beta(\sigma)+\beta(\sigma')$ and
$\alpha(\sigma+\sigma')\geq\alpha(\sigma)+\alpha(\sigma') \geq
\min\{\gamma(\sigma),\gamma(\sigma')\}
[\beta(\sigma)+\beta(\sigma')]$. In view of $\gamma(\pm\om)=\pm 1$,
$\cos^{-1}\bigl(\gamma(\sigma)\bigr)\in[0,\pi]$ can be interpreted
as the ``angle'' between $\om$ and $\sigma$.

\begin{lemma}\label{lem:angle2}
(a) If $\sigma$ is a symplectic form on $V$ and $W\subset V$ a codimension
2 $\om$-symplectic subspace whose K\"ahler angle $\theta(W)$ (with
respect to $\om$) satisfies
$$
   \theta(W) < \gamma(\sigma),
$$
then $W$ and $W^\perp$ are $\sigma$-symplectic.

(b) If $\sigma$ is a symplectic form taming J
and $K:V\to V$ an almost complex structure such that
$$
   \|K-J\| < \gamma(\sigma),
$$
then $\sigma$ tames $K$.

(c) For every codimension two $\om$-symplectic subspace $W\subset V$
there exists an almost complex structure
$K:V\to V$ which is compatible with $\om$, leaves $W$ and $W^\om$
invariant, and satisfies
$$
   \|K-J\| \leq \theta(W).
$$
\end{lemma}

\begin{remark}
In fact, in (a) the actual hypothesis is
$\sin\theta(W)<\gamma(\sigma)$ and in (c) the actual conclusion is
$\|K-J\| \leq 2\sin\left(\theta(W)/2\right)$.
\end{remark}

\begin{proof}
(a) For $0\neq x\in W$ write
$$
   Jx = \pi_W(Jx) + y
$$
with $y\perp W$. By Lemma~\ref{lem:angle} (d) we have
$\angle_M(W,Jx)\leq\theta(W)$, hence
$|\pi_W(Jx)|\geq|x|\cos\theta(W)$, or equivalently,
$|y|\leq|x|\sin\theta(W)$. Then
\begin{align*}
   \sigma\bigl(x,\pi_W(Jx)\bigr) &= \sigma(x,Jx) - \sigma(x,y) \cr
   &\geq \alpha(\sigma)|x|^2 - \beta(\sigma)|x|\,|y| \cr
   &\geq |x|^2\Bigl(\alpha(\sigma)-\beta(\sigma)\sin\theta(W)\Bigr)
   \cr
   &> 0
\end{align*}
by hypothesis and $\sin\theta\leq\theta$. Since $\pi_W(Jx)\in W$,
this proves that $W$ is $\sigma$-symplectic. The statement for
$W^\perp$ follows by the same argument because
$\theta(W^\perp)=\theta(W)$ by Lemma~\ref{lem:angle} (e).

(b) By hypothesis, we find for every $0\neq x\in V$:
\begin{align*}
   \sigma(x,Kx) &= \sigma(x,Jx) + \sigma(x,(K-J)x) \geq
   \alpha(\sigma)|x|^2 - \beta(\sigma)\|K-J\|\,|x|^2 \cr
   &= \beta(\sigma)\Bigl(\gamma(\sigma)-\|K-J\|\Bigr)|x|^2 > 0.
\end{align*}

(c) Let $W_0:=W\cap JW$. If $\dim W_0=2n-2$, then $W$ is $J$-complex
and we simply set $K:=J$. If $\dim W_0=2n-4$, we choose $K$ such that
$K=J$ on $W_0$ and $K$ leaves $W_0^\om=W_0^\perp$ invariant. After
replacing $V$ by $W_0^\perp$, it thus remains to consider the case
$\dim V=4$.

For $\dim V=4$ pick an oriented orthonormal basis $x,y$ of $W$ and
denote by $\pi_W:V\to W$ the orthogonal projection. Since $\la
Jx,x\ra=0$, we have $\pi_WJx=\cos\theta\,y$ for some
$\theta\in[0,\pi]$. Now $\cos\theta=\la Jx,y\ra=\om(x,y)>0$ implies
$\theta\leq\pi/2$, and it follows from Lemma~\ref{lem:angle} (d) that
$\theta\leq\theta(W)$. Define $K:W\to W$ as positive rotation by $90$
degrees, i.e.~$Kx:=y$ and $Ky:=-x$. Then clearly $K^2=-\id$. Moreover,
$$
   |Jy-Ky|^2 = |Jx-Kx|^2 = |Jx-y|^2 = 2-2\cos\theta =
   4\sin^2\frac{\theta}{2}
$$
and $\theta\leq\theta(W)$ shows $\|K-J\|\leq
2\sin\left(\theta(W)/2\right)$. Finally, $\om(x,Kx)=\om(x,y)>0$ shows
that $\om$ tames $K$ and $K$-invariance of $\om$ is automatic because
$W$ is 2-dimensional, so $K:W\to W$ is compatible with $\om$. Defining
$K:W^\om\to W^\om$ in the analogous way we obtained the desired
$K:V\to V$.
This concludes the proof of Lemma~\ref{lem:angle2}. 
\end{proof}

For an $\R$-linear map $T:V\to V'$ between Euclidean vector spaces
define $\nu(T)\geq 0$ as follows: If $T$ is not surjective set
$\nu(T):=0$; otherwise let $\nu(T)^{-1}$ be the minimal norm of a
right inverse of $T$. Define the {\em minimal angle} between two
subspaces $X,Y\subset V$ of a Euclidean vector space as follows: If
$X$ and $Y$ are not transverse set $\angle_m(X,Y):=0$; otherwise let
$X',Y'$ be the orthogonal complements of $X\cap Y$ in $X,Y$ and set
$$
   \angle_m(X,Y) := \inf_{0\neq y\in Y'}\angle(X,y)
   = \inf_{0\neq x\in X',0\neq y\in Y'}\angle(x,y).
$$

\begin{lemma}\label{lem:angle3}
(a) Let $T:V\to V'$ be a nonzero $\R$-linear map between Euclidean
vector spaces and let $W\subset V$ be a linear subspace. Then 
$$ 
   \angle_m(W,\ker T) \geq \frac{\nu(T|_W)}{\|T\|}. 
$$  
(b) For every $\eps>0$ and $\theta_2>0$ there exists a $\theta_3>0$
with the following property: For any pair of codimension two
$\om$-symplectic subspaces $W,W'$ of a Hermitian vector space
$(V,\om,J)$ with $\angle_m(W,W')\geq\eps$ and
$\max\{\theta(W),\theta(W')\}\leq\theta_3$ 
there exists an almost complex structure $K:V\to V$ which is
compatible with $\om$, leaves $W$ and $W'$ invariant, and satisfies
$$
   \|K-J\| < \theta_2.
$$
In particular, $W$ and $W'$ intersect positively. 
\end{lemma}

\begin{proof}
(a) Suppose that $T|_W$ is surjective (otherwise both sides are zero).
Set $\eta:=\nu(T|_W)$ and let $W'\subset W$ be the orthogonal
complement of $W\cap\ker T$. Consider $w\in W'$ and $v\in\ker T$ with
$|v|=|w|=1$. After possibly replacing $v$ by $-v$, we may assume $\la
v,w\ra\geq 0$. The definition of $\nu(T|_W)$ implies $|Tw|\geq\eta$,
and therefore 
$$
   \eta^2\leq |Tv-Tw|^2 \leq \|T\|^2|v-w|^2 = 2\|T\|^2(1-|\la
   v,w\ra|). 
$$
So the angle $\theta:=\angle(v,w)\in[0,\pi/2]$ satisfies
$$
   1-\frac{\eta^2}{2\|T\|^2} \geq \cos\theta \geq
   1-\frac{\theta^2}{2}, 
$$
thus $\theta\geq\eta/\|T\|$ and the Part (a) follows.  

(b) is obvious if we take $\theta_3=0$, and for small $\theta_3>0$ it
follows by a continuity argument. More precisely, note first that for
every codimension two $\om$-symplectic 
subspace $W\subset V$ there exists a (non-unique) codimension two
$J$-complex subspace $W_c$ with $\angle_M(W,W_c)\leq\theta(W)$. (In the
notation of the proof of Lemma~\ref{lem:angle2} (c) we can take for
$W_c$ the subspace generated by 
$W_0$, $x$ and $Jx$). Pick such complex subspaces $W_c,W_c'$ for
$W,W'$ so that $\angle_M(W,W_c)\leq\theta_3$ and
$\angle_M(W',W_c')\leq\theta_3$. By choosing $\theta_3$ small we can
achieve that $\angle_m(W_c,W_c')\geq\eps/2$
and $\angle_M(W\cap W',W_c\cap W_c')$ is arbitrarily small. Pick a
complex basis $v_1,v_2=Jv_1,\dots,v_{2n-1},v_{2n}=Jv_{2n-1}$ of $V$
such that $v_1,\dots,v_{2n-4}$ is a unitary basis of $W_c\cap W_c'$, 
$v_1,\dots,v_{2n-2}$ is a unitary basis of $W_c$, and 
$v_1,\dots,v_{2n-4},v_{2n-1},v_{2n}$ is a unitary basis of $W_c'$.
This choice and the condition $\angle_m(W_c,W_c')\geq\eps/2$ ensure
that $\angle(v_i,v_j)\geq\eps/2$ for all $i\neq j$. 
Denote by $w_1,\dots,w_{2n}$ the orthogonal projections of
$v_1,\dots,v_{2n-4}$ onto $W\cap W'$, $v_{2n-3},v_{2n-2}$ onto $W$ and
$v_{2n-1},v_{2n}$ onto $W'$. Then $w_1,\dots w_{2n}$ is a symplectic
basis close to $v_1,\dots v_{2n}$, so the linear map $\phi$ sending
$v_i$ to $w_i$ is symplectic and close to the identity. It follows
that $K:=\phi J\phi^{-1}$ is an almost complex structure compatible
with $\om$ which leaves $W$ and $W'$ invariant. Moreover, $K$ is
arbitrarily close to $J$, so we can achieve $\|K-J\| \leq \theta_2$. 
This concludes the proof of Lemma~\ref{lem:angle3}. 
\end{proof}

Recall from~\cite{MS-intro} that any pair of $\om$-compatible complex
structures $J_0,J_1$ can be connected by a canonical path $J_t$ of
almost complex structures. Moreover, the construction
in~\cite{MS-intro} shows

\begin{lemma}\label{lem:contraction}
For every $\theta_0>0$ there exist $\theta_1>\theta_2>0$ such that
$\|J_1-J_0\|<\theta_1$ implies $\|J_t-J_0\|<\theta_0$ for all
$t\in[0,1]$, and $\|J_1-J_0\|<\theta_2$ implies $\|J_t-J_0\|<\theta_1$
for all $t\in[0,1]$.
\end{lemma}

{\bf A priori estimates. }
Consider now a closed symplectic manifold $(X,\om)$ with a fixed
compatible almost complex structure $J$. For a codimension 2
submanifold $Y\subset X$ define
\begin{align*}
   \theta(Y;\om,J) := \sup_{y\in Y}\theta(T_yY;\om,J),
\end{align*}
where $\theta(T_yY;\om,J)$ is the K\"ahler angle with respect to the
Hermitian structure $(\om,J)$ on $TX$.
Denote by $\JJ(X)$ the space of almost complex
structures on $X$ compatible with $\om$. Fix constants
$$
   0<\theta_2<\theta_1<\theta_0<1
$$
as in Lemma~\ref{lem:contraction}.

Now assume that $\om$ represents an integral cohomology class
$[\om]\in H^2(X;\Z)$. By Theorem~\ref{thm:Donaldson}, for every
sufficiently large integer $D>0$ there exists a closed submanifold
$Y\subset X$ of codimension 2 which is Poincar\'e dual to
$D[\om]$ and approximately $J$-holomorphic in the sense that
$$
   \theta(Y;\om,J)<\theta_2.
$$
In view of the
hypothesis $\theta_2<1<\pi/2$, $Y$ is symplectic for $\om$.

By Lemma~\ref{lem:angle2} (c) there exists an $\om$-compatible
almost complex structure $K$ on the vector bundle $TX|_Y$ with
$$
   \|K-J\|_Y<\theta_2,\qquad K(TY)=TY,\qquad K(TY^\om)=TY^\om.
$$
By Lemma~\ref{lem:contraction} and the choice of
$\theta_2<\theta_1$, we can extend $K$ to an $\om$-compatible
almost complex structure on $X$ satisfying
$$
   \|K-J\|<\theta_1.
$$
Indeed, first extend $K$ to some $K_1$ on $\{x\in X\mid
dist(x,Y)\le r\}$ in the metric $\om(.,J.)$ such that
$\|K_1-J\|<\theta_2$, which is possible because this is an open
condition. By the lemma there is a path $K_t$ on that tube
connecting $K_1$ to $K_0=J$ and satisfying
$\|K_t-J\|<\theta_1$ for all $t$. Finally, pick  a smooth function
$\phi:[0,r]\to 1$ which is equal to $0$ near $r$ and equal to $1$ near
$0$ and set $K_x:= K_{\phi(dist(x,Y))}$ for $x\in X$.
The upshot of this discussion is that the space
$$
   \JJ(X,Y) := \{K\in\JJ(X)\mid\ K(TY)=TY\}
$$
contains a nonempty open subspace
$$
   \JJ(X,Y;J,\theta_1) := \{K\in\JJ(X,Y)\mid\ \|K-J\|<\theta_1\}.
$$
Moreover, by the preceding argument and Lemma~\ref{lem:contraction},
any two elements in $\JJ(X,Y;J,\theta_1)$ can be connected by a path in 
$$
   \JJ(X,Y;J,\theta_0) := \{K\in\JJ(X,Y)\mid\ \|K-J\|<\theta_0\}.
$$

The following lemma and proposition follow a suggestion of D.~Auroux.
Fix a closed 2-form $\alpha$ representing the first Chern class
$c_1=c_1(TX)$ and denote by $\|\alpha\|$ its norm with respect to
the metric $\om(\cdot,J\cdot)$.

\begin{lemma}\label{lem:chern}
Let $K$ be an $\om$-compatible almost complex structure with
$\|K-J\|<\theta_0$. Then for every homology class $A\in H_2(X;\Z)$
which contains a closed $K$-holomorphic curve we have
$$
   \la c_1(TX),A\ra\leq D_*\om(A),\qquad\text{with
   }D_*=D_*(X,\om,J):=\frac{1+\theta_0}{1-\theta_0}\|\alpha\|.
$$
\end{lemma}

\begin{proof}
For every tangent vector $v\in T_xX$ we have $\|v\|^2=\om_0(v,Jv)$
and the estimates
$$
   \alpha(v,Kv) \leq \|\alpha\|(1+\|K-J\|)\|v\|^2,\qquad
   \om(v,Kv) \geq (1-\|K-J\|)\|v\|^2,
$$
which combine with $\|K-J\|<\theta_0$ to
$$
   \alpha(v,Kv)\leq D_*\om(v,Kv).
$$
For a closed $K$-holomorphic curve $f:\Sigma\to X$ in the class $A$
this implies
$$
   \la c_1(TX),A\ra = \int_\Sigma f^*\alpha \leq D_*\int_\Sigma
   f^*\om = D_*\om(A).
$$
\end{proof}

\begin{remark}
For a given cohomology class $c\in H^2(X,\R)$, one may wonder how much
the quantity  
$$
   \|c\|_{g_J}:=\inf\{\|\alpha\|_{g_J}\mid d\alpha=0,\ [\alpha]=c\}
$$ 
with respect to the metric $g_J=\om(\cdot,J\cdot)$ depends on the
compatible almost complex structure $J\in\JJ(X,\om)$. For example, the
quantity 
$$
   \sup_{J\in\JJ(X,\om)}\|c\|_{g_J}
$$ 
equals $1$ for $c=[\om]$. Is it finite for $c=c_1(X)$? 
\end{remark}

In the following the {\em energy} of a $K$-holomorphic sphere
$f:S^2\to X$ will always be measured with respect to $\om$,
i.e.~$E(f)=\int_{S^2}f^*\om$.

\begin{prop}\label{prop:tangency}
Let $Y\subset X$ be a hypersurface Poincar\'e dual to $D[\om]$ with
$\theta(Y;\om,J)<\theta_2$ and let $K\in\JJ(X,Y;J,\theta_0)$. Let
$E>0$ and suppose: 

(i) All moduli spaces of simple $K$-holomorphic spheres in $Y$ of
energy $\leq E$ are smooth manifolds of the expected dimension.

(ii) All moduli spaces of nonconstant simple $K$-holomorphic spheres
of energy $\leq E$ in $X$ with prescribed tangency of order $\ell\leq
D_*E+n$ to $Y$ are smooth manifolds of the expected dimension.

Then the following holds:

(a) If $D>\max(D_*,D_*+n-4)$, then all $K$-holomorphic spheres of
energy $\leq E$ contained in $Y$ are constant.

(b) If $D>2\max(D_*,D_*+n-2)$, then every nonconstant $K$-holomorphic
sphere of energy $\leq E$ in $X$ intersects $Y$ in at least 3 distinct
points in the domain.
\end{prop}

\begin{proof}
(a) Suppose $f$ is a nonconstant $K$-holomorphic sphere of energy
$\leq E$ contained in $Y$. After replacing $f$ by its underlying
simple curve, we may assume that $f$ is simple. Let $A\in H_2(Y;\Z)$
be its homology class. The expected dimension of the moduli space of
simple $K$-holomorphic spheres in $Y$ in the homology class $A$ equals
$$
   \ind(A) = 2(n-1)-6+2\la c_1(TY),A\ra.
$$
Now note that
$$
   c_1(TY)=c_1(TX|_Y)-c_1\Bigl(N(Y,X)\Bigr),
$$
and the first Chern class of the normal bundle $N(Y,X)$ is
represented by $D\om|_Y$. In view of Lemma~\ref{lem:chern}, this
implies
\begin{align*}
   \ind(A) &= 2n-8+2\la c_1(TX),A\ra - 2D\om(A) \cr
   &\leq 2n-8+2(D_*-D)\om(A).
\end{align*}
Now $\om(A)\geq 1$ due to the integrality of $[\om]$ and the
assumption that $f$ is nonconstant. Hence $\ind(A)<0$ if
$D>\max(D_*,D_*+n-4)$, so by hypothesis (i) $f$ cannot exist.

(b) Suppose $f$ is a nonconstant $K$-holomorphic sphere of energy
$\leq E$ in $X$ which intersects $Y$ in at most 2 distinct points in
the domain. Again
after replacing $f$ by its underlying simple curve, we may assume that
$f$ is simple. Let $A\in H_2(Y;\Z)$ be its homology class.
For $I\leq D_*E+n+1$ let $\tilde\MM^s(Y,I;A,K)$ be the space of all simple
$K$-holomorphic spheres in the class $A$ that have local intersection
number $\geq I$ (i.e.~a tangency of order $\geq I-1$) with $Y$ at one
point. By hypothesis (ii) and  Proposition~\ref{prop:jets}
$\tilde\MM^*(Y,I;A,K)$ is a smooth manifold of dimension
$$
   \dim\tilde\MM^*(Y,I;A,K) = 2n - 4 + 2c_1(A) -
   2I \geq 0.
$$
Hence the local intersection number of $f$ with $Y$ at each point
satisfies
$$
   I\leq c_1(A)+n-2\leq D_*\om(A)+n-2.
$$
On the other hand, if $D>2\max(D_*,D_*+n-2)$, then the global
intersection number of $f$ with $Y$ satisfies
$$
   [f]\cdot[Y] = D\om(A) > 
   2\bigl(D_*\om(A)+n-2\bigr)\geq 2I
$$
(where we have used again $\om(A)\geq 1$). Hence $f$ must have at
least 3 intersection points with $Y$, which contradicts our
assumption.
\end{proof}

We will also need the following version of the previous proposition
for families of almost complex structures. As above let $Y\subset X$
be a hypersurface Poincar\'e dual to $D[\om]$ with
$\theta(Y;\om,J)<\theta_2$. Let $\KK\subset\JJ(X,Y;J,\theta_0)$ be a
family of almost complex structures smoothly depending on a parameter
$\tau\in P$ varying in a $k$-dimensional manifold $P$. Moduli spaces
of $\KK$--holomorphic spheres are 
spaces of pairs $(u,\tau)$ where $\tau\in P$ and $u$ is a
$\KK_\tau$-holomorphic sphere (possibly with constraints). The
corresponding transformation groups act on $u$ only. 
The (expected) dimension of the moduli spaces is increased by $k$
compared to the expected dimension of the corresponding moduli spaces
with fixed almost complex structure. The proof of the following
proposition is literally the same as the one of
Proposition~\ref{prop:tangency}, each time adding $k$ to the term
representing the dimension of the moduli space. 

\begin{prop}\label{prop:tang-family}
Let $E>0$ and suppose: 

(i) All moduli spaces of simple $\KK$-holomorphic spheres in $Y$ of
energy $\leq E$ are smooth manifolds of the expected dimension.

(ii) All moduli spaces of nonconstant simple $\KK$-holomorphic spheres
of energy $\leq E$ in $X$ with prescribed tangency of order $\ell\leq
D_*E+n+k$ to $Y$ are smooth manifolds of the expected dimension.

Then the following holds:

(a) If $2D>\max(2D_*,2D_*+2n-8+k)$, then all $\KK$-holomorphic spheres of
energy $\leq E$ contained in $Y$ are constant.

(b) If $D>2\max(D_*,D_*+n-2)+k$, then every nonconstant
$\KK$-holomorphic sphere of energy $\leq E$ in $X$ intersects $Y$ in
at least 3 distinct points in the domain.
\end{prop}

\begin{definition}\label{def:JJ*}
For $(J,Y)$ as above and $E>0$ let
$$
   \JJ^*(X,Y;J,\theta_1,E)\subset \JJ(X,Y;J,\theta_1)
$$
be the subset of those $K$ for which the conclusions (a) and (b) of
Proposition~\ref{prop:tangency} hold, i.e.

(a) All $K$-holomorphic spheres of
energy $\leq E$ contained in $Y$ are constant.

(b) Every nonconstant $K$-holomorphic sphere of
energy $\leq E$ in $X$ intersects $Y$ in at least 3 distinct points in
the domain.
\end{definition}

Define
$$
   D^*:=D^*(W,\om,J):=2(D_*+n). 
$$
Note that each degree $D\geq D^*$ satisfies the conditions in (a) and
(b) of Proposition~\ref{prop:tangency}.

\begin{cor}\label{cor:open-dense}
Let $(J,Y)$ be as in Proposition~\ref{prop:tangency} and suppose
$D\!\!\geq\!\! D^*\!(X,\!\om,\!J)$. Then 
$\JJ^*(X,Y;J,\theta_1,E)$ is open and dense in $\JJ(X,Y;J,\theta_1)$
for every $E>0$. Moreover, for $\theta_2<\theta_1$ any two elements in
$\JJ^*(X,Y;J,\theta_2,E)$ can be connected by a path in
$\JJ^*(X,Y;J,\theta_1,E)$.
\end{cor}

\begin{proof}
Condition (a) is open by Gromov compactness. Indeed, if (a) were not
open at some $K\in\JJ^*(X,Y;J,\theta_1,E)$ there would exist a sequence
$K^\nu\to K$ in 
$\JJ(X,Y;J,\theta_1)$  with nonconstant $K^\nu$-holomorphic
spheres $f^\nu$ of energy $\leq E$ contained in $Y$; Gromov
compactness would yield a nonconstant $K$-holomorphic sphere of energy
$\leq E$ contained in $Y$, contradicting
$K\in\JJ^*(X,Y;J,\theta_1,E)$. Openness of condition (b) is proved
similarly: Suppose by contradiction that there is a sequence of
$K^\nu\in\JJ(X,Y;J,\theta_1)$ not satisfying (b) and converging to
$K\in\JJ^*(X,Y;J,\theta_1,E)$. Then there
exists a sequence of nonconstant $K^\nu$-holomorphic spheres
$f^\nu:S^\nu\to X$ of energy $\leq E$ intersecting $Y$ in at most 2
points. Then Gromov compactness yields a nonconstant $K$-holomorphic
sphere of energy $\leq E$ intersecting $Y$ in at most 2
points, contradicting $K\in\JJ^*(X,Y;J,\theta_1,E)$.

For density of $\JJ^*(X,Y;J,\theta_1,E)$, observe that the set of
$K\in\JJ(X,Y;J,\theta_1)$ satisfying conditions (i) and (ii) of
Proposition~\ref{prop:tangency} is dense: For condition (i) this
is Theorem 3.1.5 in~\cite{MS} (cf.~the proof of
Lemma~\ref{lem:surj-nonconst} above). Density of condition (ii)
follows from Proposition~\ref{prop:jets} (or
Lemma~\ref{lem:jets}) with $k=1$, $Z_1:=Y$ and
$V:=X\setminus Y$. (Here we assume that (i) already holds, hence by
conclusion (a) every nonconstant $K$-holomorphic sphere of energy
$\leq E$ passes through $V$.)

For the last statement let $J_0,J_1\in \JJ^*(X,Y;J,\theta_2,E)$.
By Lemma~\ref{lem:contraction} they can be connected by a smooth  
path $(J_\tau)_{\tau\in[0,1]}$ in $\JJ(X,Y;J,\theta_1,E)$. By 
arguments similar to the previous ones, conditions (i) and
(ii) of Proposition~\ref{prop:tang-family} can be achieved for a path
$J_\tau'$ arbitrarily close to $J_\tau$. In particular $J'(0)$ and
$J'(1)$ are arbitrarily close to $J_0$ and $J_1$, respectively. Since
$\JJ^*(X,Y;J,\theta_1,E)$ is open, that means that we may assume that
$J_0$ and $J_1$ can be connected in this set.   

This concludes the proof.
\end{proof}

{\bf Extension to two hypersurfaces. }
The preceding discussion will suffice for the proof of
Theorem~\ref{thm:moduli} and~\ref{thm:pseudo}. For the proof of
Theorem~\ref{thm:indep} we need to extend it to the case of two
symplectic hypersurfaces. As most arguments immediately carry over, we
will only indicate the required modifications. 

As above, consider a closed symplectic manifold $(X,\om)$ with
integral class $[\om]$ 
Fix constants
$$
   0<\theta_3<\theta_2<\theta_1<\theta_0<1
$$
as in Lemma~\ref{lem:contraction} and Lemma~\ref{lem:angle3} (b) and a 
constant $\eps>0$. For the remainder of this section, we fix a
quadruple $(J_0,Y_0,J_1,Y_1)$ consisting of $\om$-compatible almost
complex structures $J_i$ and hypersurfaces $Y_i$, Poincar\'e dual to
$D_i[\om]$, satisfying 
$$
   \|J_0-J_1\|_0<\theta_3,\qquad \angle_m(Y_0,Y_1)\geq\eps,\qquad
   \theta(Y_i)<\theta_3 \text{ for } i=0,1.  
$$
Here $\|\ \|_0$ denotes the norm with respect to
$\om(\cdot,J_0\cdot)$ and we have set
$$
   \angle_m(Y_0,Y_1) := \inf_{y\in Y_0\cap
   Y_1}\angle_m(T_yY_0,T_yY_1). .
$$ 
\begin{remark}
The existence of such quadruples follows from the (Stabilization)
property in Theorem~\ref{thm:Donaldson}: Given a $J_0$-holomorphic
hypersurface $Y_0$, there exists a hypersurface $Y_1$ such that the
above conditions are satisfied with $J_0=J_1$. We allow the case
$J_0\neq J_1$ for later use. 
\end{remark}
By Lemma~\ref{lem:angle3} (b) there exists an $\om$-compatible
almost complex structure $K$ on the vector bundle $TX|_{Y_0\cap Y_1}$
with 
$$
   \|K-J_0\|_{Y_0\cap Y_1}<\theta_2,\qquad K(TY_0)=TY_0,\qquad
   K(TY_1)=TY_1. 
$$
Arguing as above, using Lemma~\ref{lem:contraction}, we can extend $K$
to an $\om$-compatible almost complex structure on $X$ satisfying
$$
   \|K-J_0\|<\theta_1,\qquad K(TY_0)=TY_0,\qquad
   J(TY_1)=TY_1. 
$$
Thus the space
$$
   \JJ(X,Y_0\cup Y_1) := \{K\in\JJ(X)\mid\ K(TY_0)=TY_0, K(TY_1)=TY_1,
   \} 
$$
contains a nonempty open subspace
$$
   \JJ(X,Y_0\cup Y_1;J_0,\theta_1) := \{K\in\JJ(X,Y_0\cup Y_1)\mid\
   \|K-J_0\|<\theta_1\}. 
$$
Again, any two elements in $\JJ(X,Y_0\cup Y_1;J_0,\theta_1)$ can be
connected in the space 
$$
   \JJ(X,Y_0\cup Y_1;J_0,\theta_0) := \{K\in\JJ(X,Y_0\cup Y_1)\mid\
   \|K-J_0\|<\theta_0\}. 
$$
We need the following variation of Proposition~\ref{prop:tangency}. 

\begin{lemma}\label{lem:tangency2}
Let $(J_0,Y_0,J_1,Y_1)$ be as above, $K\in\JJ(X,Y_0\cup
Y_1;J_0,\theta_0)$, and $E>0$. Suppose that all moduli spaces of
simple $K$-holomorphic spheres in the intersection $Y_0\cap Y_1$ of
energy $\leq E$ are smooth manifolds of the expected dimension
and $D_0+D_1>\max(D_*,D_*+n-5)$. Then all $K$-holomorphic
spheres of energy $\leq E$ contained in $Y_0\cap Y_1$ are constant.
\end{lemma}

\begin{proof}
We have
$$
   c_1\bigl(T(Y_0\cap Y_1)\bigr) = c_1(TX|_{Y_0\cap Y_1}) -
   D_0[\om] - D_1[\om]. 
$$
Hence the expected dimension of the moduli space of simple
$K$-holomorphic spheres in $Y$ in a homology class $A$ satisfies
\begin{align*}
   \ind(A) &= 2(n-2)-6+2\la c_1\bigl(T(Y_0\cap Y_1)\bigr),A\ra \cr
   &= 2n-10+2\la c_1(TX),A\ra - 2(D_0+D_1)\om(A) \cr
   &\leq 2n-10+2(D_*-D_0-D_1)\om(A).
\end{align*}
If $A$ is nontrivial, then $\om(A)\geq 1$ and the hypothesis
$D_0+D_1>\max(D_*,D_*+n-5)$ implies $\ind(A)<0$. In view of the
regularity hypothesis (and again passing to underlying simple
spheres), this shows that all $K$-holomorphic spheres of energy $\leq
E$ contained in the intersection $Y_0\cap Y_1$ are constant. 
\end{proof}

\begin{definition}\label{def:JJ*2}
For $E>0$ let
$$
   \JJ^*(X,Y_0\cup Y_1;J_0,\theta_1,E)\subset \JJ(X,Y_0\cup
   Y_1;J_0,\theta_1) 
$$
be the subset of those $K$ for which the conclusions (a) and (b) of
Proposition~\ref{prop:tangency} hold for each $Y_i$, i.e.

(a) All $K$-holomorphic spheres of
energy $\leq E$ contained in $Y_0\cup Y_1$ are constant.

(b) Every nonconstant $K$-holomorphic sphere of energy $\leq E$ in $X$
intersects $Y_i$ in at least 3 distinct points in the domain for
$i=0,1$. 
\end{definition}

Note that if $D_i\geq D^*=2(D_*+n)$ for $i=0,1$, then each $D_i$
satisfies the conditions in (a) and (b) of 
Proposition~\ref{prop:tangency} and $D_0+D_1$ satisfies the hypothesis
of Lemma~\ref{lem:tangency2}.
We have the following analogue of Corollary~\ref{cor:open-dense}. 

\begin{cor}\label{cor:open-dense2}
Let $(J_0,Y_0,J_1,Y_1))$ be as above and suppose $D_i\!\geq\!
D^*\!(X,\om,\!J)$ for $i=0,1$. Then 
$\JJ^*(X,Y_0\cup Y_1;J_0,\theta_1,E)$ is open and dense in
$\JJ(X,Y_0\cup Y_1;J_0,\theta_1)$ for every $E>0$.  
Moreover, for $\theta_2<\theta_1$ any two elements in
$\JJ^*(X,Y_0\cup Y_1;J_0,\theta_2,E)$ can be connected by a path in 
$\JJ^*(X,Y_0\cup Y_1;J_0,\theta_1,E)$.
\end{cor}

\begin{proof}
Openness follows from Gromov compactness as in the proof of
Corollary~\ref{cor:open-dense}. For density, it suffices to prove
density of the set of $K\in\JJ(X,Y_0\cup Y_1;J_0,\theta_1)$ satisfying
conditions (i) and (ii) of Proposition~\ref{prop:tangency} for both
hypersurfaces $Y_i$. For this, we first apply Theorem 3.1.5
of~\cite{MS} in $Y_0\cap Y_1$ to conclude density of the set of
$K\in\JJ(X,Y_0\cup Y_1;J_0,\theta_1)$ for 
which all moduli spaces of simple $K$-holomorphic spheres in $Y_0\cap
Y_1$ of energy $\leq E$ are smooth manifolds of the expected dimension. 
Then by Lemma~\ref{lem:tangency2} all $K$-holomorphic spheres
contained in $Y_0\cap Y_1$ are constant. Based on this, we can apply
Theorem 3.1.5 of~\cite{MS} in $Y_0$ and $Y_1$ (varying $K$ on $Y_0\cup
Y_1\setminus(Y_0\cap Y_1$)) to conclude density of the set of
$K\in\JJ(X,Y_0\cup Y_1;J_0,\theta_1)$ satisfying condition (i) of
Proposition~\ref{prop:tangency} for both $Y_i$. Finally, we apply 
Proposition~\ref{prop:jets} (or Lemma~\ref{lem:jets}) with $k=1$,
$Z_1:=Y_i$ and $V:=X\setminus(Y_0\cup Y_1)$ (varying $K$ outside
$Y_0\cup Y_1$) to conclude density of condition (ii).
This finishes the proof. 
\end{proof}

{\bf Addendum. }
We conclude this section with a finiteness result for homology classes
represented by $J$-holomorphic curves, although this result will not
be needed in this paper. For a compatible almost complex structure $J$
and $0<\theta_0<1$ set
$$
   \JJ(X;J,\theta_0) := \{K\in\JJ(X)\mid\ \|K-J\|<\theta_0\}.
$$

\begin{lemma}\label{lem:AA}
For each $E>0$ and $0<\theta_0<1$, all homology classes $A\in
H_2(X;\Z)$ with $\la\om,A\ra\leq E$ that can be represented by a
$K$-holomorphic sphere for some $K\in\JJ(X;J,\theta_0)$ are contained in a
finite set $\AA(E,\theta_0)\subset H_2(X;\Z)$ depending only on
$\om$, $J$, $E$ and $\theta_0$.
\end{lemma}

\begin{proof}
For a 2-form $\sigma$ on $X$ define
$$
   \gamma(\sigma) := \inf_{x\in X}\gamma(\sigma_x),
$$
where $\gamma(\sigma_x)$ is defined as above with
respect to $(\om,J)$ and the induced metric on $TX$. Fix a constant
$\gamma_0$ satisfying
$$
   0<\theta_0<\gamma_0<1
$$
and define
$$
   \Om := \{\sigma\in\Om^2(X)\mid d\sigma=0,\ \gamma(\sigma)>\gamma_0\}.
$$
By the properties of $\gamma$, $\Om$ is a cone consisting of
symplectic forms taming $J$, and $\om\in\Om$. Hence its image
$[\Om]\subset H^2(X;\R)$ in cohomology is an open cone containing
$[\om]$. Denote by
$$
   [\Om]^* := \{A\in H_2(X;\R)\mid \la\sigma,A\ra\geq 0 \text{ for all
   }\sigma\in\Om\}
$$
the dual cone in homology. Then for given $E>0$ the set
$\{A\in [\Om]^*\mid \la\om,A\ra\leq E\}$ is bounded. (This is a
general property of cones in finite dimensions which can be seen as
follows: Write $A=\sum_jA_jv_j$ in a basis of cohomology whose dual
basis satisfies $[\om]\pm v^*_j\in[\Om]$ for all $j$; then
$0\leq\la\om\pm v^*_j,A\ra\leq E\pm A_j$ implies $|A_j|\leq E$ for
all $j$). Therefore, the set
$$
   \AA:=\{A\in [\Om]^*\mid \la\om,A\ra\leq E\}\cap
   H_2(X;\Z)
$$
is finite. Since $\theta_0<\gamma_0$, Lemma~\ref{lem:angle2} (b)
and the definition of $\Om$ show that $\JJ(X;J,\theta_0)$ consists
of almost complex structures which are tamed by all $\sigma\in\Om$.
Now suppose that a class $A\in H_2(X;\Z)$ with $\la\om,A\ra\leq E$
can be represented by a $K$-holomorphic sphere $f:S^2\to X$ for
some $K\in\JJ(X;J,\theta_0)$.
Then $\la\sigma,A\ra = \int_{S^2}f^*\sigma\geq 0$ for all
$\sigma\in\Om$, which implies $A\in\AA$.
\end{proof}

\section{Proof of transversality}\label{sec:proof}

In this section we prove Theorem~\ref{thm:moduli} and
Theorem~\ref{thm:pseudo} from the introduction, as well as their
extensions to the case of two hypersurfaces. Throughout this section,
$(X,\om)$ is a closed symplectic manifold of dimension $2n$. 

{\bf Summary of transversality results. }
Before beginning the proofs, let us summarize the transversality
results proved in Sections~\ref{sec:nodal} and~\ref{sec:tangencies} in
the form we will use them. Here conditions (a) follow from
Corollaries~\ref{cor:surj-tree} and~\ref{cor:surj-tree2}, and
conditions (b) follow from Corollaries~\ref{cor:jets-tree}
and~\ref{cor:jets-tree2}.  

\begin{prop}\label{prop:trans}
Let $J_0\in\JJ$. Let $\ZZ$ be a countable collection of
$J_0$-complex manifolds not equal to $X$ and set
$V:=X\setminus\cup_{Z\in\ZZ}Z$. Then for all $\ell\geq 3$ there exist
Baire sets $\JJ_{\ell+1}^\reg(V;\ZZ)\subset\JJ_{\ell+1}(V)$ such that 
for each $I$-stable $k$-labelled tree $T$, homology classes
$A_\alpha\in H_2(X;\Z)$ for $\alpha\in T$, $Z_i\in\ZZ\cup\{X\}$ and
integers $\ell_i\geq -1$ for $i\in R$ (the reduced index set), and
$K\in\JJ_{|I|+1}^\reg(V;\ZZ)$ the following conditions are satisfied: 

(a) The moduli space $\MM^*_T(\{A_\alpha\},K;\{Z_i\})$ of stable 
$K$-holomorphic maps modelled over $T$, passing through $Z_i$ at
the points $z_i$, $i\in R$ and without a nonconstant component
entirely contained in a $Z\in\ZZ$, is a smooth manifold of dimension
$$
   \dim\MM^*_T(\{A_\alpha\},K;\{Z_i\}) = 2n -6 +
   \sum_{\alpha\in T}2c_1(A_\alpha) + 2k - 2e(T) - 2\sum_{i\in
   R}\codim_\C Z_i. 
$$
Moreover, the evaluation map $\ev^R$ factors as 
$$
   \ev^R:\MM^*_T(\{A_\alpha\},K;\{Z_i\})\to
   \MM^*_{\pi_R(T)}(\{A_\alpha\},K;\{Z_i\})\to \prod_{i\in R}Z_i 
$$
through a smooth manifold of dimension 
$$
   \dim\MM^*_{\pi_R(T)}(\{A_\alpha\},K;\{Z_i\}) = 2n -6 +
   \sum_{\alpha\in T}2c_1(A_\alpha) + 2|R| - 2e\bigl(\pi_R(T)\bigr) -
   2\sum_{i\in R}\codim_\C Z_i.  
$$
(b) Suppose that $A_\alpha\neq 0$ for a
unique vertex $\alpha$ (so the reduced index set $R$ labels 
the special points on $\alpha$). Then the moduli space
$\MM^*_T(A_\alpha,J;\{Z_i,\ell_i\})$ of (not necessarily
stable $J$-holomorphic maps modelled over $T$, 
tangent to $Z_i$ of order $\ell_i$ at the special points
$z_{\alpha i}$ on the component $\alpha$ and not entirely contained in 
a $Z\in\ZZ$, is a smooth manifold of dimension
$$
   \dim\MM^*_T(A_\alpha,J;\{Z_i,\ell_i\}) = 2n - 6 +
   2c_1(A_\alpha) + 2k - 2e(T) -
   2\sum_{i\in R}(\ell_i+1)\codim_\C Z_i.
$$
Moreover, the evaluation map $\ev^R$ factors as 
$$
   \ev^R:\MM^*_T(A_\alpha,J;\{Z_i,\ell_i\})\to
   \MM^*_{|R|}(A_\alpha,J;\{Z_i,\ell_i\})\to \prod_{i\in R}Z_i
$$
through a smooth manifold of dimension 
$$
   \dim\MM^*_{|R|}(A_\alpha,J;\{Z_i,\ell_i\}) = 2n - 6 +
   2c_1(A_\alpha) + 2|R| -
   2\sum_{i\in R}(\ell_i+1)\codim_\C Z_i.
$$
\end{prop}

{\bf Proof of Theorem~\ref{thm:moduli}. }
From now on we suppose that $[\om]\in H^2(X;\Z)$, and we have fixed
constants $\eps>0$ and $0<\theta_3<\theta_2<\theta_1<\theta_0<1$ as in
the previous section. Also recall the definition 
$$
   D^*=D^*(X,\om,J_0)=2(D_*+n),
$$
where $D_*=D_*(X,\om,J_0)$ is the constant from Lemma~\ref{lem:chern}. 

\begin{definition}\label{def:D-pair}
A {\em Donaldson pair} is a pair $(J,Y)$ consisting of an
$\om$-compatible almost complex structure $J$ and a hypersurface
$Y^{2n-2}\subset X$ Poincar\'e dual to $D[\om]$ satisfying 
$$
   \theta(Y;\om,J)<\theta_3,\qquad D\geq D^*(X,\om,J).
$$
\end{definition}

\begin{remark}
The proof of Theorems~\ref{thm:moduli} and~\ref{thm:pseudo} only uses
the weaker hypothesis $\theta(Y;\om,J)<\theta_2$. The stronger
hypothesis is needed in the proof of Theorem~\ref{thm:indep}. 
\end{remark}

Let $(J,Y)$ be a Donaldson pair of degree $D$ on $(X,\om)$. Fix an
integer $E>0$ such that  
$$
   \ell := DE \geq 3. 
$$
Let $\JJ^*(X,Y;J,\theta_1,E)\subset\JJ(X,Y;J,\theta_1)$ be the set
of almost complex structures from Definition~\ref{def:JJ*}, which
is open and dense by Corollary~\ref{cor:open-dense}.
Let $\JJ_{\ell+1}$ be the set of coherent almost
complex structures from Section~\ref{sec:coh}. 

\begin{definition}\label{def:JJ*-ell}
We call an open subset $B^*\subset \JJ^*(X,Y;J,\theta_1,E)$ {\em
$\theta_2$-contractible} if it is contractible in
$\JJ^*(X,Y;J,\theta_1,E)$ to a point of
$\JJ^*(X,Y;J,\theta_2,E)$. Then we define 
\begin{align*}
   \JJ_{\ell+1}(X,&Y;J,\theta_1) \\
     := &\{K\in\JJ_{\ell+1}\mid
   K(\zeta)\in\JJ(X,Y;J,\theta_1)\text{ for all
   }\zeta\in\bar\MM_{\ell+1}, \cr
   &\ \ K(\zeta)|_Y \text{ independent of
   }\zeta\} \cr
   \JJ_{\ell+1}^*(X,&Y;J,\theta_1) \\
   := &\{K\in\JJ_{\ell+1}\mid
   K(\zeta)\in B^*\text{ for all
   }\zeta\in\bar\MM_{\ell+1},\cr
   &\ \  B^*\subset\JJ^*(X,Y;J,\theta_1,E) \text{
   $\theta_2$-contractible, }K(\zeta)|_Y \text{ independent of
   }\zeta\}.
\end{align*}
\end{definition}

Note that $\JJ_{\ell+1}(X,Y;J,\theta_1)$ is a Banach manifold, and
$\JJ_{\ell+1}^*(X,Y;J,\theta_1)\subset$\\ $\JJ_{\ell+1}(X,Y;J,\theta_1)$
is open and nonempty by Corollary~\ref{cor:open-dense} and
Lemma~\ref{lem:contraction}. 

\begin{proof}[Proof of Theorem~\ref{thm:moduli}]
Let $(J,Y)$ and $E,\ell$ be as above. 
Consider a homology class $A\in H_2(X;\Z)$ with $\om(A)=E$, so we
have 
$$
   \ell = D\om(A) = [Y]\cdot A.
$$
For $K\in\JJ_{\ell+1}^*(X,Y;J,\theta_1)$ and $k\geq 0$ define
the moduli space
\begin{align*}
   \MM_{k+\ell}(A,K;Y) := \{(\z,f)\mid & z_i\in S^2 \text{ distinct},
   f:S^2\to X,\ \pb_Kf=0,\ [f]=A, \cr
   & f(z_i)\in Y \text{ for } i=k+1,\dots,k+\ell\}/\Aut(S^2)
\end{align*}
of $K$-holomorphic spheres in the class $A$ passing at the last
$\ell$ points through $Y$. More generally, for an {\em
$\ell$-stable} $(k+\ell)$-labelled tree $T$ and
$K\in\JJ^*_{\ell+1}(X,Y;J,\theta_1)$ consider the moduli space
\begin{align*}
   \MM_{T}(A,K;Y) &:= \coprod_{\sum
   A_\alpha=A}\MM_{T}(\{A_\alpha\},K;Y) \cr
   &:= \{(\z,\f)\in\MM_T(A,K)\mid f_{\alpha_i}(z_i)\in Y \text{ for
   }i=k+1,\dots,k+\ell\}
\end{align*}
of stable maps passing at the last $\ell$ marked points through $Y$.
Observe for elements of $\MM_{T}(A,K;Y)$ every nonconstant component
passes through $V:=X\setminus Y$: This follows from property (a) in
Definition~\ref{def:JJ*} because  
$K(\zeta)|_Y$ is independent of $\zeta$ and belongs to
$\JJ^*(X,Y;J,\theta_1)$. 

For any $J_0\in \JJ^*(X,Y;J,\theta_1,E)$ we apply
Proposition~\ref{prop:trans} to the collection $\ZZ=\{Y\}$. Let  
$\JJ_{\ell+1}^\reg(V;\{Y\};J_0)\subset\JJ_{\ell+1}(V)$ be the Baire
set provided by Proposition~\ref{prop:trans} and set 
$$
   \JJ^\reg_{\ell+1}(X,Y;J,\theta_1) := \bigcup_{J_0\in
   \JJ^*(X,Y;J,\theta_1,E)}\JJ_{\ell+1}^\reg(V;\{Y\};J_0) \cap
   \JJ^*_{\ell+1}(X,Y;J,\theta_1). 
$$
Then for $K\in\JJ^\reg_{\ell+1}(X,Y;J,\theta_1)$ the
following holds: Let $A\in H_2(X;\Z)$ with $\ell=D\om(A)$, and $T$
be an $\ell$-stable $(k+\ell)$-labelled tree with weights $A_\alpha$,
$\sum A_\alpha=A$, such that every ghost tree contains at most one of
the last $\ell$ marked points. This condition implies that the reduced
index set $R$ satisfies $\{k+1,\dots,k+\ell\}\subset R$. So by
Proposition~\ref{prop:trans} (a), the moduli space $\MM_{T}^(A,K;Y)$
of stable $K$-holomorphic maps modelled over $T$ and passing through
$Y$ at the last $\ell$ marked points is a smooth manifold of dimension 
\begin{align*}
   \dim\MM_{T}(A,K;Y) &= 2n - 6 + 2c_1(A) + 2(k+\ell) - 2e(T) -
   2\ell\codim_\C(Y) \cr
   &= 2n - 6 + 2c_1(A) + 2k - 2e(T) \cr
   &= 2d -2e(T).
\end{align*}
This proves Theorem~\ref{thm:moduli}.
\end{proof}

In particular, applying Theorem~\ref{thm:moduli} to the tree $T$ with
one vertex (which has no ghost trees), it follows that the moduli
space $\MM_{k+\ell}(A,K;Y)$ is a smooth manifold of dimension
$$
   \dim_\R\MM_{k+\ell}(A,K;Y) = 2n-6+2k+2c_1(A) = 2d.
$$

{\bf Proof of Theorem~\ref{thm:pseudo}. }
We need to show that the evaluation map at the first $k$ points
$$
   \ev^k:\MM_{k+\ell}(A,K;Y)\mapsto X^k
$$
is a $2d$-dimensional pseudocycle. The following result is the key to
the proof.

\begin{prop}\label{prop:limit-stable}
Let $K\in\JJ^*_{\ell+1}(X,Y;J,\theta_1)$, $k\geq 0$, and $A\in
H_2(X;\Z)$ with $\ell=D\om(A)$. Suppose that a sequence
$[f^\nu]\in\MM_{k+\ell}(A,K;Y)$ converges in the sense of Gromov to
a stable map $[\z,\f]$. Then the nodal curve $\z$ is $\ell$-stable.
\end{prop}

\begin{proof}
Arguing by contradiction, suppose that $\z$ is not $\ell$-stable. Then
there exists a 
nonconstant component $f_\alpha$ of $\f$ whose domain $S_\alpha$
carries at most 2 special points after removing the first $k$ marked
points. According to condition (3) after
Theorem~\ref{thm:Gromov-comp}, $\f$ is
$K_{\pi_\ell(\z)}$-holomorphic, where $\pi_\ell$ is defined by
removing the first $k$ marked points and then stabilizing. As
$\pi_\ell$ maps $S_\alpha$ to a point, this means that
$\pb_{K_\alpha}f_\alpha=0$ for some $K_\alpha\in\JJ^*(X,Y;J,\theta_1)$
(not depending on points on $S_\alpha$).

By condition (a) in Definition~\ref{def:JJ*} (since $f_\alpha$ is
nonconstant), $f_\alpha$ is not
contained in $Y$, and by condition (b) it intersects $Y$ in at least 3
points in the domain. Thus $f_\alpha$ must intersect $Y$ at a point
that is neither a node nor one of the last $\ell$ marked points. By
Proposition~\ref{prop:int}, intersections persist under small
perturbations. So the definition of Gromov convergence implies that
for $\nu$ sufficiently large, $\f^\nu$ intersects $Y$ at a point
which is not one of the last $\ell$ marked points. On the other
hand, by Proposition~\ref{prop:int} each of the last $\ell$ marked
points contributes at least $1$ to the intersection number (since
$f_\nu:S^2\to X$ is a smooth $K$-holomorphic curve with distinct
marked points). Hence we conclude $[Y]\cdot[\f^\nu]>\ell=[Y]\cdot
A$, in contradiction to $[\f^\nu]=A$.
\end{proof}


\begin{prop}\label{prop:covering}
For all $K\in\JJ^\reg_{\ell+1}(X,Y;J,\theta_1)$ the following holds:
Let $k\geq 1$, $A\in H_2(X;\Z)$ with $\ell=D\om(A)$, and $T$ be an
$\ell$-stable $(k+\ell)$-labelled tree with more than one vertex. Then
the image under the evaluation map of the moduli space $\MM_T(A,K;Y)$
in $X^k$ is contained in the image of a manifold of dimension at most
$2d-2$. 
\end{prop}

\begin{proof}
If the tree $T$ satisfies condition (ii) in Definition~\ref{def:B}
(after removing the first $k$
marked points) the conclusion of the
proposition follows from Theorem~\ref{thm:moduli} because
$e(T)\geq 1$ if $T$ has more than one vertex.

It remains to analyze the case that there is a stable map
$(\z,\f)\in\MM_T(A,K;Y)$ not
satisfying condition (ii) in Definition~\ref{def:B} after removing
the first $k$
marked points. So there exists a ghost tree
$T'\subset T$ containing more than one of the last $\ell$ marked
points. On each such maximal ghost tree $T'$, remove all but one of the
last $\ell$ marked points and stabilize. Denote by $\bar T$ the resulting
$(k+\bar\ell)$-labelled tree, so $\bar T$ is $\bar\ell$-stable.  
According to Remark~\ref{rem:collapsing}, $K$ induces a coherent
almost complex structure $\bar K$ on $\pi^{-1}(\bar\MM_{\bar T})$
depending only on $z_0$ and the last $\bar\ell$ points. 
The resulting $\bar K$-holomorphic stable map $(\bar\z,\bar\f)$
modelled over $\bar T$  satisfies both conditions of
Definition~\ref{def:B} after removing the first $k$
marked points. Moreover, it has the
same image as  $(\z,\f)$ and represents the same homology class
$A$. Now we distinguish two cases.

Case 1: $\bar T$ has more than one vertex. 
Then by Theorem~\ref{thm:moduli}, $(\bar\z,\bar\f)$ belongs to a
moduli space of dimension $2d-e(\bar T)\leq 2d-2$.

Case 2: $\bar T$ has only one vertex. Then $(\bar\z,\bar\f)$ 
belongs to the moduli space
$\MM_{k+\bar\ell}(A,\bar K;Y)$ of smooth $\pi_{\bar\ell}^*\bar
K$-holomorphic maps $S^2\to X$ in the class $A$ mapping the last 
$\bar\ell$ marked points to
$Y$, where $\bar K\in\JJ^*_{\bar\ell+1}(X,Y;J,\theta_1)$ is the
almost complex structure induced by $K$. Now recall
that $(\z,\f)$ had at least two of the last $\ell$ marked points 
on the same ghost tree.  
So Lemma~\ref{lem:lim-int} shows that $\bar f$ has local intersection
number with $Y$ at least $2$ at the marked point at which this
ghost tree was attached, say at $\bar z_{k+1}$. Hence by
Proposition~\ref{prop:int} $\bar f$ actually belongs to the moduli space
$$
   \MM_{k+\bar\ell}(A,\bar K;\{Y,1\},\{Y,0\},\dots,\{Y,0\})
   \subset \MM_{k+\bar\ell}(A,\bar K;Y)
$$
of those maps that are tangent to $Y$ of order $1$ at $\bar
z_{k+1}$. 
By Corollary~\ref{cor:coh2} we may assume that 
$\bar K$ belongs to the Baire set $\JJ^\reg_{\bar
\ell+1}(V;\{Y,\dots,Y\})$ in Proposition~\ref{prop:jets1}. 
Then by Proposition~\ref{prop:jets1},
$\MM_{k+\bar\ell}(A,\bar K;\{Y,1\},\{Y,0\},\dots,\{Y,0\})$
is a smooth manifold of dimension $2d-2$. 

This shows that the image of each stable curve in a moduli space
$\MM_T(A,K;Y)$ with $|T|>1$ is contained in the image of some moduli
space of dimension at most $2d-2$. Since only finitely many moduli
spaces appear in this process, this proves
Proposition~\ref{prop:covering}.
\end{proof}

\begin{proof}[Proof of Theorem~\ref{thm:pseudo}]
Theorem~\ref{thm:pseudo} follows immediately from
Proposition~\ref{prop:limit-stable},
Proposition~\ref{prop:covering}, and the definition of a
pseudo-cycle.
\end{proof}

{\bf Extension to two hypersurfaces. }
Now we extend the preceding discussion to the case of two
hypersurfaces. 

\begin{definition}\label{def:D-quadruple}
A {\em Donaldson quadruple} is a quadruple $(J_0,Y_0,J_1,Y_1)$
consisting of two Donaldson pairs $(J_i,Y_i)$ satisfying 
$$
   \angle_m(Y_0,Y_1)\geq\eps,\qquad \|J_0-J_1\|_0<\theta_3. 
$$
\end{definition}

The following is the analogue of Theorem~\ref{thm:moduli} in this
case. 

\begin{thm}\label{thm:moduli2}
Let $(J_0,Y_0,J_1,Y_1)$ be a Donaldson quadruple on $(X,\om)$ of
bidegree $(D_0,D_1)$. Then there exist nonempty 
sets $\JJ^\reg_{\ell+1}(X,Y_0\cup Y_1;J_0,\theta_1)\subset$ \\
$\JJ_{\ell+1}(X,Y_0\cup Y_1;J_0,\theta_1)$, 
$\ell\geq 3$, such that for any subset $I\subset\{1,\dots,\ell\}$ with 
$|I|\geq \ell\min(D_0,D_1)/(D_0+D_1)\geq 3$ and
$K\in\JJ_{|I|+1}^\reg(X,Y_0\cup Y_1;J_0,\theta_1)$ the following 
holds: Let $A\in H_2(X;\Z)$ be a homology class with
$\ell_i=D_i\om(A)$ and $\ell=\ell_0+\ell_1$. For $k\geq 0$ let $T$ be
an $(k+\ell)$-labelled tree with $e(T)$ edges such that every
ghost tree contains at most one of the last $\ell$ marked points. Then
the moduli space 
$\MM_T(A,K;Y_0\cup Y_1)$ of stable $K$-holomorphic spheres in the class $A$,
modelled over $T$ and mapping the middle $\ell_0$ marked points to
$Y_0$ and the last $\ell_1$ marked points to $Y_1$, is a
smooth manifold of real dimension
$$
   \dim_\R\MM_T(A,K;Y_0\cup Y_1) = 2\Bigl(n-3+k+c_1(A)-e(T)\Bigr).
$$
\end{thm}

\begin{proof}
To each integer $\ell\geq 3$ we associate the energy 
$$
   E_\ell:=\ell/\min(D_0,D_1).
$$ 
For this energy, let $\JJ^*(X,Y_0\cup
Y_1;J_0,\theta_1,E_\ell)\subset\JJ(X,Y_0\cup Y_1;J_0,\theta_1)$ be 
the set of almost complex structures from Definition~\ref{def:JJ*2},
which is open and dense by Corollary~\ref{cor:open-dense2}.
Let $\JJ_{\ell+1}$ be the set of coherent almost
complex structures from Section~\ref{sec:coh}.

\begin{definition}\label{def:JJ*-ell2}
We call an open subset $B^*\subset \JJ^*(X,Y_0\cup
Y_1;J_0,\theta_1,E_\ell)$ {\em $\theta_2$-contractible} if it is
contractible in $\JJ^*(X,Y_0\cup Y_1;J_0,\theta_1,E_\ell)$ to a point
in \\ $\JJ^*(X,Y_0\cup Y_1;J_0,\theta_2,E_\ell)$. Then we
define 
\begin{align*}
   \JJ_{\ell+1}(X,&Y_0\cup Y_1;J_0,\theta_1)\\
   := &\{K\in\JJ_{\ell+1}\mid
   K(\zeta)\in\JJ(X,Y_0\cup Y_1;J_0,\theta_1)\text{ for all
   }\zeta\in\bar\MM_{\ell+1}, \cr
   &\ \ K(\zeta)|_{Y_0\cup Y_1} \  \ \text{ independent of
   }\zeta\} \cr
   \JJ_{\ell+1}^*(X,&Y_0\cup Y_1;J_0,\theta_1)\\
   := &\{K\in\JJ_{\ell+1}\mid
   K(\zeta)\in B^*\text{ for all
   }\zeta\in\bar\MM_{\ell+1},\cr
   &\ \  B^*\subset \JJ^*(X,Y_0\cup Y_1;J_0,\theta_1,E_\ell) \text{
   $\theta_2$-contractible, }\\
   &\ \ K(\zeta)|_{Y_0\cup Y_1}
   \  \ \text{ independent of
   }\zeta\}.
\end{align*}
\end{definition}
Note that $\JJ_{\ell+1}(X,Y_0\cup Y_1;J_0,\theta_1)$ is a Banach
manifold, and
$\JJ_{\ell+1}^*(X,Y_0\cup Y_1;J_0,\theta_1)\subset
\JJ_{\ell+1}(X,Y_0\cup Y_1;J_0,\theta_1)$ 
is open and by Corollary~\ref{cor:open-dense2} and
Lemma~\ref{lem:contraction}. 

Now consider a homology class $A\in H_2(X;\Z)$ with $0<\om(A)$
and set
$$
   \ell_i := D_i\om(A) = [Y_i]\cdot A,\qquad \ell:=\ell_0+\ell_1. 
$$
For $k\geq 0$ we denote points in the Deligne-Mumford space
$\bar\MM_{k+\ell_0+\ell_1+1}$ by
$$
   \bar z = (z_0,\dots,z_k,z_{k+1},\dots,z_{k+\ell_0},z_{k+\ell_0+1},
   \dots,z_{k+\ell_0+\ell_1}), 
$$
so the first components after $z_0$ correspond to $k$, the middle ones
to $\ell_0$, and the last ones to $\ell_1$. Fix a subset
$I\subset\{1,\dots,\ell\}$ with 
$$
   |I|\geq\frac{\min(D_0,D_1)}{D_0+D_1}\ell = \min(D_0,D_1)\om(A) \geq
   3,
$$
hence $\om(A)\leq E_{|I|}$. Denote by  
$$
   \pi_I:\bar\MM_{k+\ell_0+\ell_1+1}\to\bar\MM_{|I|+1}
$$
the obvious projection (forgetting the marked points outside the
shifted set $I+k$ and stabilizing).  
Recall that a $k+\ell_0+\ell_1$-labelled tree is called $I$-stable if
it is still stable after removing the marked points not corresponding
to $I$. 

For $K\in\JJ_{|I|+1}^*(X,Y_0\cup Y_1;J_0,\theta_1)$ (as defined above
for the energy $E_{|I|}$) define the moduli space
\begin{align*}
   \MM_{k+\ell_0+\ell_1}(A,K; &Y_0\cup Y_1) \\
   := \{(\z,f)\mid & z_i\in S^2 \text{ distinct},
   f:S^2\to X,\ \pb_{\pi_I^*K}f=0,\ [f]=A, \cr
   & f(z_i)\in Y_0 \text{ for }k+1\leq i\leq k+\ell_0, \cr
   & f(z_i)\in Y_1 \text{ for }k+\ell_0+1\leq i\leq k+\ell_0+\ell_1\}/\Aut(S^2)
\end{align*}
of $K$-holomorphic spheres in the class $A$ passing at the middle 
$\ell_0$ points through $Y_0$ and at the last $\ell_1$ points through
$Y_1$. More generally, for an 
$I$-stable $(k+\ell_0+\ell_1)$-labelled tree $T$ and
$K\in\JJ^*_{|I|+1}(X,Y_0\cup Y_1;J_0,\theta_1)$ consider the moduli space
\begin{align*}
   \MM_T(A,K;Y_0\cup Y_1) &:= \coprod_{\sum
   A_\alpha=A}\MM_T(\{A_\alpha\},K;Y_0\cup Y_1) \cr
   &:= \{(\z,\f)\in\MM_T(A,K)\mid 
   f_{\alpha_i}(z_i)\in Y_0 \text{ for }k+1\leq i\leq k+\ell_0, \cr
   & \ \ \ \ \ f_{\alpha_i}(z_i)\in Y_1 \text{ for }k+\ell_0+1\leq i\leq
   k+\ell_0+\ell_1 \} 
\end{align*}
of stable maps passing at the middle $\ell_0$ marked points through
$Y_0$ and at the last $\ell_1$ points through $Y_1$. 
Observe that elements of $\MM_T(A,K;Y_0\cup Y_1)$
satisfy Condition (i) in Definition~\ref{def:B} with
$V:=X\setminus(Y_0\cup Y_1)$: This
follows from property (a) in Definition~\ref{def:JJ*2} because
$K(\zeta)|_{Y_0\cup Y_1}$ is independent of $\zeta$ and belongs to
$\JJ^*(X,Y_0\cup Y_1;J_0,\theta_1,E_{|I|})$, and $\om(A)\leq
E_{|I|}$).

For any $J_0\in \JJ^*_{|I|+1}(X,Y_0\cup Y_1;J_0,\theta_1,E)$,  
we apply Corollary~\ref{cor:surj-tree2} with $Z=(Y_0\cup Y_1)^\ell$,
$V:=X\setminus (Y_0\cup Y_1)$, $I=\{k+1,\dots,k+\ell\}$, 
and with $\JJ_{\ell+1}(V)$ (the set of coherent almost complex
structures that agree with $J_0$ outside $V$) replaced by its open
subset $\JJ_{|I|+1}^*(X,Y_0\cup
Y_1;J_0,\theta_1)\cap \JJ_{\ell+1}(V)$ according
to Remark~\ref{rem:add}.   
Then Corollary~\ref{cor:surj-tree2} yields a Baire set
$\JJ^\reg_{|I|+1}(X,Y_0\cup Y_1;J_0,\theta_1)\subset
\JJ^*_{|I|+1}(X,Y_0\cup Y_1;J_0,\theta_1)$ such that the following holds: Let
$K\in\JJ^\reg_{|I|+1}(X,Y_0\cup Y_1;J_0,\theta_1)$, 
$A\in H_2(X;\Z)$ with $\ell=D\om(A)$, and $T$
be an $I$-stable $(k+\ell_0+\ell_1)$-labelled tree satisfying condition
(ii) in Definition~\ref{def:B} (after removing the first $k$
marked points). Then the moduli space $\MM_T^{\ell*}(A,K;Y_0\cup Y_1)$ is a
smooth manifold of dimension
\begin{align*}
   \dim\MM_T^{\ell*}(A,K; &Y_0\cup Y_1)\\
    &= 2n - 6 + 2c_1(A) +
   2(k+\ell_0+\ell_1) \!-\! 2e(T) \!-\! \codim_\R(
   Y_0^{\ell_0}\!\times\! Y_1^{\ell_1}) \cr
   &= 2n - 6 + 2c_1(A) +2k - 2e(T) \cr
   &= 2d -2e(T).
\end{align*}
This proves Theorem~\ref{thm:moduli2}.
\end{proof}

In particular, applying Theorem~\ref{thm:moduli2} to the tree $T$ with
one vertex (which automatically satisfies condition (ii) in
Definition~\ref{def:B}), it follows that the moduli space
$\MM_{k+\ell_0+\ell_1}(A,K;Y_0\cup Y_1)$ is a smooth manifold of dimension
$$
   \dim_\R\MM_{k+\ell_0+\ell_1}(A,K;Y_0\cup Y_1) = 2n-6+2k+2c_1(A) = 2d.
$$
For Theorem~\ref{thm:pseudo} (extended to the case of two hypersurfaces),
we need to show that the evaluation map
at the first $k$ points
$$
   \ev^k:\MM_{k+\ell_0+\ell_1}(A,K;Y_0\cup Y_1)\mapsto X^k
$$
is a $2d$-dimensional pseudocycle. The key to the proof is the following
result which we will apply in the cases $I=\{k+1,\dots,k+\ell_0\}$ and
$I=\{k+1,\dots,k+\ell_0+\ell_1\}$. 

\begin{prop}\label{prop:limit-stable2}
For $\{k+1,\dots,k+\ell_0\}\subset
I\subset\{k+1,\dots,k+\ell_0+\ell_1\}$  
let $K\in\JJ^*_{|I|+1}(X,Y_0\cup Y_1;J_0,\theta_1)$, $k\geq 0$, $0<\om(A)$
and $\ell_i=D_i\om(A)$. Suppose that a sequence
$[f^\nu]\in\MM_{k+\ell_0+\ell_1}(A,K;Y_0\cup Y_1)$ converges in the sense of Gromov to
a stable map $[\z,\f]$. Then the nodal curve $\z$ is $I$-stable.
\end{prop}

\begin{proof}
Arguing by contradiction, suppose that $\z$ is not $I$-stable. Then there exists a
nonconstant component $f_\alpha$ of $\f$ whose domain $S_\alpha$
carries at most 2 special points after removing the marked points not
belonging to $I$. According to condition (3) after
Theorem~\ref{thm:Gromov-comp}, $\f$ is
$K_{\pi_I(\z)}$-holomorphic, where $\pi_I$ is defined by
removing the marked points outside $I$ and then stabilizing. As
$\pi_I$ maps $S_\alpha$ to a point, this means that
$\pb_{K_\alpha}f_\alpha=0$ for some $K_\alpha\in\JJ^*(X,Y_0\cup
Y_1;J_0,\theta_1,E_{|I|})$ (not depending on points on $S_\alpha$),
where $E_{|I|}=|I|/\min(D_0,D_1)$.

By condition (a) in Definition~\ref{def:JJ*2} (since $f_\alpha$ is
nonconstant), $f_\alpha$ is not
contained in $Y_0\cup Y_1$, and by condition (b) it intersects $Y_0$ in at least 3
points in the domain. Thus $f_\alpha$ must intersect $Y_0$ at a point
that is neither a node nor one of the middle $\ell_0$ marked points. By
Proposition~\ref{prop:int}, intersections persist under small
perturbations. So the definition of Gromov convergence implies that
for $\nu$ sufficiently large, $\f^\nu$ intersects $Y_0$ at a point
which is not one of the middle $\ell_0$ marked points. On the other
hand, by Proposition~\ref{prop:int} each of the middle $\ell_0$ marked
points contributes at least $1$ to the intersection number (since
$f_\nu:S^2\to X$ is a smooth $K$-holomorphic curve with distinct
marked points). Hence we conclude $[Y_0]\cdot[\f^\nu]>\ell_0=[Y_0]\cdot
A$, in contradiction to $[\f^\nu]=A$.
\end{proof}

\begin{prop}\label{prop:covering2}
For $\{k+1,\dots,k+\ell_0\}\subset
I\subset\{k+1,\dots,k+\ell_0+\ell_1\}$ there exists a Baire set
$\JJ^\reg_{|I|+1}(X,Y_0\cup Y_1;J_0,\theta_1)\subset
\JJ^*_{|I|+1}(X,Y_0\cup Y_1;J_0,\theta_1)$ such that 
for all $K\in\JJ^\reg_{|I|+1}(X,Y_0\cup Y_1;J_0,\theta_1)$ the following holds: Let $k\geq
1$, $0<\om(A)$, $\ell_i=D_i\om(A)$, and $T$ be an $I$-stable
$(k+\ell_0+\ell_1)$-labelled tree with more than one vertex. Then the image
under the evaluation map of the moduli space $\MM_T(A,K;Y_0\cup Y_1)$ in
$X^k$ is contained in the image of a manifold of dimension at most $2d-2$.
\end{prop}

\begin{proof}
Set $\ell:=\ell_0+\ell_1$. 
If the tree $T$ satisfies condition (ii) in Definition~\ref{def:B}
(after removing the last $\ell$ marked points) the conclusion of the
proposition follows from Theorem~\ref{thm:moduli2} because
$e(T)\geq 1$ if $T$ has more than one vertex.

It remains to analyze the case that there is a stable map
$(\z,\f)\in\MM_T(A,K;Y_0\cup Y_1)$ not
satisfying condition (ii) in Definition~\ref{def:B} after removing
the last $\ell$ marked points. So there exists a ghost tree
$T'\subset T$ containing more than one of the last $\ell$ marked
points. On each such maximal ghost tree $T'$, remove all but one of the
last $\ell$ marked points and stabilize. Denote by $\bar T$ the resulting
$(k+\bar\ell_0+\bar\ell_1)$-labelled tree and by $\bar I$ the set of
points on $\bar T$ remaining from $I$, so $\bar T$ is $\bar I$-stable.  
According to Remark~\ref{rem:collapsing}, $K$ induces a coherent
almost complex structure $\bar K$ on $\pi^{-1}(\bar\MM_{\bar T})$
depending only on $z_0$ and points in $\bar I$.  
The resulting $\bar K$-holomorphic stable map $(\bar\z,\bar\f)$
modelled over $\bar T$  satisfies both conditions of
Definition~\ref{def:B} after removing the last $\bar\ell$ marked
points, where $\bar\ell:=\bar\ell_0+\bar\ell_1$. Moreover, it has the
same image as  $(\z,\f)$ and represents the same homology class
$A$. Now we distinguish two cases.

Case 1: $\bar T$ has more than one vertex. 
Then by Theorem~\ref{thm:moduli2}, $(\bar\z,\bar\f)$ belongs to a
moduli space of dimension $2d-e(\bar T)\leq 2d-2$.

Case 2: $\bar T$ has only one vertex. Then $(\bar\z,\bar\f)$ 
belongs to the moduli space
$\MM_{k+\bar\ell_0+\bar\ell_1}(A,\bar K;Y_0\cup Y_1)$ of smooth 
$\pi_{\bar I}^*\bar K$-holomorphic maps $S^2\to X$ in the class $A$
mapping the middle 
$\bar\ell_0$ marked points to $Y_0$ and the last $\bar\ell_1$ ones to
$Y_1$, where $\bar K\in\JJ^*_{|I|+1}(X,Y_0\cup Y_1;J_0,\theta_1)$ is the
almost complex structure induced by $K$. Now recall
that $(\z,\f)$ had at least two of the last $\ell$ marked points 
on the same ghost tree. There are three subcases. 

Case 2a: $(\z,\f)$ had two of the middle $\ell_0$ marked
points on the same ghost tree. 
Then Lemma~\ref{lem:lim-int} shows that $\bar f$ has local intersection
number with $Y_0$ at least $2$ at the marked point at which this
ghost tree was attached, say at $\bar z_{k+1}$. Hence by
Proposition~\ref{prop:int} $\bar f$ actually belongs to the moduli space
$$
   \MM_{k+\bar\ell_0+\bar\ell_1}(A,\bar K;\{Y_0,1\},\{Y_0,0\},\dots,\{Y_1,0\})
   \subset \MM_{k+\bar\ell_0+\bar\ell_1}(A,\bar K;Y_0\cup Y_1)
$$
of those maps that are tangent to $Y_0$ of order $1$ at $\bar
z_{k+1}$. 
By Corollary~\ref{cor:coh2} we may assume that 
$\bar K$ belongs to the Baire set $\JJ^\reg_{|\bar
  I|+1}(V;\{Y_0,\dots,Y_1\})$ in Proposition~\ref{prop:jets1} with
$V:=X\setminus(Y_0\cup Y_1)$. 
Then by Proposition~\ref{prop:jets1},\\
$\MM_{k+\bar\ell_0+\bar\ell_1}(A,\bar K;\{Y_0,1\},\{Y_0,0\},\dots,\{Y_1,0\})$
is a smooth manifold of dimension $2d-2$.

Case 2b: $(\z,\f)$ had two of the last $\ell_1$ marked
points on the same ghost tree. This is analogous to Case 2a. 

Case 2c: $(\z,\f)$ had one of the middle $\ell_0$ and one of the last
$\ell_1$ marked points on the same ghost tree. 
Then $\bar f$ maps the marked point at which this ghost tree was
attached, say $\bar z_{k+1}$, to $Y_0\cap Y_1$. Hence $\bar f$
actually belongs to the moduli space 
$$
   \MM_{k+\bar\ell_0+\bar\ell_1}(A,K;Y_0\cap Y_1,Y_0,\dots,Y_1)
   \subset \MM_{k+\bar\ell_0+\bar\ell_1}(A,K;Y_0\cup Y_1)
$$
of those maps that map $\bar z_{k+1}$ to $Y_0\cap Y_1$. 
Again by Corollary~\ref{cor:coh2} we may assume that $\bar K$ belongs to
the Baire set $\JJ^\reg_{|I|+1}(V;\{Y_0\cap Y_1,Y_0,\dots,Y_1\})$ in
Proposition~\ref{prop:jets1}. 
Then by Proposition~\ref{prop:jets1},
$\MM_{k+\bar\ell_0+\bar\ell_1}(A,K;Y_0\cap Y_1,Y_0,\dots,Y_1)$
is a smooth manifold of dimension $2d-2$.

This shows that the image of each stable curve in a moduli space
$\MM_T(A,K;Y_0\cup Y_1)$ with $|T|>1$ is contained in the image of some moduli
space of dimension at most $2d-2$. Since only finitely many moduli
spaces appear in this process, this proves
Proposition~\ref{prop:covering2}.
\end{proof}

Proposition~\ref{prop:limit-stable2},
Proposition~\ref{prop:covering2}, and the definition of a
pseudo-cycle immediately imply the following version of
Theorem~\ref{thm:pseudo} for two hypersurfaces. 

\begin{thm}\label{thm:pseudo2}
In the notation of Theorem~\ref{thm:moduli2}, for every $k\geq 1$ the
evaluation map at the first $k$ points
$$
   \ev^k:\MM_{k+\ell_0+\ell_1}(A,K;Y_0\cup Y_1)\mapsto X^k
$$
represents a pseudocycle of dimension
$$
   2d := 2\Bigl(n-3+k+c_1(A)\Bigr).
$$
\end{thm}

\section{Independence of the auxiliary data}\label{sec:indep}

In this section we prove that different choices of the auxiliary
data $(Y,J;K)$, where $(J,Y)$ is a Donaldson pair and
$K\in\JJ^\reg_{\ell+1}(X,Y;J,\theta_1)$, give rise to
rational pseudo-cycles which are rationally cobordant.

Throughout this section we fix the following data:
\begin{itemize}
\item a symplectic manifold $(X,\om)$ with integral class $[\om]$;
\item 
a homology class $A\in H_2(X;\Z)$ with $0<\om(A)$;
\item an integer $k\geq 1$;
\item a constant $\eps>0$ and constants
  $0<\theta_3<\theta_2<\theta_1<\theta_0<1$ as in
  Lemma~\ref{lem:contraction} and Lemma~\ref{lem:angle3} (b); 
\item a closed 1-form $\alpha$ representing the first Chern class of
  $X$. 
\end{itemize}

As in the previous section, by a {\em Donaldson pair} we understand a
pair 
$(J,Y)$, where 
\begin{itemize}
\item $J$ is an $\om$-compatible almost complex structure on $X$;
\item $Y\subset X$ is a hypersurface whose K\"ahler angle with respect
to $(\om,J)$ satisfies $\theta(Y;\om,J)<\theta_3$ and such that
$[Y]=D[\om]$ with $D\geq D^*(X,\om,J)$, where $D^*(X,\om,J)$ is the
constant defined 
in Section~\ref{sec:proof}. 
\end{itemize}

We always set $\ell:=D\om(A)$. Theorem~\ref{thm:pseudo} assigns to
each Donaldson pair $(J,Y)$ and domain-dependent perturbation
$K\in\JJ^\reg_{\ell+1}(X,Y;J,\theta_1)$ in the Baire set of
Theorem~\ref{thm:moduli} a pseudo-cycle 
$$
   \ev^k(J,Y;K):\MM_{k+\ell}(A,K;Y)\to X^k. 
$$
We first prove independence of the perturbation $K$. 

\begin{prop}\label{prop:uniqueness1}
For a Donaldson pair $(J,Y)$ and
$K_0,K_1\in\JJ^\reg_{\ell+1}(X,Y;J,\theta_1)$, the pseudo-cycles
$$
   \ev^k_i(J,Y;K_i):\MM_{k+\ell}(A,K_i;Y)\to X^k, \qquad i=0,1
$$
are cobordant.
\end{prop}

\begin{proof}
Note that we can connect $K_0,K_1$ by a smooth
path $K_t$ in $\JJ^*_{\ell+1}(X,Y;J,\theta_1)$. Indeed, by the
$\theta_2$-contractibility condition in Definition~\ref{def:JJ*-ell}
we can connect the $K_i$ in $\JJ^*_{\ell+1}(X,Y;J,\theta_1)$ to
domain-independent almost complex structures $K_i^*\in
\JJ^*(X,Y;J,\theta_2,E)$, and by Corollary~\ref{cor:open-dense} we can
connect $K_0^*$ and $K_1^*$ by a path in $\JJ^*(X,Y;J,\theta_1,E)$. 
Now the proof works exactly as in~\cite{MS}, repeating the arguments
in Section~\ref{sec:proof} for a generic path $K_t$ in
$\JJ^*_{\ell+1}(X,Y;J,\theta_1)$.  
\end{proof}

As in the previous section, we denote by {\em Donaldson quadruple} a quadruple
$(J_0,Y_0,J_1,Y_1)$ consisting of two Donaldson pairs $(J_i,Y_i)$ such
that $Y_0,Y_1$ are $\eps$-transverse and $\|J_0-J_1\|_0<\theta_3$. 
Theorem~\ref{thm:pseudo2} associates a pseudo-cycle 
$$
   \ev^k(J_0,Y_0,J_1,Y_1;K):\MM_{k+\ell_0+\ell_1}(A,K;Y_0\cup Y_1)\to
   X^k 
$$
to each Donaldson quadruple $(J_0,Y_0,J_1,Y_1)$ and domain-dependent
perturbation $K\in\JJ^\reg_{\ell_0+\ell_1+1}(X,Y_0\cup
Y_1;J_0,\theta_1)$ in the Baire set of Theorem~\ref{thm:moduli2} with
$\ell_i:=D_i\om(A)$. Again, this pseudo-cycle does not depend on the
perturbation $K$:

\begin{prop}\label{prop:uniqueness2}
For a Donaldson quadruple $(J_0,Y_0,J_1,Y_1)$ and
$K_0,K_1\in\JJ^\reg_{\ell_0+\ell_1+1}(X,Y_0\cup Y_1;J_0,\theta_1)$,
the pseudo-cycles 
$$
   \ev^k_i(J_0,Y_0,J_1,Y_1;K_i):\MM_{k+\ell_0+\ell_1}(A,K_i;Y_0\cup
   Y_1)\to X^k, \qquad i=0,1 
$$
are cobordant.
\end{prop}

\begin{proof}
As in the previous proof, by Definition~\ref{def:JJ*-ell2} and
Corollary~\ref{cor:open-dense2} we can connect $K_0,K_1$
by a smooth path $K_t$ in $\JJ^*_{\ell_0+\ell_1+1}(X,Y_0\cup
Y_1;J_0,\theta_0)$ and the proof works as in~\cite{MS}, repeating the
arguments in Section~\ref{sec:proof} for a generic path $K_t$. 
\end{proof}

In view of Propositions~\ref{prop:uniqueness1}
and~\ref{prop:uniqueness2}, we can drop the perturbation $K$ from our
notation and denote by $\ev^k(J,Y)$ resp.~$\ev^k(J_0,Y_0,J_1,Y_1)$ the
induced rational cobordism classes of pseudo-cycles. 

The following result relates the rational pseudo-cycle associated to a
Donaldson quadruple $(J_0,Y_0,J_1,Y_1)$ to those associated to the
Donaldson pairs $(J_i,Y_i)$.

\begin{prop}\label{prop:uniqueness3}
For each Donaldson quadruple $(J_0,Y_0,J_1,Y_1)$ there exist rational
cobordisms 
$$
   \frac{1}{\ell_0!}\ev^k(J_0,Y_0)
   \sim \frac{1}{\ell_0!\ell_1!}\ev^k(J_0,Y_0,J_1,Y_1)
   \sim \frac{1}{\ell_1!}\ev^k(J_1,Y_1).  
$$
\end{prop}

\begin{proof}
We prove the first equivalence, the second one being analogous. 
Pick a perturbation 
$K_0\in\JJ^\reg_{\ell_0+1}(X,Y_0\cup Y_1;J_0,\theta_1)$ 
in the Baire set of Theorem~\ref{thm:moduli2}
depending only on the middle $\ell_0$ marked points. (This corresponds
to taking the index set $I=\{1,\dots,\ell_0\}$ in
Theorem~\ref{thm:moduli2}). Then in particular we have 
$K_0\in\JJ^\reg_{\ell_0+1}(X,Y_0;J_0,\theta_1)$, and 
its pullback
under the projection $\pi_{\ell_0}$ forgetting the last $\ell_1$
points satisfies
$\pi_{\ell_0}^*K_0\in\JJ^\reg_{\ell_0+\ell_1+1}(X,Y_0\cup
Y_1;J_0,\theta_1)$.  
Thus $\pi_{\ell_0}$ induces a commutative diagram of pseudo-cycles
\begin{equation*}
\begin{CD}
   \MM_{k+\ell_0+\ell_1}(A,\pi_{\ell_0}^*K_0;Y_0\cup Y_1) @>\ev^k>>
   X^k \\
   @VV\pi_{\ell_0}V @VV\id V \\
   \MM_{k+\ell_0}(A,K_0;Y_0) @>\ev^k>> X^k.
\end{CD}
\end{equation*}
Now note that the moduli space
$\MM_{k+\ell_0+\ell_1}(A,\pi_{\ell_0}^*K_0;Y_0\cup Y_1)$ is obtained
from $\MM_{k+\ell_0}(A,K_0;Y_0)$ by simply adding the $\ell_1$
intersection points with $Y_1$ as additional marked points. (All
these intersection points are distinct in the top stratum because
any coincidence leads to a tangency to $Y_1$ and thus to a
codimension 2 stratum). Hence the map
$$
   \pi_{\ell_0}: \MM_{k+\ell_0+\ell_1}(A,\pi_{\ell_0}^*K_0;Y_0\cup
   Y_1) \to \MM_{k+\ell_0}(A,K_0;Y_0)
$$
is a covering map of degree $\ell_1!$ and preceding commutative
diagram yields the equality as currents 
$$
    \ev^k(J_0,Y_0;K_0) \sim
    \frac{1}{\ell_1!}\ev^k(J_0,Y_0,J_1,Y_1;\pi_{\ell_0}^*K_0). 
$$
In view of Propositions~\ref{prop:uniqueness1}
and~\ref{prop:uniqueness2}, this proves the first equivalence and
hence Proposition~\ref{prop:uniqueness3}. 
\end{proof}

The final ingredient is invariance under homotopies of Donaldson
pairs. 

\begin{prop}\label{prop:uniqueness4}
Let $(J_t,Y_t)_{t\in[0,1]}$ be a smooth family of Donaldson
pairs. Suppose that $Y_t=\phi_t(Y_0)$ for a smooth family of
symplectomorphisms $\phi_t$ with $\phi_0=\id$. Then we have a rational
cobordism of pseudo-cycles
$$
   \ev^k(J_0,Y_0)\sim \ev^k(J_1,Y_1).
$$
\end{prop}

\begin{proof}
After cutting the interval $[0,1]$ into smaller ones we may assume
without loss of generality that $\|J_0-J_1\|$ is arbitrarily small
and $\phi:=\phi_1$ is arbitrarily $C^1$-close to the identity. Then
$\|\phi_*J_0-J_1\|$ will also be arbitrarily small. 

Let $D$ be the degree of the $Y_t$ and $\ell:=D\om(A)$. 
Pick $K_0\in\JJ^\reg_{\ell+1}(X,Y_0;J_0,\theta_1/2)$, so 
$\phi_*K_0\in\JJ^\reg_{\ell+1}(X,Y_1;\phi_*J_0,\theta_1/2)$. 
For $\|\phi_*J_0-J_1\|$ sufficiently small this implies 
$\phi_*K_0\in\JJ^\reg_{\ell+1}(X,Y_1;J_1,\theta_1)$.
Composition with $\phi$ induces a commutative diagram of pseudo-cycles 
\begin{equation*}
\begin{CD}
   \MM_{k+\ell}(A,K_0;Y_0) @>\ev^k>> X^k \\
   @V\cong V\phi V @VV\phi^k V \\
   \MM_{k+\ell}(A,\phi_*K_0;Y_1) @>\ev^k>> X^k.
\end{CD}
\end{equation*}
Since $\phi$ is smoothly homotopic to the identity, this implies the
following rational cobordisms of pseudo-cycles:
$$
   \ev^k(J_0,Y_0;K_0) \sim \phi^k\circ\ev^k(J_0,Y_0;K_0) \sim 
   \ev^k(J_1,Y_1;\phi_*K_0). 
$$
In view of Proposition~\ref{prop:uniqueness1}, 
this proves Proposition~\ref{prop:uniqueness4}. 
\end{proof}

\begin{remark}
Proposition~\ref{prop:uniqueness4} can also be proved, without using
the symplectic isotopy $\phi_t$, by considering the moduli space
$\cup_{t\in[0,1]}\MM_{k+\ell}(A,K_t;Y_t)$ for a generic family
$(K_t)_{t\in[0,1]}$. 
\end{remark}

\begin{proof}[Proof of Theorem~\ref{thm:indep}]
Let $(J_i,Y_i)$, $i=0,1$ be Donaldson pairs of degrees $D_i$ and set
$\ell_i:=D_i\om(A)$.  
Since $\theta(Y_i;\om,J_i)<\theta_3$, Lemma~\ref{lem:angle2} (c) provides
compatible almost complex structures $\bar J_i$ preserving $TY_i$
with $\|\bar J_i-J_i\|<\theta_3$.  
By the (Stabilization) statement in
Theorem~\ref{thm:Donaldson}, for each sufficiently large integer
$\bar D$ we find hypersurfaces $\bar Y_i$ of degree
$\bar D$ such that $(J_i,Y_i,\bar J_i,\bar Y_i)$ are Donaldson
quadruples for $i=0,1$. Set $\bar\ell:=\bar D\om(A)$. 
Connect $\bar J_0$ and $\bar J_1$ by a smooth path of compatible
almost complex structures $\bar J_t$. 
By the (Uniqueness) statement in Theorem~\ref{thm:Donaldson}, for
$\bar D$ sufficiently large there exist smooth paths of hypersurfaces
$\bar Y_t$ and symplectomorphisms $\phi_t$ with $\phi_0=\id$ such that
the $(\bar J_t,\bar Y_t)$ are Donaldson pairs and $\bar
Y_t=\phi_t(\bar Y_0)$ for all $t\in[0,1]$. Now
Propositions~\ref{prop:uniqueness3} and~\ref{prop:uniqueness4} yield
the following rational cobordisms of rational pseudo-cycles:
$$
   \frac{1}{\ell_0!}\ev^k(J_0,Y_0) 
   \sim \frac{1}{\bar\ell!}\ev^k(\bar J_0,\bar Y_0) 
   \sim \frac{1}{\bar\ell!}\ev^k(\bar J_1,\bar Y_1) 
   \sim \frac{1}{\ell_1!}\ev^k(J_1,Y_1).  
$$
This concludes the proof of Theorem~\ref{thm:indep}. 
\end{proof}


\end{document}